\documentclass[11pt]{article} 

\usepackage{type1cm}        
\usepackage[export]{adjustbox}
\usepackage{makeidx,}
\usepackage{graphicx}        
\usepackage{multicol,amsthm}        
\usepackage[bottom]{footmisc}

\usepackage{newtxtext}       %
\usepackage{newtxmath,color}       

\usepackage{subfigure}
\newcommand{\dpar}[2]{\dfrac{\partial #1}{\partial #2}}

 \newcommand{\R}{\mathbb R}
 \newcommand{\Z}{\mathbb Z}

\newcommand{\PP}{\mathbb P}

\newcommand{\bbf}{{\mathbf {f}}}
\newcommand{\bbg}{{\mathbf {g}}}
\newcommand{\bbs}{{\mathbf {s}}}
\newcommand{\bn}{{\mathbf {n}}}

\newcommand{\bu}{\mathbf{u}}
\newcommand{\bv}{\mathbf{v}}
\newcommand{\bw}{\mathbf{w}}

\newcommand{\bbB}{\mathbf{B}}

\newcommand{\bF}{\mathbf{f}}

\newcommand{\bbu}{\mathbf{u}}
\newcommand{\hbbg}{\hat{\mathbf{g}}}
\newcommand{\hbbf}{\hat{\mathbf{f}}}
\newcommand{\hf}{\hat{{f}}}
\newcommand{\hg}{\hat{{g}}}

\newcommand{\bx}{\mathbf{x}}
\newcommand{\by}{\mathbf{y}}

\newcommand{\normal}{\mathbf{N}}

\newcommand{\bbv}{\pmb{v}}
\newcommand{\bX}{\pmb{X}}
\newcommand{\bbw}{\pmb{w}}

\newcommand{\jump}[1]{[\![ #1]\!] }
\newcommand{\psiT}{\omega}

\newcommand{\remi}[1]{\textcolor{black}{#1}}

\newcommand{\mario}[1]{\textcolor{black}{#1}}

\newtheorem{theorem}{Theorem}[section]
\newtheorem{definition}{Definition}[theorem]
\newtheorem{lemma}{Lemma}[theorem]
\newtheorem{proposition}{Proposition}[theorem]
\newtheorem{remark}{Remark}[theorem]

\makeindex             

\begin{document}

\title{Hyperbolic balance laws: \\residual distribution, local  and global fluxes}
\author{R\'emi Abgrall\\
Institute of Mathematics, University of Z\"urich, Switzerland,\\ {\tt remi.abgrall@math.uzh.ch}\\[10pt]
Mario Ricchiuto \\Team CARDAMOM, Inria Bordeaux Sud-Ouest, Talence, France,\\ {\tt mario.ricchiuto@inria.fr}}
%
%

\maketitle


\abstract{This \remi{review} paper describes a class of scheme named "residual distribution schemes" or "fluctuation splitting schemes". They are \mario{a generalization} 
of Roe's numerical flux \cite{roe1981}  in fluctuation form. The  so-called multidimensional fluctuation schemes have historically first been  developed for steady homogeneous hyperbolic systems. Their application to   unsteady problems and conservation laws has been really understood 
only relatively   recently. This understanding has allowed to make of the residual distribution framework a powerful playground to  develop numerical discretizations  embedding  some prescribed constraints.  
This paper describes in some detail
 these techniques, with several examples, ranging from the compressible Euler equations to the Shallow Water equations.}


\section{Introduction}
We are interested in the numerical approximation of \mario{partial differential equations relevant in fluid dynamics. For the objectives of the present paper, we will focus on}
the Euler and Navier-Stokes equations on complex domains, \mario{as well as on the}
shallow water equations. These \mario{models are particular cases of  the system of balance laws:}  
$$\dpar{\bu}{t}+\text{ div }\bbf(\bu)=\text { div }\bbf_v(\bu, \nabla \bu)+ S( \bu,\bx)$$
with initial and boundary conditions. The vector $\bu$ denotes \mario{a set of \emph{conserved} variables, which are often (but not always) densities of conserved macroscopic quantities (mass,energy, etc). For the Euler and Navier-Stokes equations we have}
$$\bu=(\rho, \rho \bbv, E)^T$$ 
where as usual, $\rho$ is the \mario{mass} density, $\bbv\in \R^d$ ($d=1,2,3$) the velocity, $E=e+\tfrac{1}{2}\rho \bbv^2$ is the total energy \mario{density}  with $e=e(\rho, p)$ is the internal energy, $p$ is the pressure. Here also, $\bbf(\bu)$ is the inviscid flux
$$\bbf(\bbu)=\begin{pmatrix}\rho \bv\\ \rho \bbv\otimes \bbv+p\text{Id}\\ \bv(E+p)\end{pmatrix}, $$
$\bbf_v$ is the viscous flux,  $S$ is a source term\remi{, and $\text{Id}$ is the $d\times d$ identity matrix.}

In the case of the shallow water equations, we have
$\bu=(H, H\bbv)^T$ where $H$ is the water height, the invicid flux is
$$\bbf(\bu)=\begin{pmatrix}
H\bbv\\
H\bbv\otimes \bbv + p(H)\text{ Id}\end{pmatrix}$$
with $p(H)=\tfrac{g}{2} H^2$, the viscous flux vanishes and the source term is
$$S( \bu,\bx)=\begin{pmatrix}
\mathbf{0}\\	
gH \nabla  b(x)+ c_{F}(\bu) H\bbv \end{pmatrix}$$
where $b$ is the bottom topography,  and $c_{F}(\bu)$ a friction coefficient modelling the effects of the boundary layer on the sea-floor.

In this paper, we are mostly interested in the inviscid case, so we will drop the viscous term . See however \cite{abgralldeSantisSISC,AbgralldeSantisNS}.
The general form we will discuss is
\begin{equation}
\label{eq:1}
\dpar{\bu}{t}+\text{ div }\bbf(\bu)=S( \bu,\bx)
\end{equation}
with initial and boundary conditions. 

Finite volume methods, \mario{and more recently Discontinuous Galerkin (DG) methods}, are very popular because they are known to be locally conservative. Thanks to Lax-Wendroff theorem, this is a very desirable property since it guaranty that "if all goes well", the limit solution is a weak solution of the problem. \mario{The extension of this notion to the non-homogenous case, and the construction finite volumes and DG schemes properly accounting for source terms
 is still a very open research subject, and over the years  several interesting approaches have been proposed.
Note that in this case  
system \eqref{eq:1}  admits non trivial solutions (steady and time dependent), and  the consistency with these solutions is often a desired property for the schemes.}
%
%
Indeed, there is a close connection between the two: in the Lax-Wendroff theorem, one essential condition on the numerical flux is that of consistency. When all the arguments of the numerical flux are equal, we must recover the continuous flux. This is another way of saying that uniform solutions must be preserved. When source terms are present this is not necessarily the case. For example, the flat free surface state
\begin{equation}\label{lar} 
H+b(x) =\eta_0=\text{const}
\end{equation}
with  still water $\bbv=0$ is undoubtedly physically more relevant than constant $\bu$.  For channels with smooth surfaces, if friction is neglected the constant flux and constant energy steady state
\begin{equation}\label{lcef}
 H\bbv = q_0=\text{const},\quad \dfrac{\bbv^2}{2}+g \big (H+b(x)\big )=\mathcal{E}_0=\text{const}
\end{equation} 
becomes the relevant one in general. This state is also compatible with the appearance of hydraulic jumps, across which the energy level is modified
(see e.g. \cite{M2AN_2009__43_2_333_0}).  In many other  applications however  friction cannot be neglected, and  the most general form of steady state is
obtained from the solution of
$$H v = q_0=\text{const},\quad  q_0^2/H + gH^2/2 + \int_{x_0}^x\big (gH \dpar{b}{x} + c_{F}(q_0,H)q_0\big)\;dx=q_0^2/H_0 +  gH_0^2/2     $$
If the bathymetry is linear, for example $b(x)=b_0-\xi_0 x$,  one can show that a constant state $\bu_0=(H_0,\;q_0)$ \remi{satisfies} the equilibrium  \cite{R15}
\begin{equation}\label{csf} 
-gH_0\xi_0 + c_{F}(q_0,H_0)q_0 =0
\end{equation} 
 For more general  definitions of $b(x)$, it is less apparent that a set of constant states can be associated to the steady equilibrium.

From the discrete approximation point of view, the problem is, how to modify   a given numerical flux  so that these states are  preserved, possibly within machine accuracy. 
This is often an ad-hoc construction, and if one has a different problem depending on the application, and the flux correction  needs to be re-written almost from scratch.

If one looks at the literature, there are other type of schemes. For example the stabilized variational methods using continuous finite elements such as the SUPG scheme \cite{Hughes1}, or the Galerkin scheme with jump stabilisation due to Burman et al.\cite{burman}. There are also the fluctuation splitting schemes that were initially designed by Roe and co-authors\cite{struijs}, and later extended to high order, steady and unsteady problems, as well as the shallow water equations \cite{abgrall99,energie,Mario,Mario2,Mario1,abgrallLarat,abgralldeSantisSISC,AbgralldeSantisNS,AR:17,Dec:17}. None of these scheme are formulated initially in terms of local fluxes, but they are working well. 
Indeed, in the numerical folklore, these schemes are often   claimed  not to be locally conservative, despite the contrary having been shown in
several  works  \cite{BurQS:10,HUGHES,AR:17,abgrall:conservation,abg_cons}. The most interesting aspect for these  is that treating \eqref{eq:1} with or without source term 
involves no major modification, and no special tricks.

The purpose of this paper is to recall that when $S=0$, residual distribution and continuous finite elements  are locally conservative. In fact, we will also recall how to construct an \emph{equivalent} flux formulation, 
and provide some explicit examples.  When $S\neq 0$, we  show how the schemes have been naturally extended to embed  the source term.
We  explain how to link them to more recent flux based formulations, although this is not how they are designed, and their construction is much more natural.
In one space dimension, we also recall that they have some  relation to the so-called path conservative schemes (see \cite{CASTRO2017131} and references therein).

The format of this paper is as follows. First we recall the main discrete prototype    we are interested in, written in a residual distribution form. In a second part, we show for steady problem that they have an equivalent flux formulation, and we explicitly construct the flux. In a third part, we extend this to unsteady problems. In a fourth part, we show (in 1D only) why these methods are agnostic to flux. Numerical examples are also given.

Throughout the paper, and for simplicity, we will not consider boundary conditions, even-though this is of course doable.

\section{Geometrical notations}\label{sec:notations}
\mario{Let us fix here the main notation used for the mesh and related geometrical entities.}
 The computational domain $\Omega$, $d=1,2,3$, is covered by a tessellation $\mathcal{T}_h$. We denote by $\mathcal{E}_h$ the set of internal edges/faces of $\mathcal{T}_h$, \remi{and by $\mathcal{F}_h$ the set of boundary faces}.  
 \mario{Mesh elements are generically denoted by $K$, while we use $e$ for a}
 face/edge $e\in \mathcal{E}_h\cup \mathcal{F}_h$.  The mesh is assumed to be shape regular, $h_K$ represents the diameter of the element $K$. Similarly, if $e\in \mathcal{E}_h\cup \mathcal{F}_h$, $h_e$ represents its diameter.

 Throughout this paper, we follow Ciarlet's definition \cite{ciarlet,ErnGuermond} of a finite element approximation: we have a set of degrees of freedom $\Sigma_K$ of linear forms acting on the set $\PP^k$ of polynomials of degree $k$ such that the linear mapping
 $$q\in \PP^k\mapsto \big (\sigma_1(q), \ldots, \sigma_{|\Sigma_K|}(q)\big )$$
 is one-to-one. The space $\PP^k$ is spanned by the basis function $\{\varphi_{\sigma}\}_{\sigma\in \Sigma_K}$  defined by
 $$\forall \sigma,\, \sigma',  \sigma(\varphi_{\sigma'})=\delta_\sigma^{\sigma'}.$$
  We have in mind either Lagrange interpolations where the degrees of freedom are associated to points in $K$, or other type of polynomials approximation such as B\'ezier polynomials where we will also do the same geometrical identification.
 Considering all the elements covering $\R^d$, the set of degrees of freedom is denoted by $\mathcal{S}$ and a generic degree of freedom  by $\sigma$. We note that for any $K$, 
 $$\forall \bx\in K, \quad \sum\limits_{\sigma\in K}\varphi_\sigma(\bx)=1.$$
 For any element $K$, $\#K$ is the number of degrees of freedom in $K$.

 The integer $k$ is assumed to be the same for any element.  We define 
$$\mathcal{V}^h=\bigoplus_K\{ \bv\in L^2(K), \bv_{|K}\in \PP^k\}$$
\remi{where $\PP^k$ is the set of polynomials of degree less or equal to $k$.}
The solution will be sought for in a  space $V^h$ that is:
\begin{itemize}
\item Either $V^h=\mathcal{V}^h$. In that case, the elements of $V^h$ can be discontinuous across internal faces/edges of $\mathcal{T}_h$. There is no conformity requirement on the mesh.
\item Or  $V^h=\mathcal{V}_h\cap C^0(\Omega)$ in which case the mesh needs to be conformal.
\end{itemize}

Throughout the text, we need to integrate functions. This is done via quadrature formula, and the symbol $\oint$ used in volume integrals
$$\oint_K v(\bx)\; d\bx$$
or boundary integrals
$$
\oint_{\partial K} v(\bx)\; d\gamma.$$
\mario{Note that the integration domain is uniquely defined by the limits of the integral, while} \remi{the symbol $\oint$}  \mario{is used here
 to  explicitly denote discrete integration } 
via user defined \mario{quadrature formulas.} 

 If $e\in \mathcal{E}_h$, represents any  \emph{internal} edge, i.e. $e\subset K\cap K^+$ for two elements $K$ and $K^+$,  we define for any function $\psi$ the jump  $\jump{\nabla \psi }=\nabla \psi_{|K}-\nabla \psi_{| K^+}$. Here the choice of $K$ and $K^+$ is important, and defines an orientation.
  Similarly, $\{\bv\}=\tfrac{1}{2}\big (\bv_{|K}+\bv_{|K^+}\big )$.
 
 If $\mathbf{x}$ and $\mathbf{y}$ are two vectors of $\R^q$, for $q$ integer, $\langle \bx,\by\rangle$ is their scalar product. In some occasions, it can also be denoted as $\bx\cdot\by$ or $\bx^T\by$.  We also use  $\bx\cdot\by$ when $\bx$ is a matrix and $\by$ a vector: it is simply the matrix-vector multiplication.

In section \ref{flux}, we have to deal with oriented graph. Given two vertices of this graph $\sigma$ and $\sigma'$, we write $\sigma>\sigma'$ to say that $[\sigma,\sigma']$ is a direct edge.

Section \ref{sec: bien roule}  mostly deals with one dimensional problems in $\R$. Here, a mesh is defined from a increasing sequence of points $\{x_\sigma\}_{\sigma\in \Z}$, the elements are the intervals
$$K_{\sigma+1/2}:=[x_\sigma, x_{\sigma+1}].$$
The lenght of $K_{\sigma+1/2}$ is $\Delta_{\sigma+1/2}x=x_{\sigma+1}-x_\sigma$ and we set
$$\Delta_\sigma x=\frac{\Delta x_{\sigma+1/2}+\Delta x_{\sigma-1/2}}{2}.$$

     \section{Example of schemes and conservation}\label{sec:conservation}
Let us provide several examples for the approximation of \eqref{eq:1}, to begin with, without source terms and in the steady case. They are
\begin{itemize}
\item The SUPG \cite{Hughes1} variational formulation,  with $\bu^h, \bv^h\in V^h=\mathcal{V}^h\cap C^0(\R^d)$:
\begin{equation}
\label{SUPG:var}
\begin{split}
a(\bu^h,\bv^h)&:=-\int_\Omega \nabla \bv^h\cdot \bF(\bu^h)\; d\bx+\sum\limits_{K\subset \Omega}h_K\int_K\big [ \nabla \bF(\bu^h)\cdot\nabla \bv^h\big ] \; \tau_K\; \big [\nabla\bF(\bu^h)\cdot \nabla \bu^h\big ] d\bx\\
&\qquad +\text{ Boundary terms}.
\end{split}
\end{equation}
Here $\tau_K$ is a positive parameter, or a positive definite matrix in the system case\footnote{More precisely such that it is symmetric definite positive up to a symetrization matrix, $A_0$ that, for fluid problems, is related to the Hessian of the entropy.}. 
\item The Galerkin scheme with jump stabilization,  see \cite{burman} for details.  We have
\begin{equation}
\label{burman:var}
\begin{split}
a(\bu^h,\bv^h)&:=-\int_\Omega \nabla \bv^h\cdot \bF(\bu^h)\; d\bx+\sum\limits_{e \subset \Omega}\theta_e h_e^2\int_e \big [\!\!\big[ \nabla \bv^h \big ]\!\!\big]\cdot \big [\!\!\big[ \nabla \bu^h\big ]\!\! \big] \; d\gamma \\
&\qquad +\text{ Boundary terms}.
\end{split}
\end{equation}
 Here,  $\bu^h, \bv^h\in V^h=\mathcal{V}^h\cap C^0(\Omega)$, and $\theta_e$ is a positive parameter.
\item The discontinuous Galerkin formulation: we look for $\bu^h, \bv^h\in V^h=\mathcal{V}^h$ such that
\begin{equation}\label{DG:var}
\begin{split}
a(\bu^h,\bv^h)&:=\sum\limits_{K\subset \Omega}\bigg ( -\int_K\nabla\bv^h\cdot\bbf(\bu^h) d\bx+\int_{\partial K}\bv^h\cdot \hat{\bbf}_\bn(\bu^{h},\bu^{h,+}) \;d\gamma \bigg )
\end{split}
\end{equation}
In \eqref{DG:var}, \remi{$\hat{\bbf}_\bn$ is a numerical flux consistant with the matrix $\bbf$ vector $\bn$ product $\bbf_\bn:=\bbf\cdot\bn$.} The boundary integral is a sum of integrals on the faces of $K$, and here for any face of $K$
 $\bu^{h,+}$ represents the approximation of $\bu$ on the other side of that face in the case of internal elements, and $\bu_b$ when that face is on $\partial \Omega$. 
 Note that to fully comply with \eqref{RD:scheme}, we should have defined for boundary faces $\bu^{h,+}=\bu^h$, and then \eqref{DG:var} is rewritten as 
 \begin{equation}\label{DG:var:2}
 \begin{split}
a(\bu^h,\bv^h):&=\sum\limits_{K\subset \Omega}\bigg ( -\int_K\nabla\bv^h\cdot\bbf(\bu^h) d\bx+\int_{\partial K}\bv^h\, \hat{\bbf}_\bn(\bu^{h},\bu^{h,+}) \; d\gamma \bigg )\\
&+\sum\limits_{\Gamma\subset\partial\Omega}\text{ Boundary term for }\Gamma.
\end{split}
\end{equation}
\item A finite volume scheme. Even though the discontinuous Galerkin method boils down into a finite volume scheme (where the volumes are the mesh elements) for degree $0$, let us provide a different example.
The notations are defined in Figure \ref{fig:fv}.
\begin{figure}[h]
\begin{center}
\subfigure[]{\includegraphics[width=0.45\textwidth]{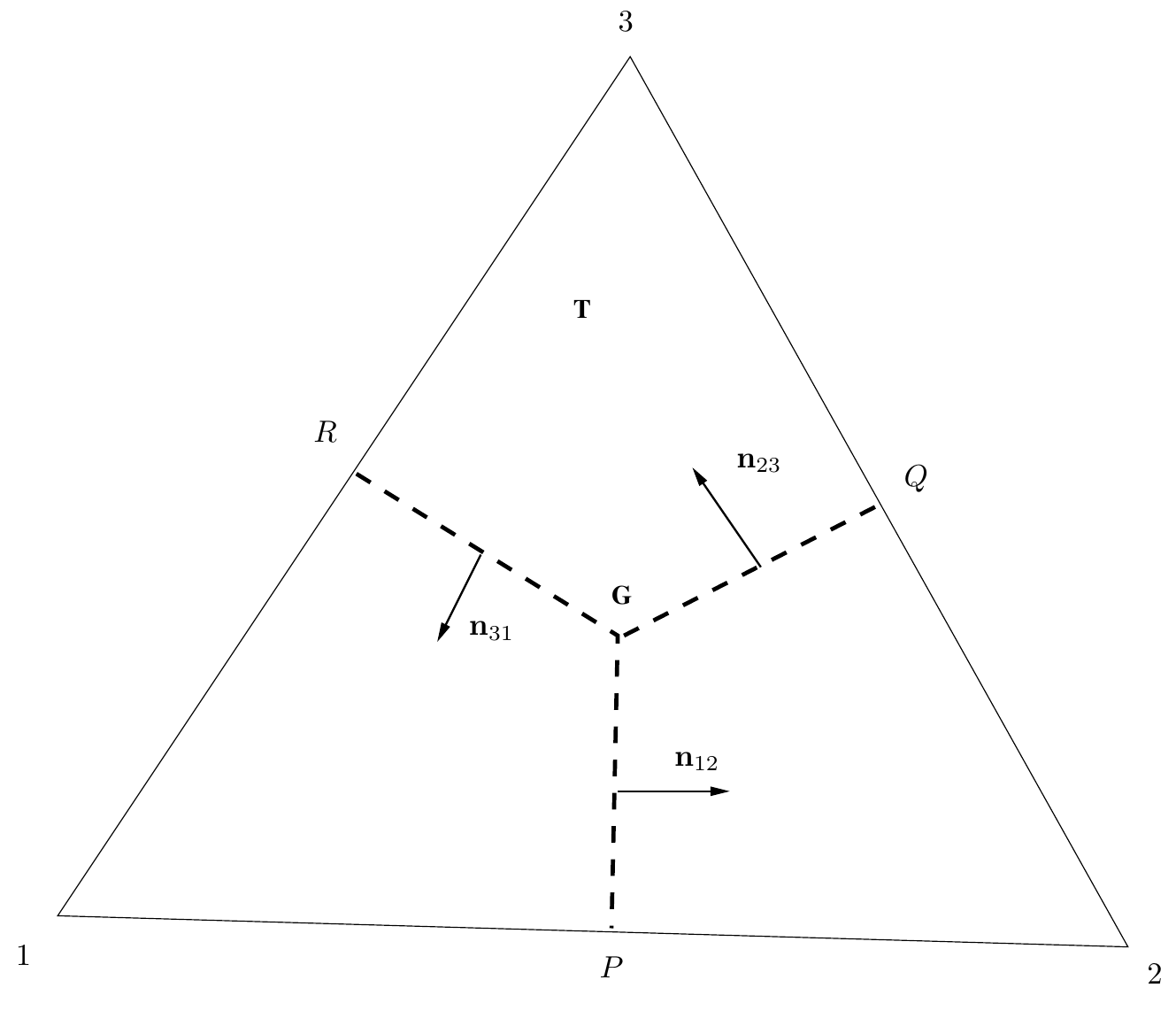}}
\subfigure[]{\includegraphics[width=0.45\textwidth]{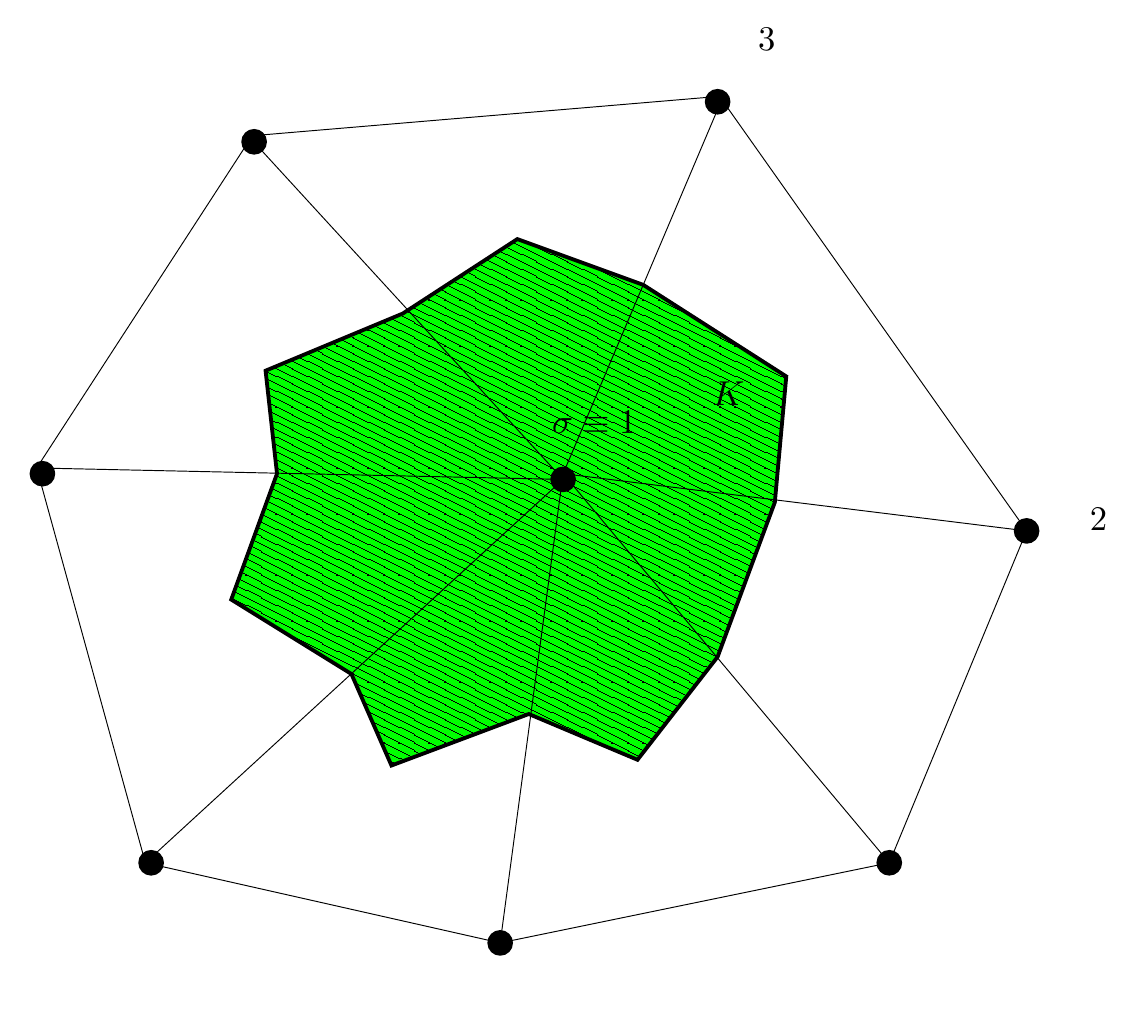}}
\end{center}
\caption{\label{fig:fv} Notations for the finite volume schemes. On the left: definition of the control volume for the degree of freedom $\sigma$.
 The vertex $\sigma$ plays the role of the vertex $1$ on the left picture for the triangle K. The control volume $C_\sigma$ associated to $\sigma=1$ is green on the right and corresponds to $1PGR$ on the left. The vectors $\bn_{ij}$ are normal to the internal edges scaled by the corresponding edge length}
\end{figure}
Again, we specialize ourselves to the case of triangular elements, but  \emph{exactly the same arguments} can be given for more general elements,
 provided a conformal approximation space can be constructed. This is  the case for
triangle elements, and we can take $k=1$.

The control volumes in this case are defined as the median cell, see figure \ref{fig:fv} and the scheme is 
\begin{equation*}
\sum_{\gamma \subset \partial C_\sigma} \hbbf_{\bn_\gamma }(\bu_\sigma , \bu^+ )= 0.
\end{equation*}
here we have taken a first order finite volume scheme, as it can be seen from the arguments of the numerical flux $ \hbbf_{\bn}$, however a high order extension with MUSCL extrapolation can equivalently be considered.
\end{itemize}

\bigskip
The interesting fact is that all these method can be rewriten in a unified manner, the residual distribution form.
{ In order to integrate the steady version of \eqref{eq:1} on a domain $\Omega\subset \R^d$, on each element $K$ and any degree of freedom $\sigma\in \mathcal{S}$ belonging to $K$,  we define residuals $\Phi_\sigma^K(\bu^h)$. Following \cite{abgrallLarat,abgralldeSantisSISC}, they are assumed to satisfy the following conservation relations:
For any element $K$, 
\begin{subequations}\label{conservation}
\begin{equation}
\label{conservation:K}
\sum\limits_{\sigma\in K}\!\Phi_\sigma^K(\bu^h)=\Phi^K(\bu^h) :=   \int_{\partial K}\!\hbbf_\bn(\bu^h,\bu^{h,+}) \;d\gamma
\end{equation}
where 
 $\Phi^K(\bu^h)$ is often referred to as  the ``element residual''.
 Note that if we denote by $\bbf^{h}$ the polynomial flux approximation within the element of maximum degree w.r.t. which the quadrature formulas  used in practice are 
 exact, we can recast   \eqref{conservation:K}  as 
\begin{equation}
\label{conservation:K-div} 
\sum\limits_{\sigma\in K}\Phi_\sigma^K(\bu^h)=\Phi^K(\bu^h)=
 \int_{\partial K}\! ( \hbbf_\bn(\bu^h,\bu^{h,+}) - \bbf_\bn^{h}) \;d\gamma  + \!\int\limits_K\!\nabla\cdot\bbf^{h} d\bx
 \end{equation}
which allows to write the element residual as the integral of the PDE plus  a boundary fluctuation.

In the case of a conformal mesh and  with  continuous elements, the  conservation condition becomes 
 \begin{equation}
\label{conservation:K_c}
\sum\limits_{\sigma\in K}\Phi_\sigma^K(\bu^h)=\Phi^K(\bu^h)=
 \int_{\partial K}\!  \bbf_\bn(\bu^h)   \;d\gamma  =  \!\int\limits_K\!\nabla\cdot\bbf^{h} d\bx
\end{equation} 
where we recall that the  difference between $\bbf(\bu^h)$ and  $\bbf^{h}$ is that the latter is the polynomial approximation within $K$ of the  highest possible degree
for which the   quadrature used  is exact. In general,   $\bbf(\bu^h(\bx))\ne \bbf^{h}(\bx)$.

The discretisation of \eqref{eq:1} is achieved via: for any $\sigma\in \mathcal {S}$,
\begin{equation}
\label{RD:scheme}
\sum\limits_{K\subset\Omega, \sigma\in K}\Phi_\sigma^K(\bu^h)+\text{ Boundary terms}=0.
\end{equation}
\end{subequations}
Concerning the boundary terms, they can be very naturally embedded by appropriately  embedding 
fluctuations on the boundary faces. We will omit this aspect for simplicity as it is mainly a  technical detail.


Using the fact that the basis functions that span $V_h$ have a \emph{compact} support, then each scheme can be rewritten in the form \eqref{RD:scheme} with the following expression for the residuals:
    \begin{itemize}
    \item For the  SUPG scheme \eqref{SUPG:var}, the  residual are defined by
    \begin{equation}\label{SUPG}
    \begin{split}\Phi_\sigma^K(\bu^h)&=\int_{\partial K}\varphi_\sigma \bF(\bu^h)\cdot \bn \; d\gamma -\int_K \nabla \varphi_\sigma\cdot \bF(\bu^h) \; d\bx\\
    &\qquad+h_K
    \int_K \bigg (\nabla_\bu\bF(\bu^h)\cdot \nabla \varphi_\sigma \bigg )\tau_K \bigg (\nabla_\bu\bF(\bu^h)\cdot \nabla \bu^h \bigg )\;d\bx .
    \end{split}\end{equation}
  \remi{  Note that in \eqref{SUPG} we have made an abuse of language that we will make systematically: to comply with the Gauss theorem and the form of \eqref{SUPG:var}, we have written
$$\int_\Omega \nabla \varphi_\sigma\cdot \bF(\bu^h)\; d\bx$$ where in the integral we have the 'product' of the vector $\nabla \varphi_\sigma$ with the matrix $\bF(\bu^h)$. It has to be understood
as $$\int_\Omega \bF(\bu^h)\cdot \nabla \varphi_\sigma\; d\bx.$$}
    \item For the Galerkin scheme with jump stabilization \eqref{burman:var}, the residuals are defined  by:  
        \begin{equation}\label{burman}\Phi_\sigma^K(\bu^h)=\int_{\partial K}\varphi_\sigma \bF(\bu^h)\cdot \bn\; d\gamma -\int_K \nabla \varphi_\sigma\cdot \bF(\bu^h)\; d\bx +
    \sum\limits_{\substack{e \text{ faces}\\\text{ of }K}} \frac{\theta_e}{2} h_e^2 \int\limits_{e} [\![\nabla \bu^h]\!]\cdot [\![\nabla \varphi_\sigma]\!]\; d\gamma.\end{equation}
    Here, since the mesh is conformal, any internal edge $e$ (or face in 3D) is the intersection of the element $K$ and another element denoted by $K^+$.
    \item For the discontinuous Galerkin scheme,
    \begin{equation}\label{DG}
    \Phi_\sigma^K(\bu^h)=-\int_K\nabla\varphi_\sigma\cdot\bbf(\bu^h) d\bx+\int_{\partial K}\varphi_\sigma\; \hat{\bbf}_\bn(\bu^{h},\remi{\bu^{h,+}}) \; d\gamma .
    \end{equation}
     \item For the  finite volume scheme    the fact that the boundary of the control volume is  closed implies that  the sum of the outward normals vanishes. So, we can define
    \begin{equation}
\begin{split}
\Phi_\sigma^K(\bu^h)&=\sum\limits_{\gamma\subset  \remi{\big (\partial C_\sigma\big )}\cap K} \big ( \hbbf_{\bn_\gamma }(\bu_\sigma, \bu^+ )-\bbf (\bu_\sigma)\cdot \bn_\gamma \big )\\
&=\sum\limits_{\gamma\subset \remi{\partial \big ( C_\sigma\cap K \big )}}  \hbbf_{\bn_\gamma }(\bu_\sigma, \bu^+ ).
\end{split}
\label{fv:res:sigma}
\end{equation}
We can then use the fact that on each facet separating   nodes $i, j$ local conservation implies  
$ \hbbf_{\bn_{ij}}(\bu_i,\bu_j) + \hbbf_{\bn_{ji}}(\bu_j,\bu_i)=0$, and
recover the elemental conservation relation as follows:
\begin{equation*}
\begin{split}
\sum_{\sigma\in K} \Phi_\sigma^K(\bu_h)&= \bigg ( \hbbf_{\bn_{12}}(\bu_1,\bu_2)-\hbbf_{\bn_{13}}(\bu_1,\bu_3)-\bbf(\bu_1)\cdot\bn_{12}+\bbf(\bu_1)\cdot\bn_{31}\bigg )\\
&+\bigg ( \hbbf_{\bn_{23}}(\bu_2,\bu_3)-\hbbf_{\bn_{12}}(\bu_2,\bu_1)+\bbf(\bu_2)\cdot\bn_{12}-\bbf(\bu_2)\cdot\bn_{23}\bigg )\\
&+\bigg ( -\hbbf_{\bn_{23}}(\bu_3,\bu_2)+\hbbf_{\bn_{31}}(\bu_3,\bu_1)-\bbf(\bu_3)\cdot\bn_{23}+\bbf(\bu_3)\cdot\bn_{31}\bigg )\\
&= \bbf(\bu_1)\cdot \big ( \bn_{12}-\bn_{31}\big ) +\bbf(\bu_2)\cdot \big ( -\bn_{23}+\bn_{31}\big )
+\bbf(\bu_3)\cdot \big ( \bn_{31}-\bn_{23}\big )\\
&=\bbf(\bu_1)\cdot\frac{\bn_1}{2}+\bbf(\bu_2)\cdot\frac{\bn_2}{2}+\bbf(\bu_3)\cdot\frac{\bn_3}{2}
\end{split}
\end{equation*}
where $\bn_j$ is the scaled inward normal of the edge opposite to vertex $\sigma_j$, i.e. twice the gradient of the $\PP^1$ basis function
 $\varphi_{\sigma_j}$ associated to this degree of freedom.
Thus, we can reinterpret the sum as the boundary integral of the Lagrange interpolant of the flux.
The finite volume scheme is then a residual distribution scheme with residual defined by \eqref{fv:res:sigma}
and a total residual defined by
\begin{equation}
\label{fv:tot:residu}
\Phi^K:=\int_{\partial K} \bbf^h_ \bn\;d\gamma , \qquad \bbf^h=\sum_{\sigma\in K} \bbf(\bu_\sigma)\varphi_\sigma.
\end{equation}
\item The residual distribution formalism  has also been used to build new schemes. \\
  A classical example  is  the nonlinear Lax-Friedrich's  discretization built to satisfy 
both  a high order truncation error estimate of order $\mathcal{O}(h^{p+1})$ for a polynomial approximation of degree $p$, and a positive-coefficient property 
\cite{Dec:17,AR:17}. The scheme reads 
$$\Phi_\sigma^K=\beta_\sigma^K \Phi^K$$
where the coefficients $\beta_\sigma$ are designed in such a way that the scheme is both monotonicity preserving and formally $k+1$-th order accurate if  a polynomial approximation of degree $k$ is used. This can be achieved in two steps as follows, see \cite{abgrallLarat} for the technical details
\begin{enumerate}
\item First evaluate the Rusanov (or Local lax-Friedrich residuals),
$$\Phi_\sigma^{LF, K}=\frac{\Phi^K}{N_K}+\alpha_K(\bu_\sigma-\bar \bu^K)
$$
where $\alpha_K$ is larger that the maximum on $K$ of $\Vert \nabla \bbf^h\Vert$, 
$N_K$ is the number of degree of freedom on $K$ and 
$\bar\bu^K$ is the arithmetic average of the $\bu_\sigma$ for $\sigma\in K$.
\item Define $x_\sigma$ as the ratio of ${\Phi_\sigma^{LF, K}}$ by ${\Phi^K}$, and 
$$\beta_\sigma=\dfrac{\max(x_\sigma,0)}{\sum_{\sigma'\in K} \max(x_{\sigma'}, 0)}.$$
Since\remi{$\sum\limits_{\sigma\in K}x_\sigma=1$}, we see that $\sum\limits_{\sigma'\in K} \max(x_{\sigma'}, 0)\geq 1$, so that there is no problem in the definition of this quantity as long as $\Phi^K\neq 0$. If $\Phi^K=0$, we can take any value since in the end
 the residual we are going to use is $\Phi_\sigma^K=\beta_\sigma^K \Phi^K$.
 
 Other expressions for $\beta_\sigma^K$ are feasible, but this is the one that is used in practice since it is very simple
 \end{enumerate}
This is not enough, as can be found in \cite{abgrallLarat}: the solution appears wiggly, especially in the smooth part of the solution. This is \emph{not}  a problem of stability. \remi{This occurs} because the scheme is over-compressing. One way to overcome this is to add some stabilizing/filtering  term. For example, in several papers a streamline upwind/least square term has been added to $\beta_\sigma^K \Phi^K$, namely
$$
\Phi_\sigma^{K\star}= \beta_\sigma^K \Phi^K+h_K\oint_K \big [ \nabla\bbf^h(\bu^u)\cdot\nabla\varphi_\sigma \big ]\tau_K \big [\nabla\bbf^h(\bu^u)\cdot \nabla \bu^h\big ] \; d\bx$$
The resulting scheme is referred to later on in the paper as LLFs (the first "L" for Limited, the "s" for stabilized).
In \cite{RD:quad} is discussed the choice of minimal quadrature formula for the evaluation of the integral term.
Adding the least square term destroys in principle the maximum preserving property of the method, however in practice it does not, this is why we call this essentially non oscillatory RD scheme: the least square term acts as a mild filter of the spurious modes. Another possible technique to achieve this filtering is to use jump terms as in \eqref{burman}.

In the case of system, one can extends the construction by using a characteristic decomposition of the residual, see again \cite{abgrallLarat}. Last, any monotone first order residual can be used in the step 1 of the construction, not only the Local Lax-Friedrichs one.
  \end{itemize}
All these residuals satisfy the relevant conservation relations, namely \eqref{conservation:K}, depending if we are dealing with element residuals or boundary residuals.

It can be shown, see \cite{abgrallRoe} that a scheme defined by \eqref{conservation} satisfies a Lax-Wendroff like theorem: if the mesh is regular, if the numerical sequence is bounded in $L^\infty$ and if a subsequence converges in $L^2$ (for example) to a $\bv\in L^2$, then this function is a weak solution of \eqref{eq:1}. \remi{A similar result holds on the entropy if an entropy inequality exists.}

\section{Flux formulation of Residual Distribution schemes}\label{flux}
Conversely, we show in this section that any scheme \eqref{RD:scheme} also admits a flux formulation: \emph{the method is also locally conservative}. In addition, we provide  an \emph{explicit} form of the flux. Local conservation is of course  well known for the Finite Volume and discontinuous Galerkin  approximations. It is much less understood for the  continuous finite elements methods, despite the papers \cite{Hughes1,BurQS:10}.
 This question of finding an equivalent flux formulation amounts to defining control volumes and flux functions. 

We first have to adapt the notion of consistency. We define a multidimensional flux as follows:
\begin{definition}\label{MD:consistency}
 A multidimensional flux 
$\hbbf_\bn:=\hbbf_\bn(\bu_1, \ldots , \bu_N)$
is consistent if, when $\bu_1= \bu_2= \ldots = \bu_N=\bu$ then
$\hbbf_\bn(\bu, \ldots , \bu)=\bbf(\bu)\cdot \bn.$
\end{definition}
We proceed first with the general case and show the connection with elementary fact about graphs, and then provide several examples. The results of this section apply to any 
finite element method but also to discontinuous Galerkin methods. There is no need for the exact evaluation of integral formula (surface or boundary), so that these results apply to schemes as they are implemented.

\bigskip

We consider the general case, i.e when $K$ is a polytope contained in $\R^d$ with degrees of freedoms on the boundary of $K$. The set $\mathcal{S}$ is the set of degrees of freedom. We construct  a triangulation $\mathcal{T}_K$\remi{, i.e. a graph, } of $K$ whose vertices  are exactly the elements of $\mathcal{S}$. Choosing an orientation of $K$, it is propagated on $\mathcal{T}_K$: the edges are oriented.

\remi{Inspired by \eqref{fv:tot:residu}, starting from the  graph,  we will construct control volumes and flux. Again inspired by the finite volume example, and since the shape of the control volumes is still unknown, we label the flux by the edges of the graph, and hence we slightly change notations.} The problem is to find quantities $\hbbf_{\sigma,\sigma'}$ for any edge $[\sigma,\sigma']$ of   $\mathcal{T}_K$ such that:
\begin{subequations}\label{GC:1}
\begin{equation}
\label{GC:1.1}
\Phi_\sigma=\sum_{\text{ edges }[\sigma,\sigma']} \hbbf_{\sigma,\sigma'}+\hbbf_\sigma^{b}
\end{equation}
with, \remi{since permuting two vertices amounts to change the orientation,}
\begin{equation}
\hbbf_{\sigma,\sigma'}=-\hbbf_{\sigma',\sigma}
\label{GC:1.2}
\end{equation}
and $\hbbf_\sigma^{b}$ is the 'part' of $\oint_{\partial K} \hbbf_\bn(\bu^h,\bu^{h,+}) \; d\gamma$ associated to $\sigma$. {The control volumes will be defined by their normals so that we get consistency.} \remi{The normal corresponding to the edge $[\sigma, \sigma']$ is denoted by $\bn_{\sigma,\sigma'}$.}

Note that \eqref{GC:1.2} implies the conservation relation
\begin{equation}
\label{GC:conservation}
\sum\limits_{\sigma\in K}\Phi_\sigma=\sum\limits_{\sigma\in K}\hbbf_\sigma^b.
\end{equation}
In this paper, we will take
\begin{equation}
\label{BC:1.3}
\hbbf_\sigma^b=\oint_{\partial K} \varphi_\sigma\; \hbbf_\bn (\bu^h,\bu^{h,+}) \; d\gamma,
\end{equation}
\end{subequations}
but other  examples can be considered provided the consistency \eqref{GC:conservation} relation holds true, see \cite{abgrall:conservation}.
Any edge $[\sigma,\sigma']$ is either direct or, if not, $[\sigma',\sigma]$ is direct. Because of \eqref{GC:1.2}, we only need to know $\hbbf_{\sigma,\sigma'}$ for direct edges. Thus we introduce the notation $\hbbf_{\{\sigma,\sigma'\}}$ for  the flux  assigned to  the direct edge whose extremities are $\sigma$ and $\sigma'$. We can rewrite \eqref{GC:1.1} as, for any $\sigma\in \mathcal{S}$,
\begin{equation}
\label{GC:1.1bis}
\sum_{\sigma'\in \mathcal{S}} \varepsilon_{\sigma,\sigma'} \hbbf_{\{\sigma,\sigma'\}}=\Psi_\sigma:=\Phi_\sigma-\hbbf_\sigma^b,
\end{equation}
with $$
\varepsilon_{\sigma,\sigma'}=\left \{
\begin{array}{ll}
0& \text{ if }\sigma \text{ and }\sigma' \text{ are not on the same edge of }\mathcal{T},\\
1& \text{ if } [\sigma,\sigma']\text{ is an edge and } \sigma \rightarrow \sigma' \text{ is direct,}\\
-1&  \text{ if } [\sigma,\sigma']\text{ is an edge and } \sigma' \rightarrow \sigma \text{ is direct.}
\end{array}
\right .
$$
$\mathcal{E}^+$ represents the set of direct edges.

Hence the problem is to find  a vector $\hbbf=(\hbbf_{\{\sigma,\sigma'\}})_{\{\sigma,\sigma'\} \text{ direct edges}}$ such that
$$A\hbbf=\Psi$$
where $\Psi=(\Psi_\sigma)_{\sigma\in \mathcal{S}}$ and $A_{\sigma \sigma'}=\varepsilon_{\sigma,\sigma'}$.

We have  the following lemma which shows the existence of a solution. Its proof can be found in \cite{abgrall:conservation}.
\begin{lemma}\label{lemma:flux}
For any couple $\{\Phi_\sigma\}_{\sigma\in \mathcal{S}}$ and $\{\hbbf_\sigma^{b}\}_{\sigma\in \mathcal{S}}$ satisfying the condition  \eqref{GC:conservation}, there exists numerical flux functions $\hbbf_{\sigma,\sigma'}$ that satisfy \eqref{GC:1}. Recalling that the  matrix of the Laplacian of the graph is $L=AA^T$, we have
\begin{enumerate}
\item The rank of $L$ is $|\mathcal{S}|-1$ and its image is $\big (\text{span}\{\mathbf{1}\})^\bot$. We still denote the inverse of $L$ on $\big (\text{span}\{\mathbf{1}\} )^\bot$ by $L^{-1}$,
\item 
With the previous notations, a solution is 
\begin{equation}
\label{eq:lemma}\big (\hbbf_{\{\sigma,\sigma'\}}\big )_{\{\sigma,\sigma'\} \text{ direct edges}}=A^TL^{-1} \big (\Psi_\sigma\big )_{\sigma\in \mathcal{S}}.\end{equation}
\end{enumerate}
\end{lemma}

This set of flux are consistent and we can estimate the normals $\bn_{\sigma,\sigma'}$.
In the case of a constant state, we have $\Phi_\sigma=0$ for all $\sigma\in K$. Let us assume that
\begin{equation}
\label{GC:consistency}
\hbbf_\sigma^b=\bbf(\bu^h)\cdot \mathbf{N}_\sigma
\end{equation}
with $\sum\limits_{\sigma\in K} \mathbf{N}_\sigma=0$: this is the case for all the examples we consider. If $\hbbf_\sigma$ is defined by \eqref{BC:1.3}, we see that
$$\mathbf{N}_\sigma=\oint_{\partial K} \varphi_\sigma \bn\; d\gamma.$$
The flux $\bbf(\bu^h)$ has components on the canonical basis of $\R^d$: $\bbf(\bu^h)=\big (f_1(\bu^h), \ldots , f_d(\bu^h)\big )$, so that from \eqref{GC:consistency}, we get
$$\hbbf_\sigma^b=\sum\limits_{i=1}^d f_i(\bu^h)\mathbf{N}^i_\sigma.$$
Applying this to $\big (\hbbf_{\sigma_1}^b, \ldots, \hbbf_{\sigma_{\#K}}^b\big )$, we see that the $j$-th component of $\bn_{\sigma,\sigma'}$ for $[\sigma,\sigma']$ direct, must satisfy:
$$\text{ for any }\sigma\in K, \; \mathbf{N}^j_\sigma=\sum\limits_{[\sigma,\sigma']\text{ edge }}\varepsilon_{\sigma,\sigma'}\bn_{\sigma,\sigma'}^j$$
i.e.
$$\big ( \mathbf{N}^j_{\sigma_1}, \ldots , \mathbf{N}^j_{\sigma_{\#K}}\big )^T= A \; \big ( \bn_{\sigma,\sigma'}^j\big )_{[\sigma,\sigma']\in \mathcal{E}^+}.$$
We can solve the system and the solution, with some abuse of language, is
\begin{equation}
\label{GC:normals}
\big ( \bn_{\sigma,\sigma'}\big )_{[\sigma,\sigma']\in \mathcal{E}^+}=A^TL^{-1} \big ( \mathbf{N}_{\sigma_1}, \ldots , \mathbf{N}_{\sigma_{\#K}}\big )^T
\end{equation}
This also defines the control volumes since we know their normals. We can state:
\begin{proposition}
If the residuals $(\Phi_\sigma)_{\sigma\in K}$ and the boundary fluxes $(\hbbf_\sigma^b)_{\sigma\in K}$ satisfy \eqref{GC:conservation}, and if the boundary fluxes satisfy the consistency relations \eqref{GC:consistency}, then we can find  a set of consistent flux $(\hbbf_{\sigma,\sigma'})_{[\sigma,\sigma']} $ satisfying \eqref{GC:1}. They are given by \eqref{eq:lemma}. In addition, for a constant state,
$$\hbbf_{\sigma,\sigma'}(\bu^h)=\bbf(\bu^h)\cdot\bn_{\sigma,\sigma'}$$ for the normals defined by \eqref{GC:normals}.
\end{proposition}

\bigskip

We can state a couple of general remarks:
\begin{remark} $ $
\begin{enumerate}
\item The flux $\hbbf_{\sigma,\sigma'}$ depend on the $\Psi_\sigma$ and not directly on the $\hbbf_\sigma^b$. We can design the fluxes independently of the boundary flux, and their consistency directly comes from the consistency of the boundary fluxes.
\item 
The residuals depends on more than 2 arguments. For stabilized finite element methods, or the non linear stable residual distribution
 schemes, see e.g.  \cite{Hughes1,struijs,abgrallLarat}, the residuals depend on all the states on  $K$. Thus
the formula \eqref{eq:lemma} shows that the flux depends on more than two states in contrast to the  1D case. In the finite volume case however, the support of the flux function is generally larger than the three states of $K$, think for example of an ENO/WENO method, or a simpler MUSCL ones.
\item The formula \eqref{eq:lemma} make no assumption on the approximation space $V^h$: they are valid for continuous and discontinuous approximations. The structure of the approximation space appears only in the total residual.
\end{enumerate}

\end{remark}

\bigskip
Let us give two examples, that will be valid for SUPG and the Galerkin scheme with stabilisation because the explicit form of the residual does not play any role.

Let $K$ be a fixed triangle. \remi{The degrees of freedom (the vertices) will be denoted by $\{\sigma\}_{\sigma_\in K}$ or $\{\sigma_i\}_{i=1,2,3}$ or their label in $\{1,2,3\}$.} We are given a set of residues $\{\Phi_\sigma^K\}_{\sigma\in K}$, our aim here is to define a
 flux function such that relations similar to \eqref{fv:res:sigma} hold true. The \remi{adjacency} matrix is 
$$A=\left (\begin{array}{rrr} 1&0&-1\\
-1&1&0\\
0&-1&1
\end{array}\right ).
$$
A straightforward calculation shows that the matrix $L=A^TA$ has eigenvalues $0$ and $3$ with multiplicity 2 with eigenvectors
$$R=\begin{pmatrix}
\frac{1}{\sqrt{3}} & \frac{1}{\sqrt{2}} & \frac{1}{\sqrt{6}}\\
\frac{1}{\sqrt{3}} & \frac{-1}{\sqrt{2}} & \frac{1}{\sqrt{6}}\\
\frac{1}{\sqrt{3}} & 0                         & \frac{-2}{\sqrt{6}}
\end{pmatrix}
$$
To solve $A\hbbf=\Psi$, we decompose $\Psi$ on the eigenbasis:
$$\Psi=\alpha_2 R_2+\alpha_3R_3$$ where explicitly
$$\begin{array}{l}
\alpha_2=\frac{1}{\sqrt{2}} \big (\Psi_1-2\Psi_2+\Psi_3\big )\\
\\
\alpha_3=\sqrt{\frac{3}{2}}\big (\Psi_1-\Psi_3\big )
\end{array}
$$
so that 
$$
\hbbf_{12}=\frac{1}{3}\big (\Psi_1-\Psi_2\big ), \quad
\hbbf_{23}= \frac{1}{3}\big (\Psi_2-\Psi_3\big ),\quad
\hbbf_{32}=\frac{1}{3}\big (\Psi_3-\Psi_1\big )
.$$

In order to describe the control volumes, we first have to make precise the normals $\bn_\sigma$ in that case. It is easy to see that in all the cases described above, we have 
$$\normal_\sigma=-\frac{\bn_\sigma}{2}.$$ Then a short calculation shows that
$$\begin{pmatrix}
\bn_{12} \\ \bn_{23} \\ \bn_{31} \end{pmatrix}=
\frac{1}{6}\begin{pmatrix} \bn_1-\bn_2 \\ \bn_2 -\bn_3 \\ \bn_3 -\bn_1 \end{pmatrix}.
$$
Using elementary geometry of the triangle, we see that these  are the normals of the elements of the dual mesh. For example, the normal $\bn_{12}$ is the normal of 
 $PG$, see figure \ref{fig:fv}.
 \begin{figure}[h!]
\begin{center}
\includegraphics[width=0.45\textwidth]{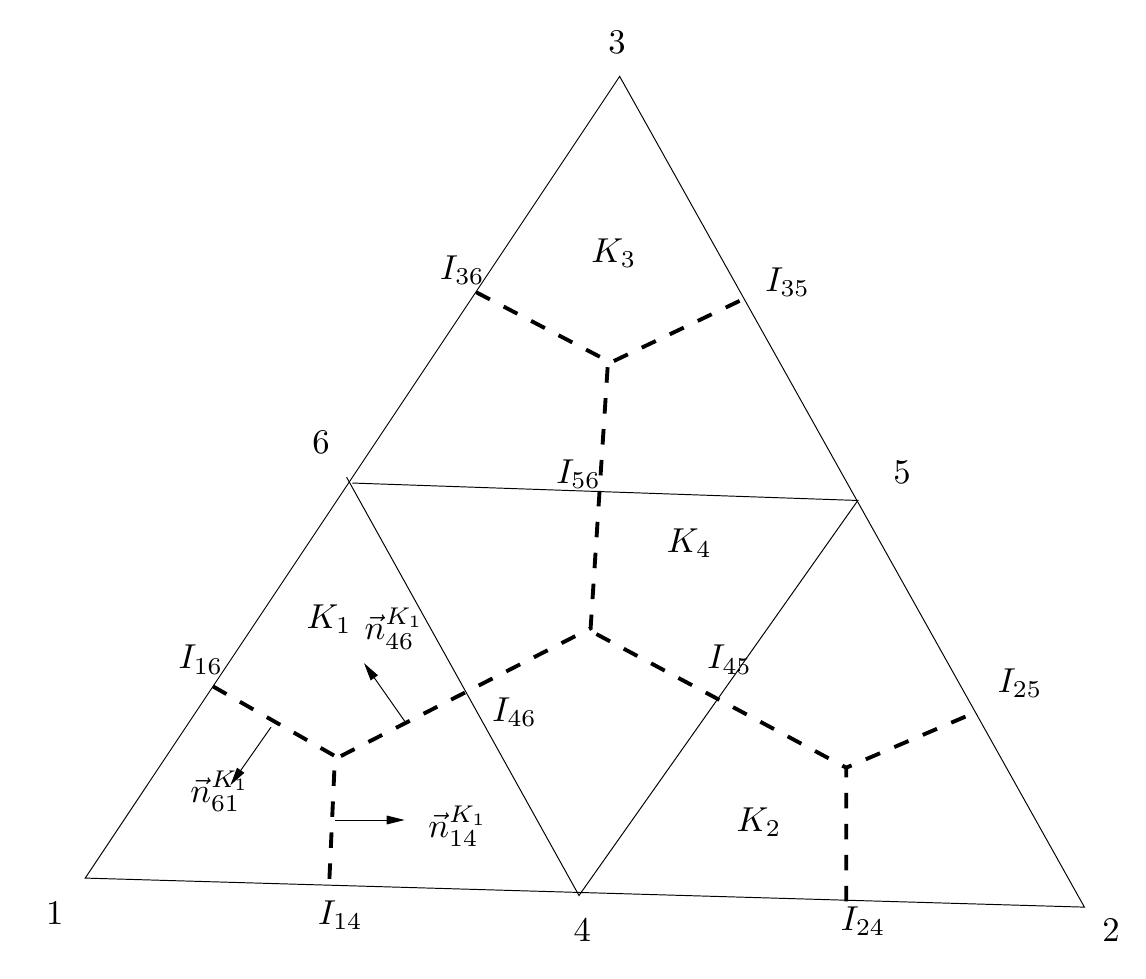}
\end{center}
\caption{\label{fig:P2} Geometrical elements for the $\PP^2$ case. $I_{ij}$ is the mid-point between the vertices $i$ and $j$. The intersections of the dotted lines are the centroids of the sub-elements.}
\end{figure}

 In the case of a quadratic approximation a similar set of formula can be given. Again this will be valid for SUPG and Galerkin with jumps.
Using a similar method, we get (see figure \ref{fig:P2} for some notations): {\color{black}
\begin{equation*}
\begin{split}
 \hbbf_{14}&=\frac{ 5 \Psi_1}{12}-\frac{ 5 \Psi_2}{36}-\frac{ \Psi_3}{36}-\frac{ 7 \Psi_4}{36}-\frac{ \Psi_5}{12}+\frac{ \Psi_6}{36}\\
\hbbf_{42}&=\frac{ 5 \Psi_1}{36}-\frac{ 5 \Psi_2}{12}+\frac{ \Psi_3}{36}+\frac{ 7 \Psi_4}{36}-\frac{ \Psi_5}{36}+\frac{ \Psi_6}{12}\\
\hbbf_{25}&=-\frac{ \Psi_1}{36}+\frac{ 5 \Psi_2}{12}-\frac{ 5 \Psi_3}{36}+\frac{ \Psi_4}{36}-\frac{ 7 \Psi_5}{36}-\frac{ \Psi_6}{12}\\
\hbbf_{53}&=\frac{ \Psi_1}{36}+\frac{ 5 \Psi_2}{36}-\frac{ 5 \Psi_3}{12}+\frac{ \Psi_4}{12}+\frac{ 7 \Psi_5}{36}-\frac{ \Psi_6}{36}\\
\hbbf_{36}&=-\frac{ 5 \Psi_1}{36}-\frac{ \Psi_2}{36}+\frac{ 5 \Psi_3}{12}-\frac{ \Psi_4}{12}+\frac{ \Psi_5}{36}-\frac{ 7 \Psi_6}{36}\\ 
\hbbf_{61}&=-\frac{ 5 \Psi_1}{12}+\frac{ \Psi_2}{36}+\frac{ 5 \Psi_3}{36}-\frac{ \Psi_4}{36}+\frac{ \Psi_5}{12}+\frac{ 7 \Psi_6}{36}\\
\hbbf_{64}&=\frac{ \Psi_1}{9}-\frac{ \Psi_3}{9}+\frac{2 \Psi_4}{9}-\frac{2 \Psi_5}{9}\\
\hbbf_{45}&=-\frac{ \Psi_1}{9}+\frac{ \Psi_2}{9}+\frac{2 \Psi_5}{9}-\frac{2 \Psi_6}{9}\\
\hbbf_{56}&=-\frac{ \Psi_2}{9}+\frac{ \Psi_3}{9}-\frac{2 \Psi_4}{9}+\frac{2 \Psi_6}{9}
\end{split}
\end{equation*}
Then we choose the boundary flux:
$$\hbbf_\sigma^b=\int_{\partial K}\varphi_\sigma\bbf(\bu^h)\cdot\bn\; d\gamma$$ and get:
$$
\begin{array}{lll}
\normal_l=-\dfrac{\bn_l}{6} & \text{if }l=1,2,3\\ &&\\
\normal_4=\dfrac{\bn_3}{3}& \normal_5=\dfrac{\bn_1}{3}& \normal_6=\dfrac{\bn_2}{3}
\end{array}
$$
The normals are given by:
\begin{equation*}
\begin{split}
 \bn_{14}&=\frac{ 5 \normal_1}{12}-\frac{ 5 \normal_2}{36}-\frac{ \normal_3}{36}-\frac{ 7 \normal_4}{36}-\frac{ \normal_5}{12}+\frac{ \normal_6}{36}\\
 \bn_{42}&=\frac{ 5 \normal_1}{36}-\frac{ 5 \normal_2}{12}+\frac{ \normal_3}{36}+\frac{ 7 \normal_4}{36}-\frac{ \normal_5}{36}+\frac{ \normal_6}{12}\\
\bn_{25}&=-\frac{ \normal_1}{36}+\frac{ 5 \normal_2}{12}-\frac{ 5 \normal_3}{36}+\frac{ \normal_4}{36}-\frac{ 7 \normal_5}{36}-\frac{ \normal_6}{12}\\
\bn_{53}&=\frac{ \normal_1}{36}+\frac{ 5 \normal_2}{36}-\frac{ 5 \normal_3}{12}+\frac{ \normal_4}{12}+\frac{ 7 \normal_5}{36}-\frac{ \normal_6}{36}\\
\bn_{36}&=-\frac{ 5 \normal_1}{36}-\frac{ \normal_2}{36}+\frac{ 5 \normal_3}{12}-\frac{ \normal_4}{12}+\frac{ \normal_5}{36}-\frac{ 7 \normal_6}{36}\\ 
\bn_{61}&=-\frac{ 5 \normal_1}{12}+\frac{ \normal_2}{36}+\frac{ 5 \normal_3}{36}-\frac{ \normal_4}{36}+\frac{ \normal_5}{12}+\frac{ 7 \normal_6}{36}\\
\bn_{64}&=\frac{ \normal_1}{9}-\frac{ \normal_3}{9}+\frac{2 \normal_4}{9}-\frac{2 \normal_5}{9}\\
\bn_{45}&=-\frac{ \normal_1}{9}+\frac{ \normal_2}{9}+\frac{2 \normal_5}{9}-\frac{2 \normal_6}{9}\\
\bn_{56}&=-\frac{ \normal_2}{9}+\frac{ \normal_3}{9}-\frac{2 \normal_4}{9}+\frac{2 \normal_6}{9}
\end{split}
\end{equation*}
}

The case of the discontinuous Galerkin can be worked out similarly.

\section{Embedding source terms: well balancing and global fluxes}\label{sec: bien roule}
\subsection{The one-dimensional case}
We look now into the approximation of solutions to  the steady limit of  \eqref{eq:1} in one dimension :
$$\dpar{\bbf(\bu)}{x}=S(\bu,x).$$
Despite the apparent simplicity of the problem, it is well known that   
some  fundamental change of paradigm is required  compared to conservation laws.  In particular, the non-autonomous character
of the problem,  associated to the presence of the source term $S(\bu,x)$ requires a more general notion of consistency.

The examples provided in the introduction for the shallow water equations show that these non-trivial states can only in some cases be
characterized by a set of physically relevant invariants. A possible way out to replace the notion of \emph{consistency with constant states}
is to  introduce an (unknown) ``source flux'' $\bbs$  as
$$\bbs(x)=\int_{x_0}^x S(\bu(x),x)\; dx$$
One can  now argue that a more relevant notion of steady states is the one associated to a constant global flux  
$$\bbg = \bbf - \bbs = \bbg_0 = \text{const}$$
Although  several works  have proposed   explicit  constructions  of the local values of  $\bbs$ \cite{DK03,DK05,CCHKT19},
this is essentially possible only in a 1D setting or by means of some dimension by dimension splitting.  
The main issue is how to construct schemes consistent with the notion of a constant global flux, without necessarily
having its explicit knowledge.

This issue is dealt with very naturally in the residual distribution setting. Let us focus for the moment on  continuous approximations on conformal meshes.
The natural way to proceed is to   generalize  the notion of conservation 
defined by  \eqref{conservation:K_c}  by including the whole PDE in it:
 \begin{equation}
\label{conservation:K_cS}
\sum\limits_{\sigma\in K}\Phi_\sigma^K(\bu^h)=\Phi^K(\bu^h)=
 \!\int_K ( \nabla\cdot\bbf^{h}  - S^h)\,dx
\end{equation} 
where $S^h$ is a discrete  approximation of the source within the element compatible with a certain quadrature strategy.  We will discuss this in aspect in some detail  shortly.
For the moment, let us consider the 1D residual distribution scheme, seeking the steady solution as the limit of
$$
 \Delta x_{\sigma}
\dfrac{\bu^{n+1}_{\sigma} - \bu^{n}_{\sigma}}{\Delta t}  + \sum\limits_{K, \sigma\in K} \Phi_\sigma^K(\bu^{h,n})=0
$$
\remi{We focus on the $\PP^1$ case to begin with, but the extension to higher polynomials can be obtained similarly to what we discussed in section \S \ref{flux}.
Each node will receive two contributions, $\Phi_\sigma^{\sigma-1/2}$ from the element on its left, $K_{\sigma-1/2}$, the second, $\Phi_\sigma^{\sigma+1/2}$ from the element on its right, $K_{\sigma+1/2}$. The conservation relation \eqref{conservation:K_cS} in one dimension and the $\PP^1$ setting simply writes}
$$\Phi_\sigma^{\sigma+1/2}+\Phi_{\sigma+1}^{\sigma+1/2}=\bbf(\bu_{\sigma+1})-\bbf(\bu_\sigma)-\int_{K_{\sigma+1/2}} S^h(x, \bbu^h)\; dx.$$
\remi{We have again adapted the notations to make them less heavy in the 1D case.}
We can proceed as follows. First we set
$$
\remi{S_{\sigma+1/2}:= \dfrac{1}{|K_{\sigma+1/2}|} \!\int_{K_{\sigma+1/2}}  S^h(x, \bbu^h)\,dx=\dfrac{1}{\Delta_{\sigma+1/2}x}\int_{x_\sigma}^{x_{\sigma+1}} S^h(x, \bbu^h)\,dx.}
$$ 
Next we define
$$
\bbs_{\sigma+1} = \bbs_{\sigma} + \Delta_{\sigma+1/2}x\; S_{{\sigma+1/2}}(x, \bbu^h) 
$$
with $ \bbs_{\sigma_0}=0$ for a given but arbitrary $\sigma_0$.  We then set $\bbg_{\sigma+1}=\bbf_{\sigma+1}-\bbs_{\sigma+1}$, and  recast the iterations as  
$$
 \Delta_{\sigma}x
\dfrac{\bu^{n+1}_{\sigma} - \bu^{n}_{\sigma}}{\Delta t}  + (\bbg_{\sigma} + \Phi_\sigma^{K_{\sigma+1/2}})  - ( \bbg_{\sigma}- \Phi_\sigma^{K_{\sigma-1/2}}  )=0
$$
Finally we set  
$$
\hbbg_{\sigma\pm1/2}:=\bbg_{\sigma}  \pm  \Phi_\sigma^{{\sigma\pm1/2}}
$$
which is a consistent numerical global flux. Note however that $\hbbg_{\sigma\pm1/2}$ is never explicitly built in the residual distribution approach!

To give a few examples, \mario{let us start from the centered splitting} 
$$
\Phi_\sigma^{{\sigma \pm1/2}} = \dfrac{1}{2}  \Phi^{{\sigma \pm1/2}}= \dfrac{1}{2}( \Delta\bbf_{\sigma\pm 1/2} - \Delta_{\sigma+1/2}x\; S_{{\sigma+1/2}})
$$
\mario{We can easily show that this splitting  leads to an equivalent global finite volume flux}  
$$
\hbbg_{\sigma\pm1/2}:=\bbg_{\sigma}  \pm \dfrac{1}{2} \Phi^{K_{\sigma\pm1/2}} = \bbf_{\sigma}  - \bbs_{\sigma} +  \dfrac{\Delta\bbf_{\sigma\pm 1/2}}{2} - \dfrac{\Delta\bbs_{\sigma\pm 1/2} }{2}= 
 \dfrac{\bbg_{\sigma}+\bbg_{\sigma\pm1}}{2}
$$
\mario{Similarly, the} Galerkin scheme can be shown to be equivalent to the \mario{finite volume scheme with global numerical flux given by}  
$$
\hbbg^{\text{Gal}}_{\sigma\pm1/2}=   \dfrac{\bbg_{\sigma}+\bbg_{\sigma\pm1}}{2} \pm \int_{K_{\sigma\pm1/2}}\big(\varphi_{\sigma}-\dfrac{1}{2}\big)\left(\dpar{\bbf^h}{x} -S^h \right)\,dx\;,
$$
and for the SUPG we have
$$
\hbbg^{\text{SUPG}}_{\sigma\pm1/2}=   \dfrac{\bbg_{\sigma}+\bbg_{\sigma\pm1}}{2} \pm \int_{K_{\sigma\pm1/2}}\left[\big(\varphi_{\sigma}-\dfrac{1}{2}\big)+
\nabla_{\bu}\bbf\,  \dpar{\varphi_{\sigma}}{x}\;\tau_{K_{\sigma\pm1/2}}
\right]
\left(\dpar{\bbf^h}{x} -S^h \right)\,dx\;,
$$
In general, \mario{we can follow  for example section \ref{flux},   consider test functions $\{\omega_\sigma\}$ defining a partition of unity for conservation porposes,  and set}
\begin{equation}\label{phi-rd-source}
\Phi_\sigma^{K}=\int_{K}\omega_{\sigma}\left(\nabla\cdot  \bbf^h  -S^h \right)\,dx,
\end{equation}
\mario{This scheme} is equivalent in 1D to the finite volume global flux method defined by
\begin{equation}\label{g-flux-RD}
\hbbg^{\text{RD}}_{\sigma\pm1/2}=   \dfrac{\bbg_{\sigma}+\bbg_{\sigma\pm1}}{2} \pm \int_{K_{\sigma\pm1/2}} \big(\omega_{\sigma}-\dfrac{1}{2}\big)
\left(\dpar{\bbf^h}{x} -S^h \right)\,dx\;,
\end{equation}
In one space dimension, all these schemes are compatible with the discrete steady state
\begin{equation}\label{g-flux-consistent}
\dpar{\bbf^h}{x} =S^h  \iff  \quad \forall \sigma \quad\bbg_{\sigma}(\bu,S(\bu,x))=\bbg_{0}=\text{const}\
\end{equation}
which is here the only relevant consistency condition.

\bigskip
\noindent\emph{Second order at steady state.}  It is important to remark the following: the residual distribution numerical flux  \eqref{g-flux-RD}
is a compact consistent flux (in the sense of \eqref{g-flux-consistent}) which takes as inputs \emph{unreconstructed} states:
$$
\hbbg^{\text{RD}}_{\sigma +1/2}=\hbbg^{\text{RD}}_{\sigma +1/2}(\bu_{\sigma},\bu_{\sigma+1};x_{\sigma},x_{\sigma+1})
$$
Despite of this fact, the residual formulation provides a framework to design the flux in a way guaranteeing  at least second-order truncation at steady state,
without any gradient reconstruction.  This can be shown for  steady balance laws following e.g.  \cite{Dec:17} by estimating  the truncation
error defined as (see also \cite{RAD07}, appendix B)
\begin{equation*}
\begin{split}
\epsilon :=& \big \Vert\sum\limits_{\sigma} v(x_{\sigma}) \sum\limits_{K, \sigma\in K} \Phi_\sigma^{K}(\bw^h_{\text{ex}})\big \Vert \\
=& \big \Vert \int_{\Omega}v^h\left(  \dpar{\bbf^h_{\text{ex}}}{x}-S^h_{\text{ex}} \right)  + \sum\limits_K\sum\limits_{\sigma,\sigma'\in K}\dfrac{v(x_{\sigma})-v(x_{\sigma'})}{2}\int_{K}(\omega_{\sigma}-\varphi_\sigma)
\left(  \dpar{\bbf^h_{\text{ex}}}{x}-S^h_{\text{ex}} \right) 
\big \Vert
\end{split}
\end{equation*}
with $v(x)$ any smooth compactly supported test function,  with $\bw_{\text{ex}}$ a regular enough steady solution,
  $\Phi_\sigma^{K}(\bw^h_{\text{ex}})$ the residual distribution \eqref{phi-rd-source} evaluated when nodally replacing the numerical solution with samples
of the exact one. The analysis shows that the main design rules for second order  fluxes of the form  \eqref{g-flux-RD} are
  the boundedness of $\omega_{\sigma}$ and the formal second   order of the spatial approximations of the flux $\bbf^h$, and of the source $S^h$, which are readily obtained by means of e.g.
  linear interpolation  between  two neighbouring states.
  
The most  classical   particular case is  the upwind fluctuation splitting of Roe \cite{Roe87}, obtained  in the $P^1$ case by setting 
  $$
  \omega_\sigma\big|_{\sigma\pm 1/2} := \dfrac{1 \mp \text{sign}(\widetilde{\nabla_{\bu}\bbf}_{\sigma\pm 1/2}) }{2}
  $$
where  the sign of a matrix is defined as usual via its  eigen-decomposition, and where   following  \cite{Roe87} $\widetilde{\nabla_{\bu}\bbf}$ denotes the exact  linearization of the flux Jacobian verifying the conservation condition
$$
\widetilde{\nabla_{\bu}\bbf}_{\sigma\pm 1/2}\Delta \bu_{\sigma\pm 1/2} =\Delta \bbf_{\sigma\pm 1/2}\,.
$$
Note that in 1D this linearization establishes a direct link between the  cell conservation relation \eqref{conservation:K_c} and the linearized  non-conservative form of the PDE.
This allows to mention another known particular case, when the initial differential problem contains non-conservative terms
$$
\dpar{\bbf(\bu)}{x}+B(\bu)\dpar{\bu}{x}=S(\bu,x).
$$
In this case, one cannot simple apply the definition of conservation according to the principles  introduced so far. 
In the residual distribution setting this can be handled by embedding the non-conservative term in the   cell residual, so that \eqref{conservation:K_cS}
becomes in 1D
$$
\sum\limits_{\sigma \in K} \Phi_{\sigma}^K(\bu^h) = \Phi^K(\bu^h) = \int_K(\dpar{\bbf^h}{x}- S^h + B(\bu^h)\dpar{\bu^h}{x})d\bx
$$
The approximation of the last term has been reduced in the residual  distribution setting simply to  a quadrature 
problem for a given (linear) variation  of $\bu^h$ (see e.g. \cite{STAEDTKE2005379,Ricchiuto2003ARD} \S \ref{conservation}.5, and \cite{vrd00,Valero20091950}).
This is exactly what is done in   path-conservative  finite volume  (see \cite{CASTRO2017131} and references therein), when the path  chosen to connect  the left and right states at a cell interface is linear.
In particular, if $A={\nabla_{\bu}\bbf}+ B$, path conservative finite volumes are equivalent to the residual scheme obtained with
  $$
  \omega_\sigma\big|_{\sigma\pm 1/2} := \dfrac{1 \mp \text{sign}(A_{\sigma\pm 1/2})) }{2}\,,\quad
  \Phi^K=\Delta \bbf_K - |K| S_K +  B_K \Delta \bu_K
  $$
where$B_K$ can be evaluated for example with a one point quadrature  over the element. The  approximation and quadrature choices made  above to evaluate 
$$(B^h\dpar{\bu^h}{x})(x)= B(\bu^h(x))\dpar{\bu^h(x)}{x}$$ 
correspond to  the choice of the path in the finite volume context.  Assuming $\bu^h(x)$ to   linearly join two  states is one possibility.
One coud  also  have  $ \bu^h(x)= \bu(\bv^h(x), f(x) )$ with $\bv=\bv(\bu)$ some array of physical states (assumed to  by a $C^1$  invertible function of $\bu$ and to vary linearly), and $f$ a given field.
Many other choices are  possible. \remi{ Note that this does not solve   the issues raised by the non-conservative nature of the system, namely the
fact that the classical characterizaation of weak solutions and the  Lax Wendroff  theorem 
cannot we applied. This leaves all the uncertainties on the
 right form for a numerical scheme, see \cite{karni} for a counter example. 
However, we also remark that the RD framework can help in correcting schemes that discretise a non conservative form of a system in conservation form to account for certain constraints: see \cite{ABGRALL201810} for an example involving   multiphase flows. }

\subsection{Multiple dimensions, beyond second order, and other extensions}

The discussion provided allows to systematically design, by means of a residual based approach,
well balanced fluxes with  a genuine  second  order truncation without the need of any reconstruction.
We consider here several extensions, with focus on the multidimensional  steady case:
\begin{equation}
\nabla\cdot\bbf(\bu) +\bbB(\bu)\cdot\nabla\bu = S(\bu,\bx)\,,
\end{equation}
although we will not dwell too much on the issues related to the  the presence of the non-conservative term \remi{for the reasons stated above.}

The main recipe behind the method considered is already contained in
equation \eqref{phi-rd-source}. As in the 1D case, without loos of generality we will  assume that 
the discrete  unknowns are obtained as the steady limit of the pseudo-time iteration
\begin{equation}\label{update}
|C_{\sigma}|\dfrac{\bu_{\sigma}^{n+1}-\bu_{\sigma}^{n}}{\Delta t} + \sum\limits_{K, \sigma\in K}\Phi_{\sigma}^K(\bu^{h,n})
\end{equation}
with the conservation/consistency constraint that
\begin{equation}\label{conservation:K_cbS}
\begin{split}
\Phi^K:=\sum\limits_{\sigma \in K}\Phi_{\sigma}^K(\bu^{h}) =& \oint_{\partial K}\bbf^h_{\bn} d\gamma - \int_{K}S^h d\bx + \int_{K}  (\bbB\cdot\nabla\bu )^h d\bx \\
= &\int_K\left\{\nabla\cdot \bbf^h -S^h +  (\bbB\cdot\nabla\bu )^h
\right\}d\bx
\end{split}
\end{equation}
where as before $\bbf^h$ is the polynomial flux approximation of the highest degree for which the quadrature employed is exact, while both $S^h$ and 
$(\bbB\cdot\nabla\bu )^h$ are appropriately defined continuous approximations of the source and non-conservative terms, consistent with the quadrature strategy adopted. 

\bigskip
\noindent\emph{Global fluxes.}   Without specifying the form of $\Phi_{\sigma}^K$  we could repeat the construction of section \S \ref{flux} and  abstractly provide definitions of local fluxes embedding all the terms of the PDE.
Differently from the 1D case however, in multiple dimensions  the presence of the source term $S$, makes it quite unclear how to define consistency in a genuinely multidimensional setting. 

Concerning the non-conservative term, the choice of  the approximation/quadrature for the term $ (\bbB\cdot\nabla\bu )^h$ can be seen as choosing the manifold along which 
solutions can evolve. In this sense one could speak of  \emph{manifold-conservative} approach. As for path-conservative schemes, the authors remain skeptical 
as to how much  specifying this   notion would allow to side-step the fact that the classical definition of weak solution does not apply here. As in 1D, several choices are possible, the most obvious being
here to take 
$$
(\bbB\cdot\nabla\bu )^h=\bbB(\bu^h)\cdot\nabla\bu^h
$$
and evaluate the integral of this term be means of some quadrature formula. In the remainder of the paper, we will omit this term as none of the examples considered contain it.

\bigskip
\noindent\emph{Consistency for general smooth  steady solutions.}    The examples provided in section \S \ref{sec:conservation}can all be cast as a particular case of the general prototype
\begin{equation}\label{proto}
\begin{split}
\Phi_{\sigma}^K(\bu^h)=\!\!\int_{K}\!\omega_{\sigma}\left\{\nabla\cdot \bbf^h -S^h  
\right\}d\bx + \oint_{\partial K}    \,
[\![    \mathcal{L}(\varphi_{\sigma})]\!]\cdot [\![ \tau_{\mathcal{L}} \mathcal{L}(\bu^h) ]\!] d\gamma
\end{split}
\end{equation}
with $\mathcal{L}(\cdot)$ some linear differential operator.
The error analysis recalled in section \S \ref{sec: bien roule}.1 can be used in this more general setting.  In particular, given a smooth exact solution $\bw$ we define
\begin{equation}\label{err}
\begin{split}
\epsilon(\bw^h):=&\Big\Vert\sum\limits_{\sigma}\sum\limits_{K, \sigma\in K}v(\bx_{\sigma})\Phi_{\sigma}^K(\bw^h)
\Big\Vert=\Big\Vert\int_{\Omega}v^h\left\{\nabla\cdot \bbf^h -S^h  
\right\}d\bx\\ 
&\quad+
\sum\limits_K\sum\limits_{\sigma,\sigma'\in K}\dfrac{v(\bx_{\sigma})-v(\bx_{\sigma'})}{N_K}
\int_K (\omega_{\sigma}-\varphi_{\sigma} )\left\{\nabla\cdot \bbf^h -S^h  
\right\}d\bx \\  &\quad+\sum\limits_K\sum\limits_{\sigma,\sigma'\in K}\dfrac{v(\bx_{\sigma})-v(\bx_{\sigma'})}{N_K} \oint_{\partial K} \!\!   \,
[\![    \mathcal{L}(\varphi_{\sigma})]\!]\cdot [\![ \tau_{\mathcal{L}} \mathcal{L}(\bu^h) ]\!] d\gamma
\Big\Vert
\end{split}
\end{equation}
Simple approximation arguments can be used to show that \cite{AR:17} the above prototype has a consistency of order $\mathcal{O}(h^{p+1})$ 
as soon as the underlying polynomial approximation  is of degree $p$, and  provided that
$\omega_{\sigma}$ is uniformly bounded (w.r.t.  solution, mesh size, and problem data), and  that   the  $\tau_{\mathcal{L}}$ scales appropriately. 
For  $\mathcal{L}=\nabla$, the appropriate scaling is    $\tau_{\mathcal{L}}=\mathcal{O}(h^2)$  as in the Galerkin with jump stabilization  \eqref{burman}.
This  generalizes the compact second order construction discussed in the previous section to the multidimensional case, and to higher  degree approximations.
The estimate essentially allows to recover the underlying finite element approximation error.  One can however to more if some knowledge of the exact solution is embedded
in this approximation.

\bigskip
\noindent\emph{Super-consistency exact preservation of steady invariants.}   
Interesting results  can be shown when the source term  depends on some given data, say a given field $f(\bx)$ as 
for example the bathymetry in the    shallow water  equations, or some geometrical parametrization
 when considering the solution of the differential problem on a manifold (see e.g. \cite{rbl04} and references therein).  We are in particular interested in 
 exact steady solutions  characterized by the existence of a  set of invariants $\bv=\bv(\bbu,f)$ constant throughout the spatial domain.   Several examples
 have been provided in the introduction.
Assuming a sufficient smoothness of $f$, of the solution, and of the mapping  $(\bv,f)\mapsto \bbu(\bv,f)$,  we can write 
$$
\nabla\cdot \bbf(\bu) = \left(\nabla_{\bu}\bbf \, \nabla_{\bv}\bu \right)\cdot \nabla \bv +  \left(\nabla_{\bu}\bbf \,\nabla_{f}\bu\right)  \cdot \nabla f
=   \nabla_{\bv}\bbf\cdot \nabla \bv +\nabla_{f }\bbf \cdot \nabla f
$$
Solutions characterized by the invariance relation $\bv=\bv_0=$const. $\forall\bx$, satisfy  
\begin{equation}\label{eqinv}
\nabla_{f }\bbf(\bv_0,f) \cdot \nabla f +S(\bv_0,f) =0\,.
\end{equation}
This shows that for these solutions the flux and source dependence on the data,
and the    approximation and quadrature of the latter will play a crucial role.
For smooth/simple enough problems, the above relation   can be reproduced quite accurately in the residual context.
An interesting result  can be obtained by analyzing the error \eqref{err} \mario{when the approximation is written directly for the } 
steady invariants, and thus $\epsilon(\bw^h) =\epsilon(\bv^h,f^h) =\epsilon(\bv_0,f^h) $. \mario{For schemes 
of the form \eqref{proto}  with  approximation/quadrature choices consistent with  exactness for $\bv$ constant,  
the following is shown in \cite{R11,R15}.}
 
\begin{proposition}[Steady invariants and superconsistency]
Under standard regularity assumptions on the mesh, provided  the  test function $\omega_\sigma$  in \eqref{proto}
 is uniformly bounded w.r.t. $h$, $\bu_h$, element residuals, data of the problem,
 and provided   $\tau_{\mathcal{L}}$   is $\mathcal{O}(h^{2d_{\mathcal{L}}})$ with 
 $d_{\mathcal{L}}$ the highest derivative order of the operator $\mathcal{L}$, then:
 \begin{itemize}
\item  for exact integration   scheme \eqref{update}-\eqref{proto}  with $\bu^h = \bu(\bv^h,f)$,  $\bbf^h = \bbf^h(\bv^h,f)$,  and  $S^h = S^h(\bv^h,f)$
 preserves exactly the  equilibrium \eqref{eqinv};
\item  for approximate integration,  assuming that a  flux quadrature exact for approximate polynomial fluxes of  degree $p_f$ is used,
and a source quadrature exact for approximate polynomial sources of degree $p_v$, and assuming that $f\in H^{p+1}$ with 
$\nabla f \in H^p$, and $p> \min (p_f, p_v)$, then the scheme  is superconsistent w.r.t.
 \eqref{eqinv}, and in particular, its consistency is of order $r=\min(p_f+2,p_v+3)$.
\end{itemize}
\end{proposition}
Independently of the  details, the meaning of this result is that if the field $f$ \mario{and its derivatives
can be approximated by a smooth enough function given analytically,}
 and if the approximation is done in terms
of steady invariants instead of conserved variables, than the consistency of the  scheme is determined by the quadrature strategy and it is in particular
independent on the order of the underlying approximation 

\bigskip
A few remarks are in order.   The numerical results will provide examples indicating that  numerical convergence  w.r.t. the order of the quadrature formulas is indeed observed in practice,
at least for simple cases. However, exact preservation is possible for some important  and physically relevant examples.  We can mention at least two for the shallow water equations :
\begin{enumerate}
\item Lake at rest state \eqref{lar}. Exact preservation has been guaranteed by choosing the same approximation   for  $h$ and $b$,
and performing the quadrature of the   hydrostatic  terms  $\omega_{sigma}( \nabla (gh^2/2)+ gh\nabla b)$ either exactly, or using  the chain rule 
  $\omega_{sigma}( \nabla (gh^2/2)+ gh\nabla b)=\omega_{sigma} gh\nabla \eta $ \cite{RAD07,RB09};
\item Constant slope equilibrium \eqref{csf}. This is a case compatible with   contant $\bu$, for which any consistent  quadrature becomes exact. However, exact preservation 
is  guaranteed only due to the residual formulation in which  the different source terms are simultaneously integrated \cite{R15}. 
\end{enumerate}
Another remark concerns the   smoothness of $f$. The proposition above is built upon estimates involving approximation estimates, and related quadrature error formulas over elements.
This suggests that one can construct higher order approximations with less quadrature points either by  means of a clever    mesh generation step, 
embedding regions contaning  jumps in $f$ of in its derivatives as mesh edges/points, or by means of an adaptive quadrature strategy, avoiding the use of quadrature formulas across such discontinuities.


\section{Time dependent problems}\label{sec:unsteady}
\subsection{Preliminaries: global fluxes, time derivative and mass matrices}
To fix some basic concepts, we start by the simplest problem: the 1D advection equation
$$
\dpar{u}{t}+a\dpar{u}{x} =0
$$
The most classical discretization we can apply is the upwind scheme  
$$
\mario{\Delta_{\sigma} x \,} \dfrac{du_\sigma}{dt} = - (\hf_{\sigma+1/2}- \hf_{\sigma-1/2}) \;,\;\; 
\hf_{\sigma+1/2}=\left\{\begin{array}{ll}
f_\sigma\quad&\text{if }a>0\\
f_{\sigma+1}\quad&\text{if }a<0\\
\end{array}\right..
$$
This scheme is known to be only first order in space, and typically high order approximations are obtained by replacing the values of $f_\sigma=f(u_\sigma)$ and $u_\sigma$ by
evaluations of appropriately reconstructed polynomials on either side of  the interfaces $\sigma \pm 1/2$. \\ 

Let us make another experiment instead: we  set $S:= -\partial u/\partial t$, and apply an  upwind global flux method. We can proceed as in section \S5.1,
and formally define a global flux consistent with 
$$
g=f -\int_{x_0}^x S 
$$
We can now for example define upwind fluxes as 
$$
\hg_{\sigma+1/2}=\left\{\begin{array}{ll}
g_\sigma\quad&\text{if }a>0\\
g_{\sigma+1}\quad&\text{if }a<0\\
\end{array}\right.,
$$
having set 
$$
g_{\sigma+1} = g_{\sigma}+ \Delta f_{\sigma+1/2} - \int_{x_\sigma}^{x_{\sigma+1}} S \approx g_{\sigma}+ \Delta f_{\sigma+1/2} + \dfrac{|K_{\sigma+1/2}|}{2}\big(\dpar{u_\sigma}{t} +\dpar{u_{\sigma+1}}{t}\big)\;.
$$
The resulting scheme reads
$$
 \remi{\hg_{\sigma+1/2}-\hg_{\sigma-1/2} = 0\,,}
$$
or equivalently, using the definition of the numerical flux, and rearranging terms: 
$$
\dfrac{g_{\sigma+1} -g_{\sigma-1} }{2} +\dfrac{\text{sgn}(a)}{2}( g_\sigma-g_{\sigma-1} )
- \dfrac{\text{sgn}(a)}{2}(g_{\sigma+1}  - g_\sigma) =0
$$
The definition of the global flux given above leads, after some manipulations,  to the  following semi-discrete evolution scheme
\begin{equation}\label{uns0}
\mario{\Delta_{\sigma} x \,} \dfrac{d\hat{u}_{\sigma}}{dt} +\hat{f}_{\sigma+1/2}- \hat{f}_{\sigma-1/2} =0
\end{equation}
where  $\hat{f}_{\sigma\pm 1/2}$ are exactly those of  the first order upwind scheme, while the nodal approximation of the time derivative is now defined as 
\begin{equation}\label{uns1}
\begin{split}
\Delta_\sigma x \dfrac{d\hat{u}_\sigma}{dt}:=&
\dfrac{|K_{\sigma-1/2}|}{4}\big(\dpar{u_{\sigma-1}}{t} +\dpar{u_{\sigma}}{t}\big)+ \dfrac{|K_{\sigma+1/2}|}{4}\big(\dpar{u_{\sigma+1}}{t}+\dpar{u_\sigma}{t} \big)\\
&+\dfrac{\text{sgn}(a)}{2} \dfrac{|K_{\sigma+1/2}|}{2}\big(\dpar{u_{\sigma+1}}{t} +\dpar{u_{\sigma}}{t}\big)  \\&\qquad-\dfrac{\text{sgn}(a) }{2} \dfrac{|K_{\sigma-1/2}|}{2}\big(\dpar{u_{\sigma-1}}{t} +\dpar{u_{\sigma}}{t}\big)
\end{split}
\end{equation}
Quite interestingly, 
this method can be checked (e.g. with a truncated Taylor series analysis)  to have a second order truncation error in space without the need of any polynomial reconstruction. \mario{This is not related to error compensation on a uniform mesh, but to the improved 
balance of the different terms for linear data within each cell.}
This simple example  shows how the notion of a global flux can be applied to other types of terms in the PDE. 
In the case of the time derivative, the global flux approach leads to the appearance of a mass matrix.

As we have shown previously, there is a direct between the upwind  finite volume method and residual based schemes which can be summarized into the equality
$$
\remi{\hf_{\sigma+1/2}- \hat{f}_{\sigma-1/2}=\sum\limits_{K, \sigma\in K}\int\limits_{K}\omega_\sigma\dpar{f^h}{x}\; dx}
$$
with $f^h$ piecewise linear. There are at least two definitions of the test function $\omega_{\sigma}$ which give back the upwind scheme,
namely
\begin{equation*}
\begin{split}
{\omega_\sigma}_{{|_{K_{\sigma\pm 1/2}}}}=&\varphi_\sigma + a\dpar{\varphi_\sigma}{x}\,\tau_{\sigma+1/2}\,,\;\;\tau_{\sigma+1/2}=\dfrac{|K_{\sigma+1/2}|}{2|a|}\\
\text{ and }&\\
{\omega_\sigma}_{|_{K_{\sigma\pm1/2}}}=& \dfrac{1\mp\text{sgn}(a)}{2}
\end{split}
\end{equation*}
with $\varphi_{\sigma}$ the linear finite element test functions.
The method \eqref{uns0}-\eqref{uns1} can be obtained in a much more natural and elegant way as a particular case of a residual method,
and in particular of the one corresponding to the second definition  of $\omega_\sigma$ above. The first definition, provides and even better variant with
a truncation error \mario{which improves to an} order $\Delta x^3$ (see e.g. \cite{rf14})  \remi{for uniform meshes} !  The benefit of this idea is to allow high order of accuracy with the most  compact   stencil.
Its drawback is that it requires inverting the mass matrix.
This analogy has also been used in other context to generate compact  high order 
finite difference schemes associated to a variational form \cite{LXZ18}.\\

The next sections discuss how to generalize this idea to multiple dimensions and, more importantly, how the issue of inverting the mass matrix has been side-stepped.

\subsection{Generalization}
As the last section has shown, following a finite element strategy for \eqref{eq:1} for the unsteady case  will always lead to a formulation of the form
$$M\big ( \bu^{n+1}-\bu^n\big ) +\Delta t\; \delta \mathbf{F}=\mathcal{S},$$
where $M$ is a mass matrix, $\delta \mathbf{F}$ contains all the spatial approximation terms, and $\mathcal{S}$ the approximation of the source term. It is possible, depending on the formulation, that several instances of $\bu$ appear. One of the biggest problem is  the mass matrix.

Things are different for  the classical formulations of finite volume and discontinous Galerkin schemes, which lead to diagonal or block diagonal matrices because of the locality of the approximation of $\bu$. These \mario{are small (however dense) matrices  which can be inverted locally on each element, and are independent on the mesh connectivity.}
This is probably one of the keys of the success of these methods, especially for genuinely   hyperbolic  and evolutionary problems.  In the case of continuous approximation, the story is not as simple. For example, the SUPG method will lead to a mass matrix that may evolve in time. This is also the case of the RD schemes developed in \cite{AbgralldeSantisNS,abgrall99,Mario,Mario1,Mario2}. The Galerkin method with jump stabilisation does not have this problem, but nevertheless, we still have a sparse positive definite matrix to invert. One of the strategies followed in the past has been to work on highly implicit variants of the schemes, trying to cover this computational overhead with
the possibility of using large time steps.  Unfortunately, despite the excellent results, the schemes obtained in this way are relatively cumbersome to code \cite{abg2001c,drd05,Mario1,hr11}.
Moreover, the advantage of using large time steps,  very useful for viscous flows and  problems with large stiffness, is 
less obvious for   wave propagation problems, even on non-uniform meshes \cite{hr11,hrs00,hrs20}.

In \cite{Mario} is explained how to approximate the solution in time, but without having to invert a mass matrix. In this reference, the method is explained for piecewise linear element (and triangular element). The method was further extended to any order (and any type of simplex) in \cite{abgrallDec}. 

In practice, for the steady version of \eqref{eq:1}, each of the known schemes can be written using test functions. This is clear for the SUPG scheme, where the test functions are defined in each element and are possibly discontinuous accross element. Please note that in the non linear case, the test functions will depend on $\bu$. The same is true for the schemes of \cite{abgrallLarat,abgralldeSantisSISC}, except that the scheme will be non linear even for a linear problem in order to enforce non oscillatory constraints. In the case of Galerkin method with jump, one can also reinterpret the method in this way, thanks to the use of  a lifting operator allowing to embed the jump terms in a numerical flux. Hence, in all cases, we write
$$\Phi_\sigma^K(\bu^h)=\int_K \psiT_\sigma \; \text{ div }\bbf(\bu^h) \; d\bx$$
where $\psiT_\sigma$ is the test function associated to the element $K$ and the degree of freedom $\sigma$. For example, for the SUPG method, this is 
$$
\psiT_\sigma=\varphi_\sigma+ h_K \nabla_{\bu}\bbf(\bu^h) \cdot\nabla \varphi_\sigma\; \tau_K.$$
Integrating \eqref{eq:1}, and using, for simplicity of exposure, the mid point rule in time,  will lead to
\begin{equation}\label{eq:time}
\int_\Omega \psiT_\sigma \big ( \bu^{n+1}-\bu^n\big ) +\frac{\Delta t}{2} \bigg ( \int_\Omega \psiT_\sigma\text{ div }\bbf(\bu^{n+1})\; d\bx+\int_\Omega \psiT_\sigma\text{ div }\bbf(\bu^{n})\; d\bx\bigg )=0,
\end{equation}
that despite its complexity, we still can rewrite in a form similar to \eqref{conservation} with:
$$\Phi_\sigma^K(\bu, \bv)=\int_K \psiT_\sigma \big ( \bu-\bv\big ) +\frac{\Delta t}{2} \bigg ( \int_K \psiT_\sigma\text{ div }\bbf(\bu)\; d\bx+\int_K \psiT_\sigma\text{ div }\bbf(\bv)\; d\bx\bigg )$$
and $\bu=\bu^{n+1}$, $\bv=\bu^n$.

What makes this complex is the term in time. If we consider the simpler set of residual, \remi{ for $|C_\sigma|>0$ to be defined, }
$$\psi_\sigma^K(\bu, \bv)=|C_\sigma^K| \big ( \bu_\sigma-\bv_\sigma)+\frac{\Delta t}{2} \bigg ( \int_K \psiT_\sigma\text{ div }\bbf(\bu)\; d\bx+\int_K \psiT_\sigma\text{ div }\bbf(\bv)\; d\bx\bigg ),$$
then the scheme defined from this is easily solvable: we know $\bu^n$, we can get $\bu^{n+1}$ explicitely. The key question are (i) how  to define the lumping parameter $|C_\sigma^K|$, and (ii) how to combine the schemes defined by this two set of residual in order to get an approximation the solution $\bu^{n+1}$ given by \eqref{eq:time} with the \emph{same} accuracy.

In \cite{Mario}, in the case of a ${\mathbb P}^1$ approximation and triangle element, it is shown that if $|C_\sigma^K|=\tfrac{|K|}{3}$ and if we define $\bu^{n+1}$ \remi{using a predictor corrector algorithm} as:
\begin{equation}\label{eq:unsteady:O2}
\begin{split}
|C_\sigma|\big (\bu_\sigma^{(1)}-\bu_\sigma^n\big )=&- \Delta t\int_\Omega \psiT_\sigma\text{ div }\bbf(\bu^n)\;d\bx\\
|C_\sigma|\big (\bu_\sigma^{(2)}-\bu_\sigma^{(1)}\big )=&- \frac{\Delta t}{2}
\bigg ( \int_\Omega \psiT_\sigma\text{ div }\bbf(\bu^{(1)})\; d\bx+\int_\omega \psiT_\sigma\text{ div }\bbf(\bu^n)\; d\bx\bigg ) \\
&- \mario{\int_\Omega \psiT_\sigma \big ( \bu^{(1)}-\bu^n\big )} \\
\bu^{n+1}=\bu^{(2)}&
\end{split}
\end{equation}
then we have a second order scheme in time, with \mario{similar stability} 
properties as the original steady scheme.

The extension to higher than second order and general simplex has been done in \cite{abgrallDec}. The main idea is to notice that \eqref{eq:unsteady:O2} can be reinterpreted as a defect correction method: if one wants to solve $L^{(2)}(U)=0$, and if one has a second operator, called $L^{(1)}$, such that in some norm,
\begin{subequations}\label{dec:conditions}
\begin{equation}\label{dec:1}
\Vert \big (L^{(1)}(U)-L^{(2)}(U)\big )- \big (L^{(1)}(V)-L^{(2)}(V)\big )\Vert \leq \Delta \Vert U-V\Vert,\end{equation}
if in addition $L^{(1)}$ satisfies a coercivity relation,
\begin{equation}\label{dec:2}
\alpha \Vert U-V\Vert \leq \Vert L^{(1)}(U)-L^{(1)}(V)\Vert, 
\end{equation}
\end{subequations}and finally, if $L^{(2)}(U)=0$ has a unique solution $U^\star$, then the solution $U^{(p)}$ of the iterative scheme:
\begin{equation}\label{eq:Dec}
\begin{split}
U^{(0)} & \text{ given }\\
\text{ do for }p\geq 0\quad L^{(1)}(U^{(p+1)})&=L^{(1)}(U^{(p)})-L^{(2)}(U^{(p)}).
\end{split}
\end{equation}
\remi{Of course the question is to know a good stopping criteria. The answer is given by the following:}
it can be shown \cite{abgrallDec} that
\begin{equation}
\label{dec:error}\Vert U^{(p)}-U^\star\Vert\leq \bigg ( \frac{\Delta}{\alpha}\bigg )^p\Vert U^{(0)}-U^\star\Vert.\end{equation}

Here, if
 the coefficients $|C_\sigma^K|$ are chosen such that
 $$\sum_{\sigma \in K} |C_\sigma^K| \bu_\sigma=\int_K \bu(x)\; d\bx,$$
 i.e
 $$ |C_\sigma^K|=\int_K \varphi_\sigma\; d\bx,$$
\remi{ then the conditions \eqref{dec:conditions} are met and }
 then a CFL-like condition $\tfrac{\Delta}{\alpha}\approx \Delta t$,  we can interpret \eqref{dec:error}, after $p$ iterations,  as
 $$\Vert U^{(p)}-U^\star\Vert\approx C \Delta t^{p}.$$
 This means that if the un-lumped formulation is of order $p$, then the lumped with the algorithm \eqref{eq:Dec} one will provide a solution with the same accuracy after $p$ iteration only. In the case of piecewise linear elements, $L^2$ is defined from the residuals $\Phi_\sigma^K$ and $L^1$ is defined from $\Psi_\sigma^K$.
 
 The problem is that often 
 $\int_K \varphi_\sigma\; d\bx$ is not positive: this is the case for quadratic Lagrange interpolant in triangles. A possible remedy to  this is to use basis functions that are positive as for example   B\'ezier polynomials  \cite{abgrallDec}.
 Note that the linear polynomials are also B\'ezier polynomials of degree 1. 
 
 {\color{black}
In practice, we split the time interval $[t_n,t_{n+1}]$ with $p$ sub-time steps $t_n=t_{p,0}<t_{p,1}=t_n+\alpha_1\Delta t<\ldots <t_{p,p-1}=t_n+\alpha_{p-1}\Delta t<t_{p,p}=t_{n+1}=t_n+\Delta t$, the vector $\bu$ contains the approximations of $\bu$ for the sub-time steps, i.e. $\bu=(\bu_0, \bu_1, \ldots, \bu_p)$ with $\bu_j\approx \bu(~\cdot~,t_{p,j})$. Then we write
\begin{enumerate}
\item Set $\bu^{(0)}=(\bu^n, \bu^n, \ldots , \bu^n)$: we initialise the vector with the state at time $t_n$,
\item Do for $l=1, \ldots , p$, do for $k=1, \ldots p$
$$
|C_\sigma| \big ( \bu_k^{(l+1)}-\bu_k^{(l)}\big ) +\sum\limits_{K, \sigma \in K}\Big (\underbrace{ \int_K \psiT_\sigma \big ( \bu_k^{(l)}-\bu^n\big ) \; d\bx}_{(I)}+
\underbrace{\oint_{t_n}^{t_{p,k}} \bigg (\int_K\psiT_\sigma \text{ div }\bbf(\bu^{(p)}) \; d\bx \bigg )\Big )}_{(II)}=0
$$
$(I)$ is evaluated by spatial quadratures and $(II)$ by space-time quadratures:
we write
$$\oint_{t_n}^{t_{p,k}} \bigg (\int_K\psiT_\sigma \text{ div }\bbf(\bu) \; d\bx \bigg )=\Delta t\sum\limits_{r=0}^p\theta_r^k \int_K \psiT_\sigma \text{ div }\bbf(\bu_r^{(l)}) \; d\bx$$
where $\theta_r^k$ is the integral over $[0,\alpha_k]$ of the (time) Lagrange interpolant at the points $\{\alpha_0=0, \alpha_1, \ldots , \alpha_{p-1}, \alpha_p=1\}$.
and $|C_\sigma|=\sum\limits_{K, \sigma\in K}|C_\sigma^K|$.
\item $\bu^{n+1}=\bu^{(p+1)}$
\end{enumerate}
The integer $p$ is equal to the expected order of accuracy. The procedure is explicit.

Using the results of section \ref{flux}, it is easy to see that the sum of  $(I)$  and $(II)$ is a sum of flux: If the quadrature formula in time needs $p$ steps, then we extrude the element $K$ adding $p$ layers, this induces natural graph that is connected, and hence the discussion of section \ref{flux} can be repeated.
}
\subsection{Unsteady problems and well balanced on dynamic meshes}   

To complete the presentation of the time dependent case, we go back to  balance laws.   The schemes developed in these pages have been \mario{initially designed having in mind} 
general adaptive meshes.  In the time dependent case, it is thus natural to envision some dynamic adaptation method. 
The literature is  filled with promising adaptive techniques, see e.g. \cite{LOSEILLE2017263,DONAT2014937,ArR:18} and references therein for a (non-comprehensive) review.  
We want underline here an important   aspect allowing to generalize some of the concepts introduced in section \S5. In particular, when dynamic meshes are employed 
a fundamental step is the operator allowing to map the solution from one mesh to another. In steady computations the error possibly introduced during the projection
from one mesh to another may be lost at convergence, provided the boundary conditions are not affected by this aspect. In time dependent simulations however, the
remap may pollute the evolution of the solution, and  affect both its accuracy and   stability, as much as the underlying discretization method. 

Design criteria for the remap are \mario{thus} consistency,  conservation, monotonicity preservation,  etc: \mario{exactly the same criteria applied when
solving the main PDE problem !}  For balance laws, if the source term depends on external data, these must also be projected,
and a well balanced condition may be  on the shopping list of desired properties. To fix ideas, we will consider here projection methods based on some kind of Lagrangian,
or rather Arbitrary  Lagrangian Eulerian (ALE) remap  (see e.g.  \cite{re:hal-01633476,ArR:18}   and references therein).
To start with, we recast \eqref{eq:1}  with $S=S(\bu, f(\bx))$ in ALE form: 
\begin{equation}
\label{eq:1-ALE}
\dpar{ (J \bu)}{t}\Big|_{\bX} + J \text{ div }\big(\bbf(\bu) - \bbw \bu \big)=J\, S( \bu,f(\bx))
\end{equation}
with the additional definitions/relations
\begin{equation}
\label{eq:1-ALEa}
\begin{split}
\dfrac{d\bx(t)}{dt}\Big|_{\bX} &= \bbw\,,\;\; \bx(0) = \bX\\[10pt]
\dpar{ J }{t}\Big|_{\bX} &- J \text{ div }  \bbw  =0\\[10pt]
\dpar{ f }{t}\Big|_{\bX} &- \bbw\cdot \nabla  f  =0
\end{split}
\end{equation}
with $J$ the determinant of the Jacobian of the mapping $M: \bX\mapsto\bx$, namely
$$
J:=\text{det}\Big( \dpar{\bx}{\bX}\Big)\,.
$$ 
Note that the  numerical discretization essentially provides a discrete equivalent of \eqref{eq:1-ALE}. The relations \eqref{eq:1-ALEa},
which are true and exact  on the continuous level, 
are not explicitly solved numerically but must be  seen as contraints to embed as much as possible in the discretization. 
As in \mario{ \cite{re:hal-01633476,ArR:18}} 
we assume that mesh operations can be represented by some continuous deformation operator, so that the projection
from one mesh to the other boils down \mario{somehow to mimic    \eqref{eq:1-ALE}, with constraints \eqref{eq:1-ALEa}}. 
This allows to \mario{update the list of design criteria for the schemes:} 
\begin{itemize}
\item  \emph{ Discrete Geometric Conservation.} This is essentially a discrete analog of the second relation in \eqref{eq:1-ALEa} and represents  the conservation of volume 
along the mapping $M$;
\item  \emph{ Mass  Conservation.}  Without loss of generality, we can assume that     the first relation  in  \eqref{eq:1-ALE} 	 is homogenous (so $S_1=0$) and represents   mass conservation:
$$
\dpar{ (J\rho)}{t}\Big|_{\bX}+J\nabla\cdot(\bbv\rho -\bbw\rho) =0\,.
$$ 
Note that one of most   classical  translations of geometric conservation is to check that uniform flows are preserved   by the discretization \cite{ThL:79}. 
This corresponds to the fact that for $\rho$ and $\bbv$ constant the last mass conservation equation becomes  the second in \eqref{eq:1-ALEa}, and is also equivalent
to the standard notion of consistency with respect to constant states; 
\item  \emph{ Well \remi{balanceness}} As already remarked in section \S5, consistency for constant states is not necessarily applicable to balance laws, as 
only non-constant, data dependent, steady states are admissible. So well balanced is in contradiction with some of the above properties, and notably conservation;
\item  \emph{ ALE remap.} The last relation in \eqref{eq:1-ALEa} represents the ALE time derivative of the data $f$, which is also often written  in conservative form by combining it with geometric conservation:
$$\dpar{ (J f) }{t}\Big|_{\bX} - J \text{ div } ( \bbw  f)  =0\,.
$$
\mario{The ALE time  derivative is essentially an advection operator. Its approximation poses very similar questions of consistency, accuracy,
stability, and bounded variations for the data, as the approximation of  \eqref{eq:1-ALE}. Moreover, }
for problems in which the data play a key role (e.g. topography inundation and coastal risk assessment),
\mario{the deterioration of such data due to the approximation error may introduce an unacceptable uncertainty in the 
predictions. So the use of standard techniques to approximate this advection problem, as e.g. proposed in  \cite{Zho:13} to guarantee
both well balancedness and mass conservation,  
may be very delicate.  For example, to avoid diverging
from reality the last reference  proposes to periodically  re-initialize  the data,  which implies 
loosing the consistency with the ALE projection which may cost well balancedness, or mass conservation.} 
%
\end{itemize}
Ideally, all of the above properties should be satisfied. However,  as remarked   
 well balanced is in general in contradiction with e.g. mass conservation, \mario{and the  ALE projection may be at odds with the preservation of the 
 accuracy of the data involved. }
We consider here a (physically relevant) example to better highlight this issue, and show a possible solution in the context of  residual distribution.\\

{\bf  Example: Lake at rest solutions on moving meshes.} 
We consider the shallow water equations  in ALE form with an arbitrary (non constant)  bathymetry $b(\bx)$.
As discussed above,  the  classical characterization of the DGCL based on the preservation of 
the state $\bu=\bu_0=\text{const}$ cannot be used here,  as constant states are not  solutions of the problem due to the presence of the source.
To better fix ideas, we will recast  \eqref{eq:1-ALEa} as follows
 \begin{equation*}
 J \underbrace{ \left( \left.\frac{\partial \bu}{\partial t}\right|_{X}-  \bbw\cdot\nabla \bu \right) }_{H_1} 
 +\bu \underbrace{ { \left( \left.\frac{\partial J}{\partial t}\right|_{X}-J \nabla\cdot\boldsymbol \bbw\right) }}_{H_2}+
 J\underbrace{ {\biggl( \nabla\cdot\bbf - S \biggr) }}_{H_3}=0             
 \end{equation*}
As amply discussed in section \S5,  we have a general framework to devise  well-balanced Eulerian discretization methods do embed integral (or even local) versions of $H_3=0$. 
Previous work has shown how to extend this framework to an ALE setting for both explicit and implicit time integration
\cite{ArR:14,doi:10.2514/6.2005-493,hrs20,Mich:2001},  embedding  discretely the   constraint of the geometric conservation law $H_2=0$. 
Unfortunately  by their nature  Eulerian methods are unable to embed the condition $H_1=0$,  which will be polluted, also in correspondence of steady exact solutions which would 
be exactly represented on fixed meshes.

A possible way out of this limitation for solutions admitting a set of steady invariants $\bv$,  is that the ALE formulation, and thus the mesh projection,
should be performed using $\bv$ as main variable. To  explain we consider the lake at rest state, but other cases can be treated in a similar way.
In this case, we can set $\bv=[\eta, h\bbv]=[H+b, H\bbv]$, and steady states are characterized by $\bv=\bv_0=[\eta_0,0]$, and thus $h=h(\bx)=\eta_0-b(\bx)$.
In the continuous case  we can invoke  the fact that the bathymetry satisfies the ALE remap (last in \eqref{eq:1-ALEa}), namely 
$$
\underbrace{\left.\frac{\partial b}{\partial t}\right|_{X} - \bbw\cdot \nabla b=0 }_{H_4=0}
$$
This can be used to modify the  ALE formulation  and write it  directly in terms of $\bv$: 
\begin{equation}\label{eq:1-ALEb}
\left.\frac{\partial (J  \bv )}{\partial t}\right|_{X}
+J\nabla\cdot\big(\boldsymbol{f}(\bv, b) -  \bbw  \bv\big) + J  S =0
\end{equation}
This formulation is equivalent to
 \begin{equation*}
 J \underbrace{ \left( \left.\frac{\partial \bv}{\partial t}\right|_{X}- \bbw \cdot\nabla\bv \right) }_{H_1+H_4}
 +\bv\underbrace{ { \left( \left.\frac{\partial J}{\partial t}\right|_{X}- J \nabla\cdot\bbw  \right) }}_{H_2}+
 J \underbrace{ {\biggl( \nabla\cdot\boldsymbol{f} + S \biggr) }}_{H_3}=0             
 \end{equation*}
 We can now use any  Eulerian scheme  which is well balanced and compatible with geometric conservation,
 and we will be able to ensure that  all the terms in the above sumamtion will be zero if $\bv$ is constant.

For completeness, we recall that the formulation  \eqref{eq:1-ALEb}, which is referred to in  \cite{ArR:18} as to the well balanced form of the equations,
 is similar to   the pre-balanced form of the Shallow Water equations of  \cite{RBT:03} which  uses a modified definition of the 
  flux and source  terms (cf. \cite{ArR:18}  for datails). 

{\it The problem of mass conservation.} We now consider the additional constraint of achieving discrete   conservation of the total water mass in the domain.
We integrate in space and in time the mass conservation equation  in  well-balanced form  \eqref{eq:1-ALEb} 
\begin{equation*}
\int_{\Omega(t)} \eta(\boldsymbol{x}(t),t) \, d \boldsymbol{x} - \int_{\Omega_{X}} \eta(\boldsymbol{X},0) \, d \boldsymbol{x} + \int_{0}^t \!\! \int_{\partial\Omega(t)} ( H\bbv - \eta \bbw  )\cdot\boldsymbol{n}\, ds \, dt  = 0
\end{equation*}
Let $V(t)=\int_{\Omega(t)} H \, d \boldsymbol{x}$ be the total volume of water in the domain at time $t$, and define $B(t)=\int_{\Omega(t)} b \, d \boldsymbol{x}$. 
We can rewrite the above   conservation statement as
\begin{equation} 
H(t) - H(0)
+ \int_{0}^t \!\!\! \int_{\partial\Omega(t)} \! H( \bbv-\bbw ) \cdot\boldsymbol{n}\, ds dt  = - \Big(  B(t) - B(0) -\int_{0}^t \!\!\!\int_{\partial\Omega(t)} \! b\bbw  \cdot\boldsymbol{n}\, ds \, dt\Big)  \label{eq:masscons}
\end{equation}
which states that, modulo the boundary conditions, we have conservation over the full domain  provided that we satisfy geometry consrvation and   the bathymetry satisfies the ALE  remap
given by the last in  \eqref{eq:1-ALEa},  namely if
\begin{equation*}
 B(t) - B(0) - \int_{0}^t \!\!\!  \int_{\partial\Omega(t)} b \bbw \cdot\boldsymbol{n}\, ds \, dt = 0 
\end{equation*}
As already said,  some work in literature propose indeed to evolve the bathymetry according to the ALE  remap  (see e.g.  \cite{Zho:13}). As discussed earlier,
the uncertainty on the topography associated to the error introduced by this approach   may not be acceptable in many applications (e.g. assessment of coastal risks).
Combining the ALE remap with some periodic reinitialization of the data would  end up breaking mass conservation unless a more clever fix is sought.
A possible one has been suggested in \cite{ArR:18,ArR:20}, and is recalled hereafter.

Assume for simplicity that the domain boundaries are not moving, or that $\bbw\cdot\boldsymbol{n}$ is verified. We can write the mass error at time $t$ as 
\begin{equation*}
E_{mass} = H(t) -H(0) +   \int_{0}^t \!\!\!  \int_{\partial\Omega} H \bbv \cdot\boldsymbol{n}\, ds \, dt = B(0) - B(t)
\end{equation*}
We now remark that the two quantities on the right hand side are in principle equal, as they are both approximations of the integral of $b(\boldsymbol{x})$ over the domain. If the domain boundaries are not moving, this quantity should remain constant in time. In practice however, these two integrals will be evaluated on a moving mesh. This means that, even if both the domain of  integration  and the data being integrated are constant, the quadrature points used will move, so the result will not be the same.  To be more precise,  the evaluation of  $B(t)$ will be expressed by a sum which will depend on the time update of the scheme.
For example, for the explicit approach discussed in section \S6.2, we will have
\begin{equation}
B(t) = \sum_{\sigma }   b_i(t)|C_\sigma(t)|\label{eq:masscons-3}
\end{equation}
with $b_i=b(\boldsymbol{x}_i(t))$.
The idea proposed in  \cite{ArR:18,ArR:20} is to replace $b(\boldsymbol{x}_i(t))$  in the last expression by some mapped value, 
allowing to minimize the overall mass error, and exploiting as much as possible the actual bathymetric data.
In particular, one way to achieve this is to set
\begin{equation}
\tilde b_\sigma:= \dfrac{1}{|C_\sigma|} \int_{C_i(t)} b(\boldsymbol{x}(t)) \, d \boldsymbol{x} \approx \sum^{N_{q}}_{f=1} \omega_q b(\boldsymbol{x}_q(t))
\end{equation}
where the right hand side defines a high order accurate quadrature formula of the real (initial/reference) data over the  current cell.
The increase in accuracy of the  quadrature within each moving cell, allows to compensate for the movement of the cells themselves.
In practice, for most problems quadrature formulas exact for degree 2 polynomials are enough to keep the mass error to machine zero levels.

Please refer to  \cite{ArR:18,ArR:20}  to for the extension to problems with dry areas and on curvilinear coordinates.

\section{Examples}

\subsection{Some examples of compressible flows simulations}
We present two cases: the first one is the well known DMR test case by Colella  and Woodward: it is the interaction of a Mach 10 shock wave in a quiescent media with wedge angle of $30^\circ$. The result is \remi{displayed} in figure \ref{remi:fig:1} for a cubic (B\'ezier approximation) and the quality of results is comparable to what can be found in the litterature for a simular resolution (estimated as $100^2$ for a Cartesian mesh).
\begin{figure}[h]
\begin{center}
\subfigure[]{\includegraphics[width=0.45\textwidth]{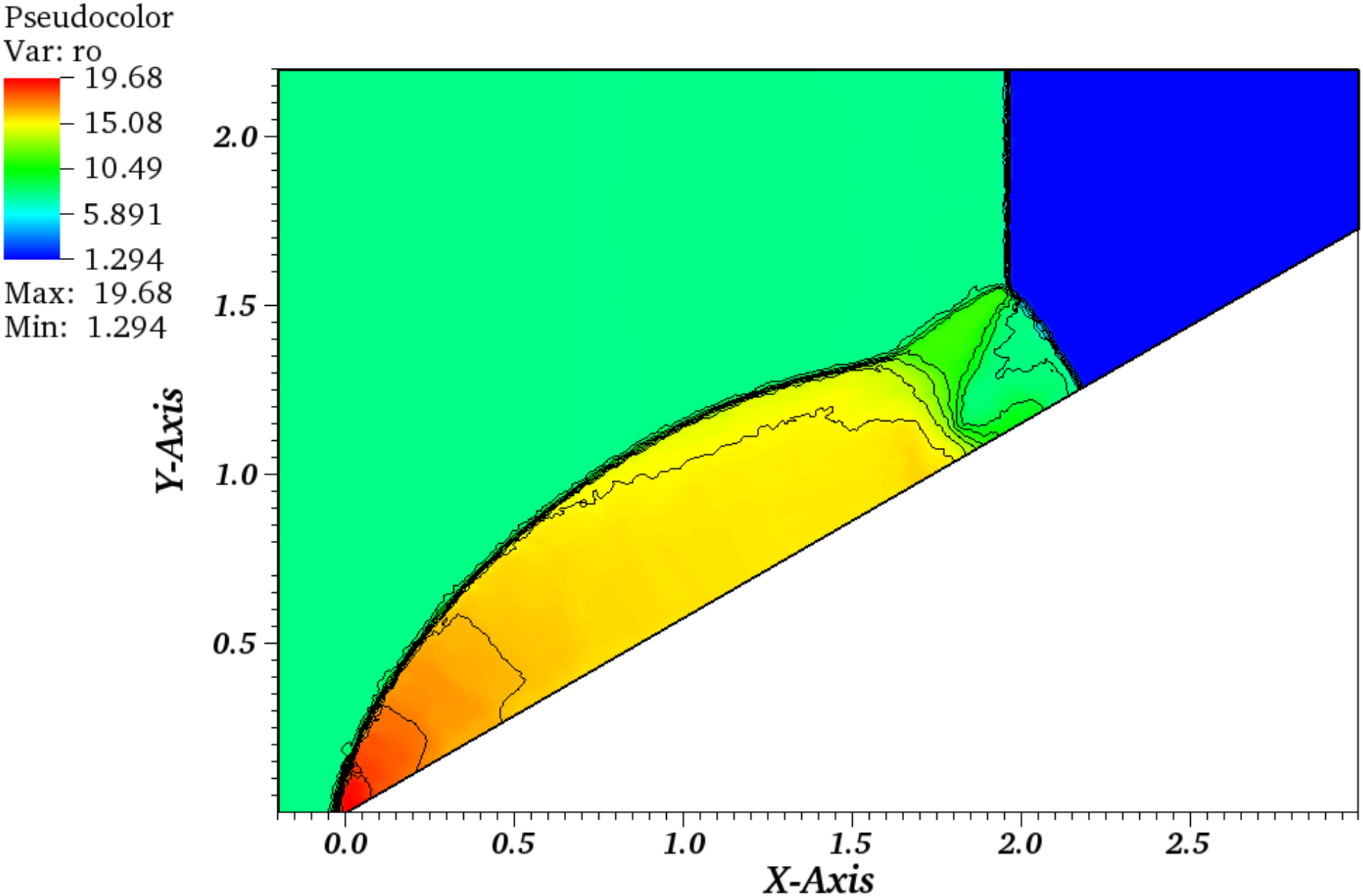}}\subfigure[]{\includegraphics[width=0.45\textwidth]{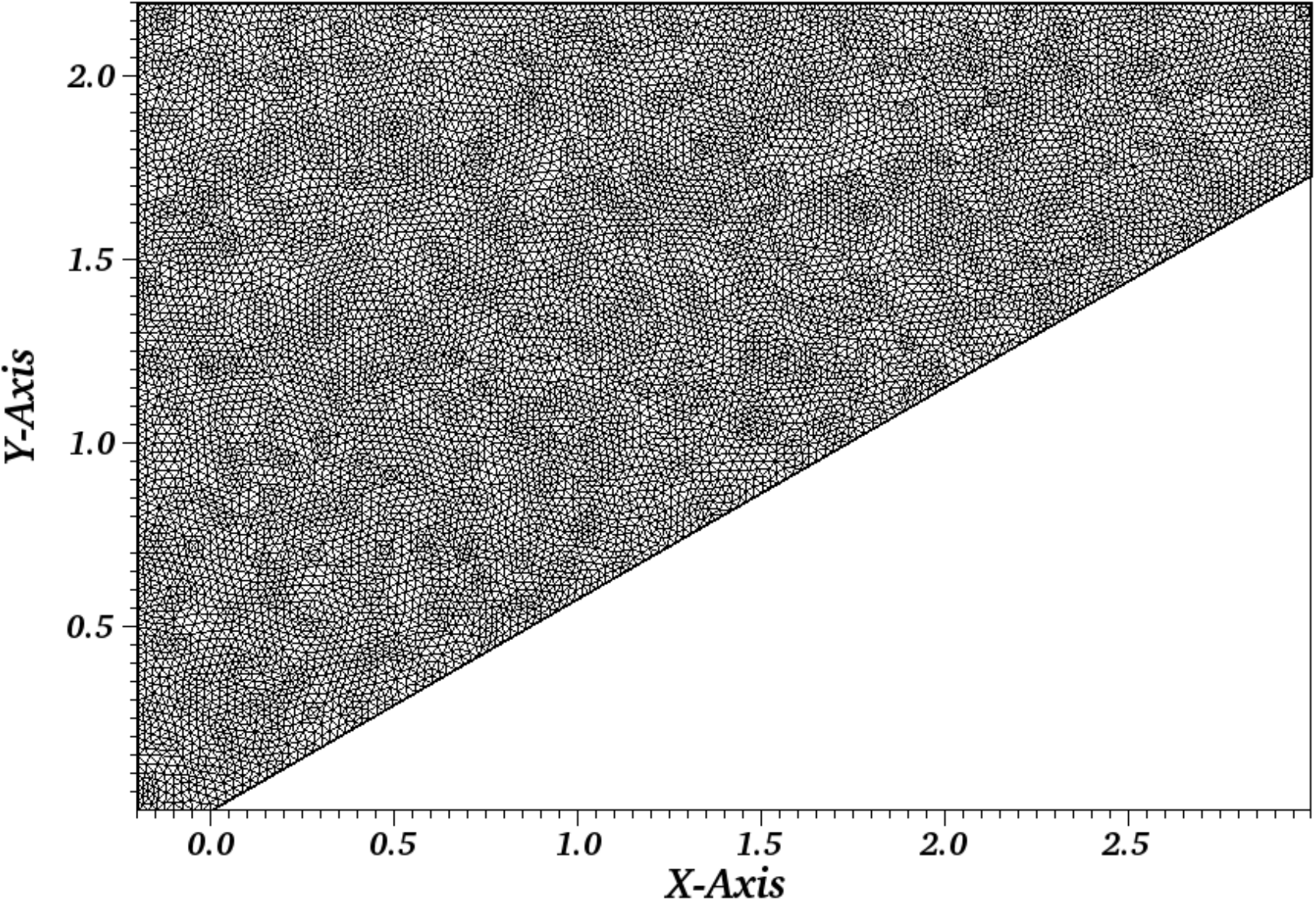}}
\end{center}
\caption{\label{remi:fig:1} DMR at $T=0.2$, solution on a triangular mesh with  19248 elements. }
\end{figure}
The second one can be seen as a 2D version of the Shu-Osher case where we have the interaction between a shock wave and a density wave. The conditions are thus: at $\bx$ in the disc of center $(0,0)$ and radius $6$,
$$
\begin{pmatrix}\rho \\ \vec{u} \\p\end{pmatrix}=\left \{
\begin{array}{ll}
(3.857143, 2.69369\frac{\bx}{\Vert \bx\Vert}, 10.333333)^T &\text{ if } \Vert \bx\Vert\leq 1\\
(1+2\sin(5\Vert x\Vert), \vec{0}, 1)^T & \text{ if } 1<\Vert \bx\Vert \leq 4\\
(1+2\sin(5\times 4), \vec{0}, 1)^T & \text{ if else.}
\end{array}
\right .$$
The mesh, the initial solution, an intermediate solution and the final one at $t=1.8$ are displayed on figure \ref{remi:fig:2}. They are obtained with a third order accurate time scheme and quadratic B\'ezier approximation. It can be observed that the sine wave is not damped, though the resolution of the mesh is not very fine (all degrees of freedom are represented, there are $8985$ dofs).
\begin{figure}[h]
\begin{center}
\subfigure[]{\includegraphics[width=0.45\textwidth]{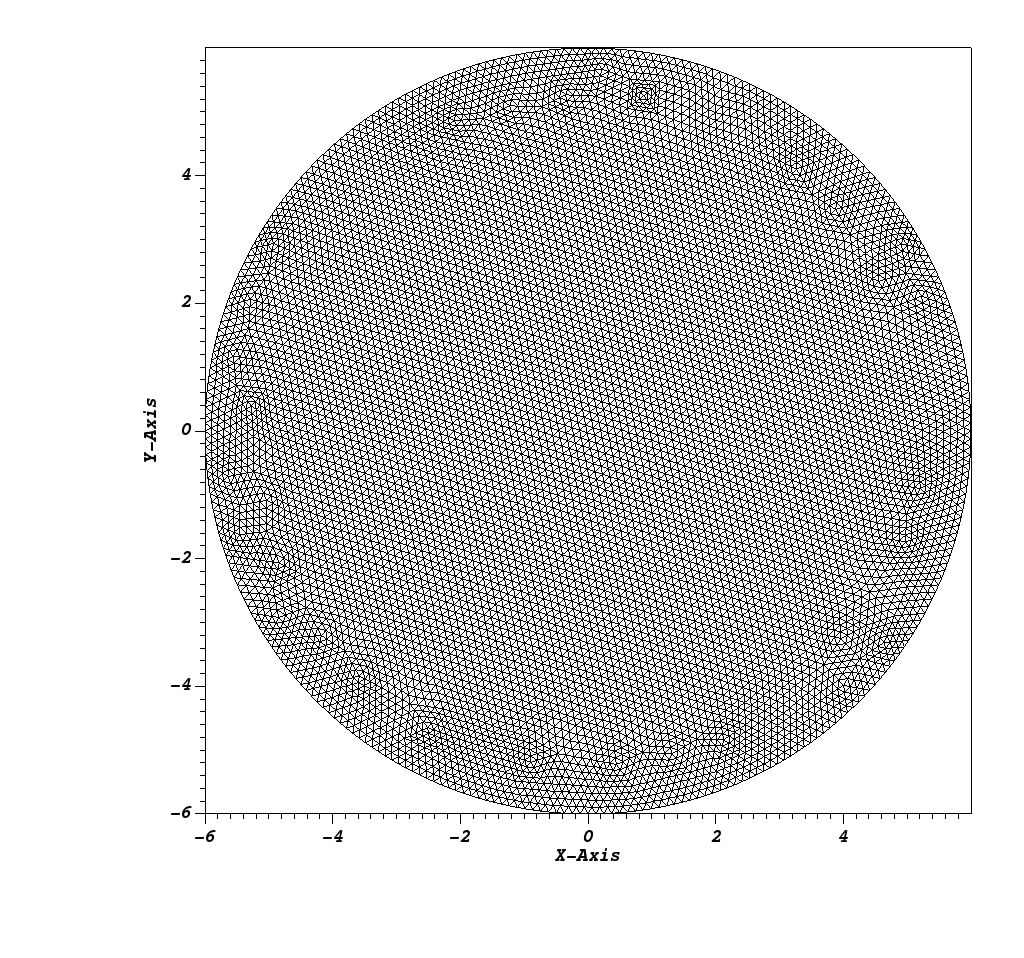}}
\subfigure[]{\includegraphics[width=0.45\textwidth]{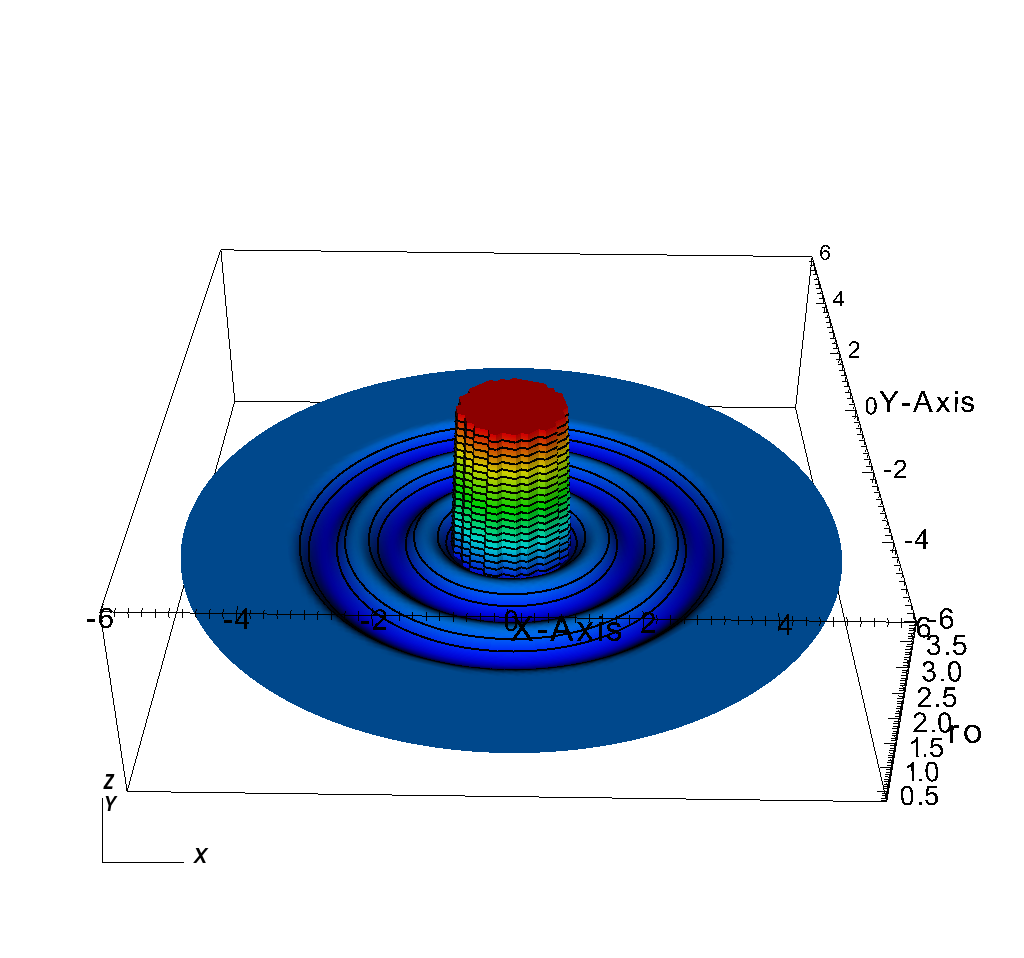}}
\subfigure[]{\includegraphics[width=0.45\textwidth]{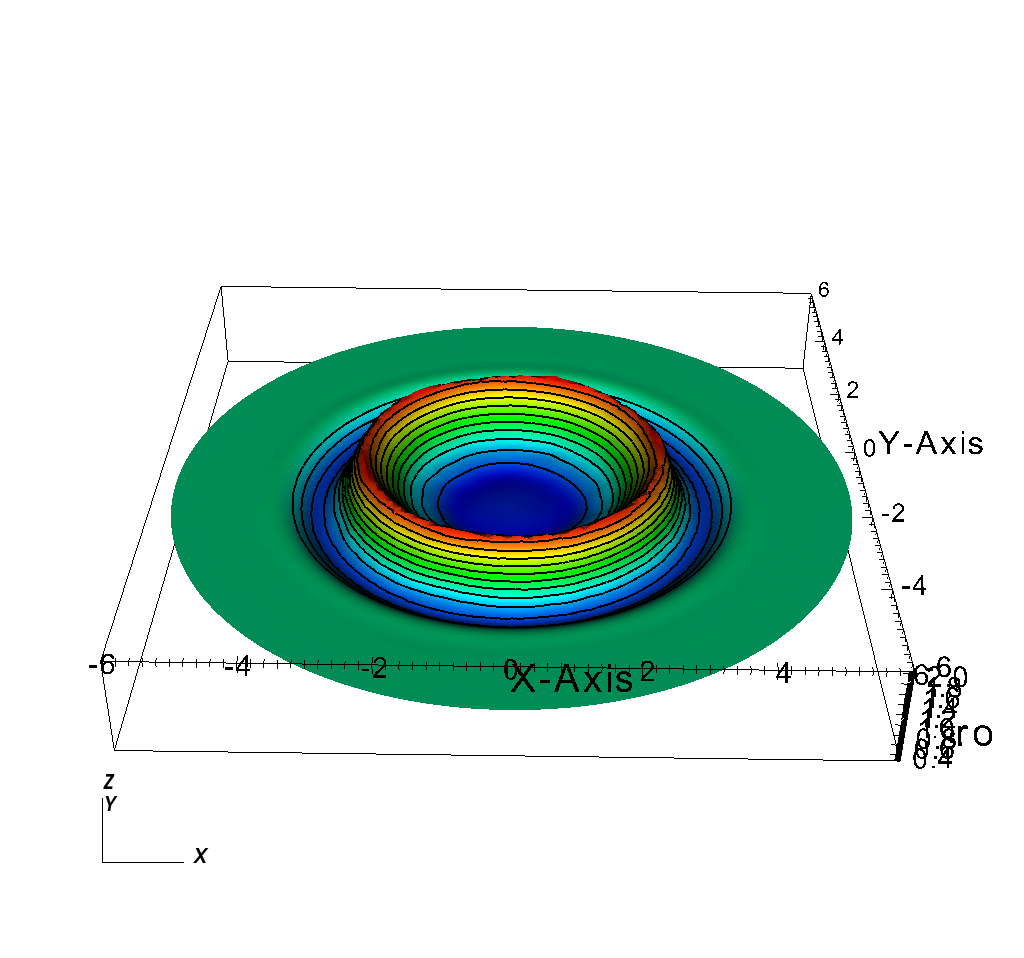}}
\subfigure[]{\includegraphics[width=0.45\textwidth]{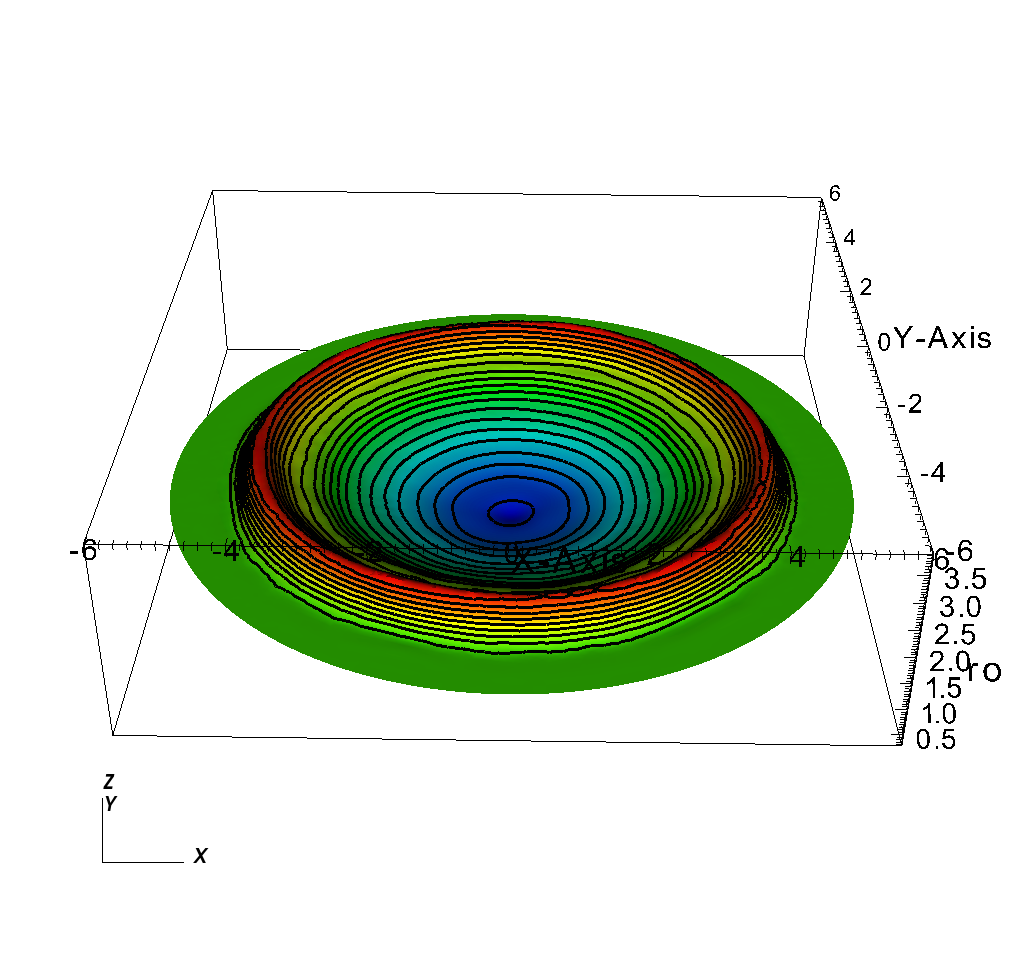}}
\end{center}
\caption{\label{remi:fig:2} (a): mesh, all the degrees of freedom are represented, (b): initial condition, (c) Solution at an intermediate time, (d) Solution at $t=1.8$. The CFL is $0.125$.}
\end{figure}
The schemes are a combination of the Galerkin scheme with jump stabilisation and the LLFs with quadratic/cubic approximation, see \cite{paola2} for more details.

\subsection{Shallow water and  the lake at rest state state}

We comment on some results \mario{of the shallow water equations} originally appeared in \cite{RB09,ArR:18}. 
The test considered is a quite  classical  
a perturbation of the lake at rest state, with a bathymetry defined by a smooth exponential hump. We refer to the original papers, and references therein, for details. 
\mario{The scheme used is the }
non-linear  limited Lax-Friedrich's residual distribution, described in section \S\ref{sec:conservation},  with \mario{appropriate modifications of the mass matrix and stabilization operators to
handle both smooth and discontnuous flows, whicle accounting for  well balanced and wet/dry transitions}
 (see  \cite{R15} for the explicit method on fixed meshes).  \mario{The same polynomial approximation is used for } 
 the conservative variables 
 $\bu^h$,  \mario{and for }
 the topography $b^h$. \mario{The quadrature 
strategy is exact on lake at rest solutions (cf.   section \S \ref{sec: bien roule}).} The main difference  between the two references is that in  \cite{RB09}  no special care is taken in handling the time derivative, and a fully implicit   (in time) approach 
is used, based on a second order trapezoidal   method.
  On the contrary, in \cite{ArR:18} the authors have combined the error correction method discussed in this paper, with  a well balanced ALE formulation on moving adaptive meshes.  
\mario{The high order discrete remapping of the initial topographic data recalled in 
section  \S \ref{sec: bien roule} is used to preserve mass conservation on moving meshes, within }
the  errors of  the quadrature formulas used in the remap. 
\begin{figure}
\begin{center}
\includegraphics[width=0.4\textwidth,valign=c]{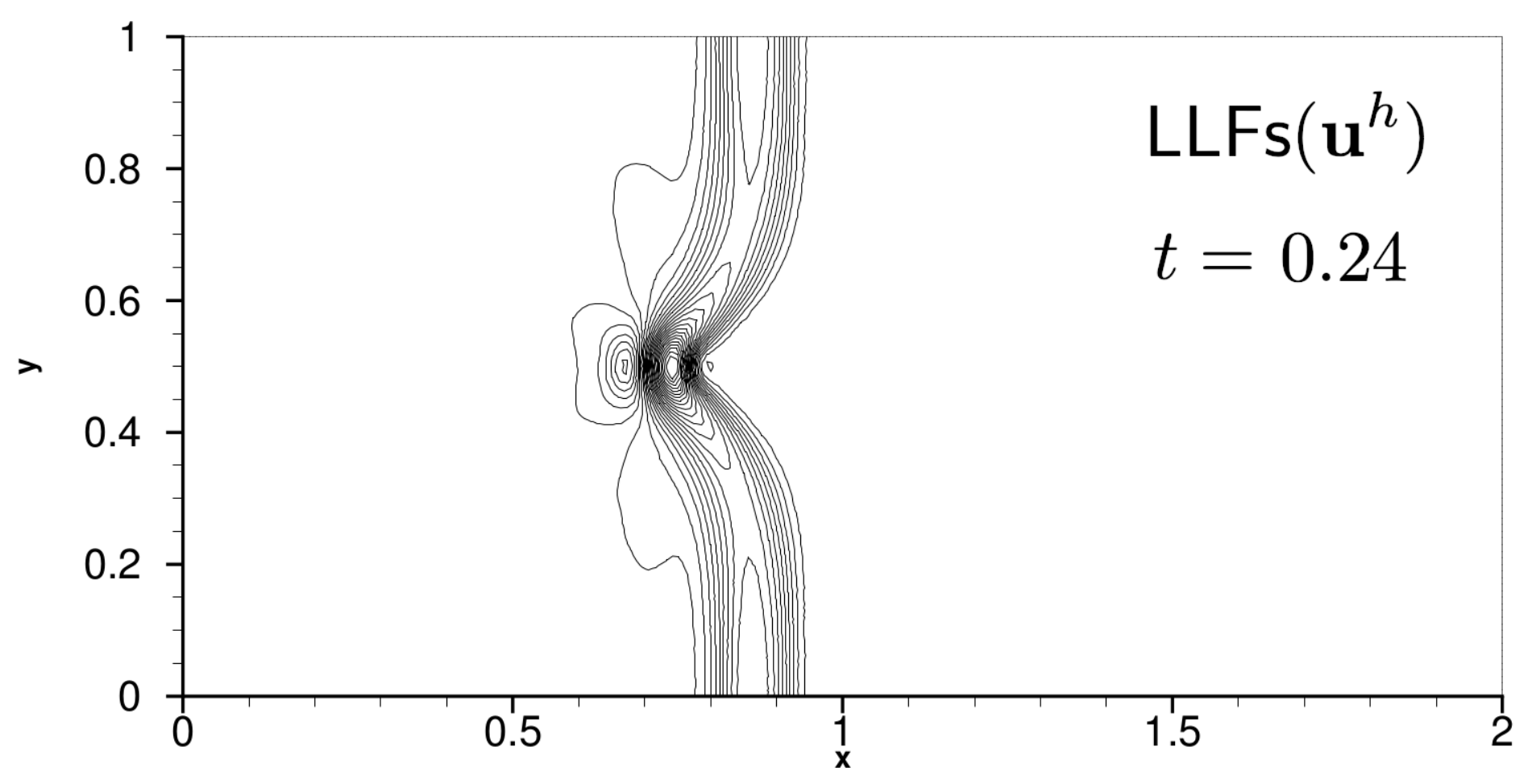}\includegraphics[width=0.4\textwidth,valign=c]{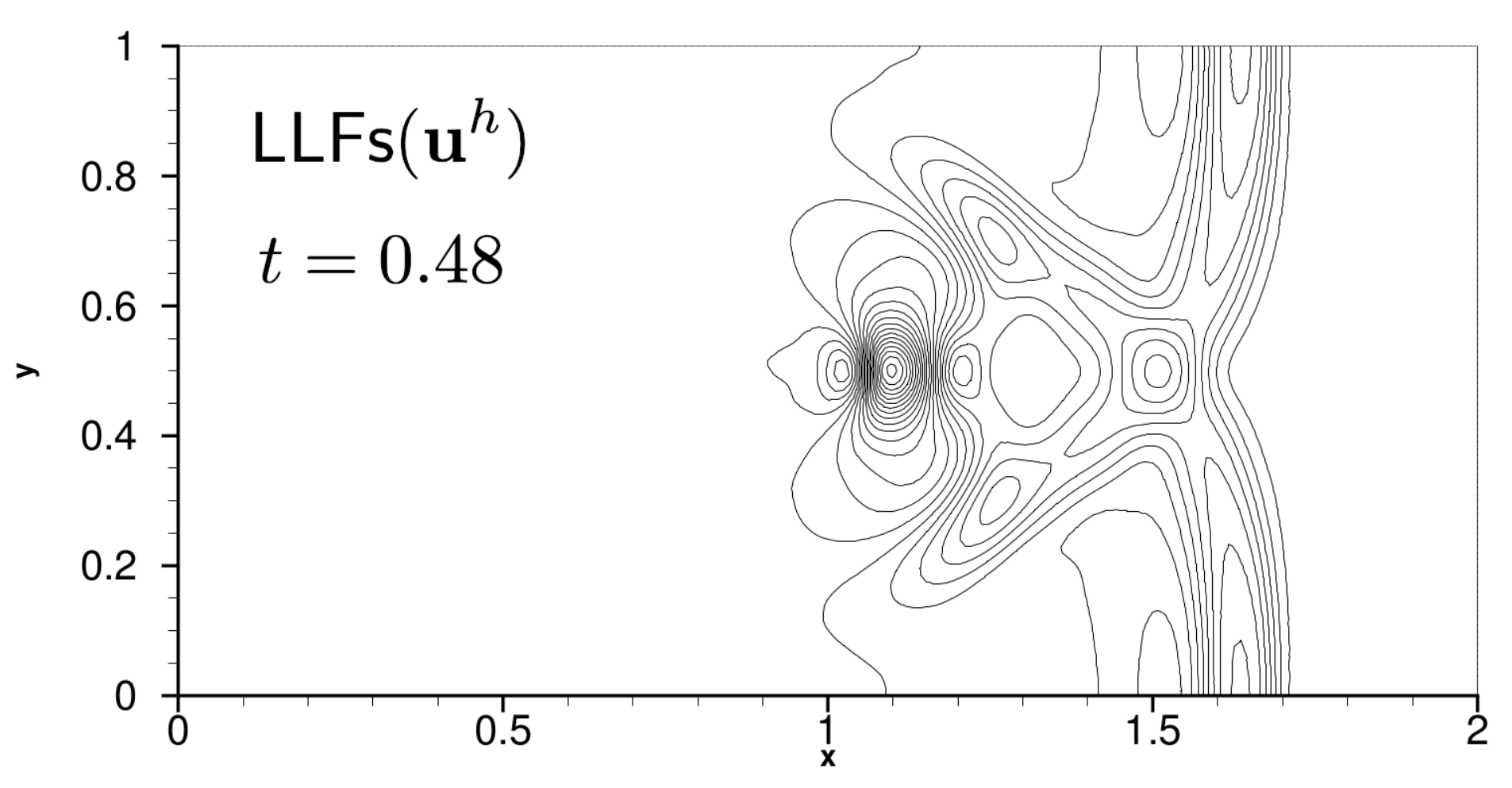}
\includegraphics[width=0.4\textwidth,valign=c]{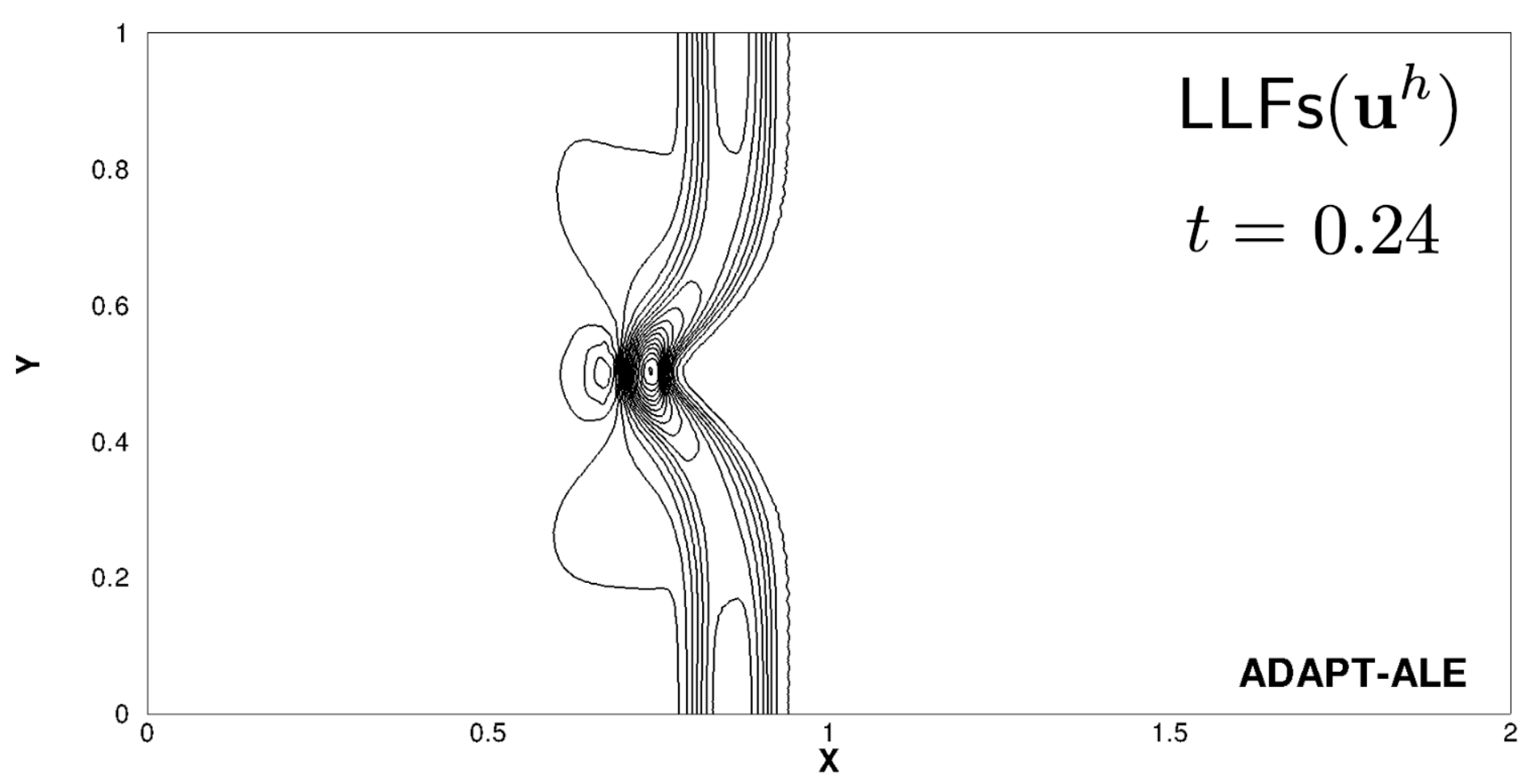}\includegraphics[width=0.4\textwidth,valign=c]{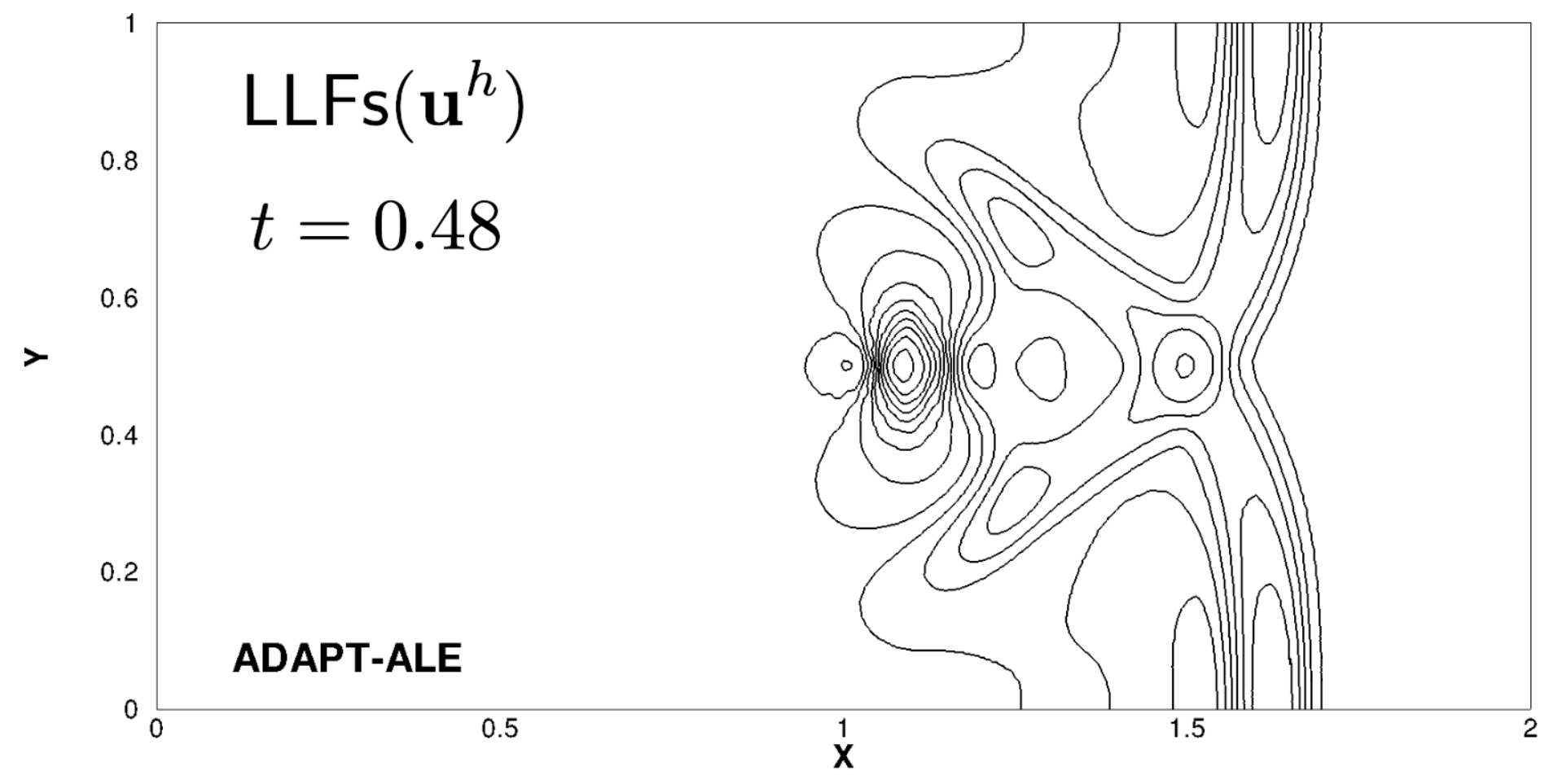}
\includegraphics[width=0.4\textwidth,valign=c]{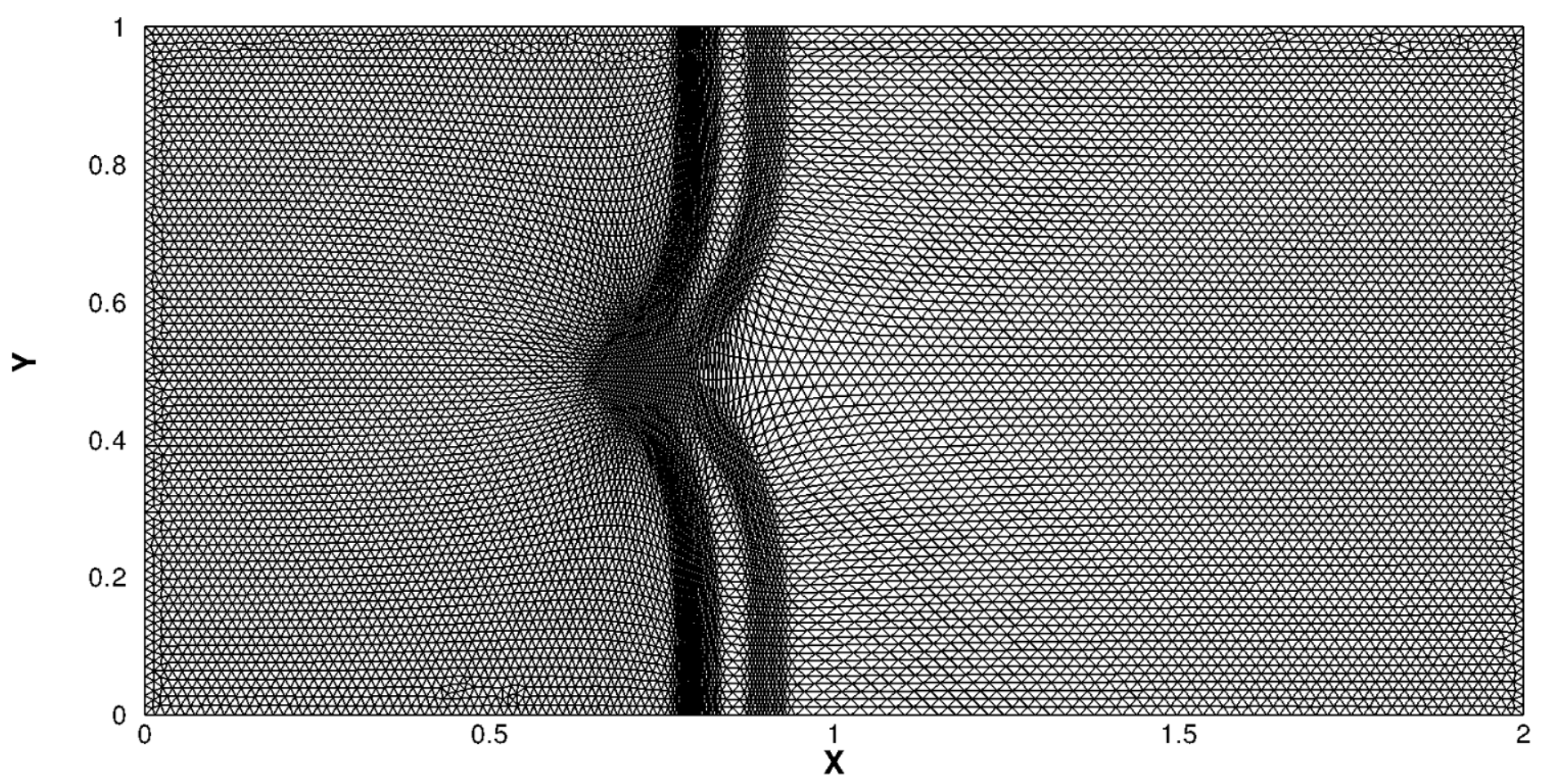}\includegraphics[width=0.4\textwidth,valign=c]{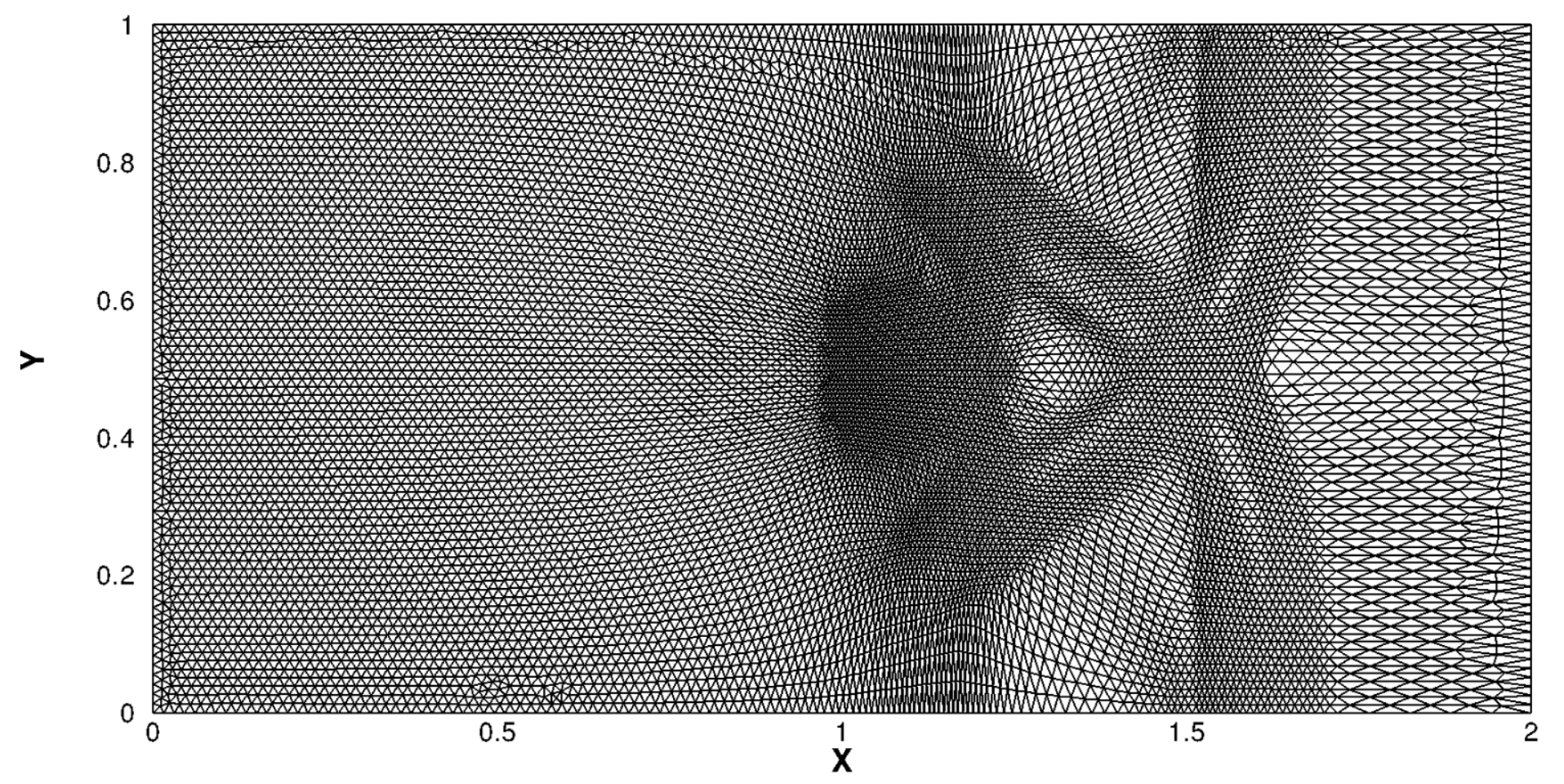}
\end{center}
\caption{\label{fig:lake1}Perturbation of  the lake at rest over a smooth topography: free surface evolution. Top: implicit scheme of \cite{RB09}. Middle and bottom : explicit adaptive moving  mesh approach by \cite{ArR:18} }
\end{figure}

We report in  the top and mid-rows of figure \ref{fig:lake1}  the contours of the free surface level obtained respectively with the implicit scheme (top), and
with the explicit adaptive approach (middle).  The results  allow to compare the evolution of the initial perturbation, in both cases unspoiled by any unphysical oscillations 
 in  the free stream relation,   which is the  main interest of  constructing  well balanced schemes. 
The bottom row shows the moving adaptive meshes produced with technique proposed in \cite{ArR:18} which follow very closely the wave pattern.

Despite of the fact that the contours of the implicit method are slightly crisper than those of the explicit one,  the water levels obtained are very similar.  This is confirmed by the cuts along the centreline, 
reported on figure  \ref{fig:lake2}. The plots in this figure show that the adaptive computations on a relatively coarse mesh  (roughly 12k nodes) compare very favourably in terms of  the peaks and troughs  of the free surface
 with those  of  the implicit scheme on a finer grid (roughly 20k nodes), and with a reference solution of obtained with the  explicit scheme on a fine mesh (roughly 50k nodes).
 The water levels on the unadapted coarse mesh are also reported to show the substantial benefit of adaptation.

\begin{figure}
\begin{minipage}{0.5\textwidth}
\begin{center}
\includegraphics[width=0.8\textwidth,valign=c]{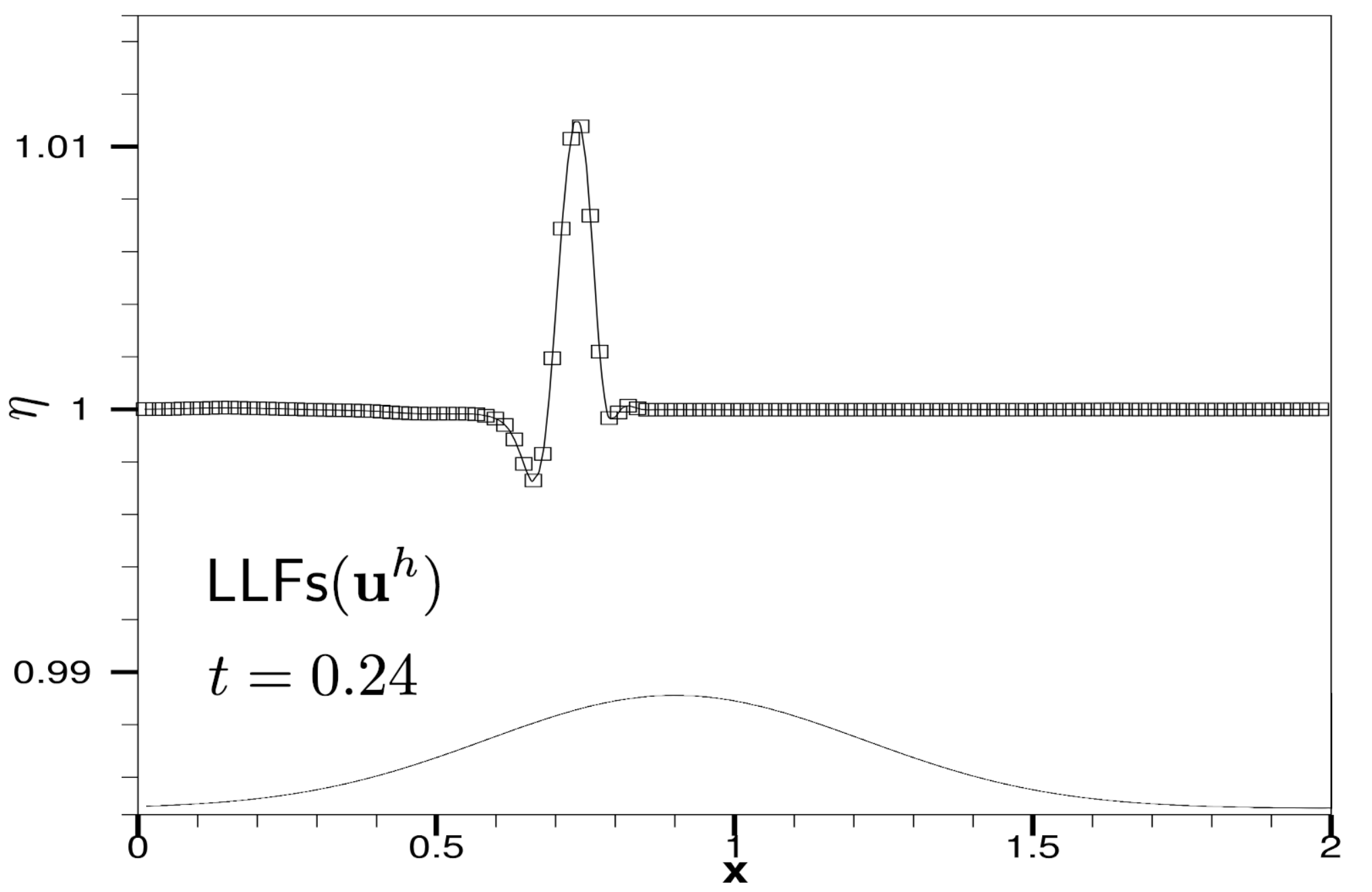}
\includegraphics[width=0.8\textwidth,valign=c]{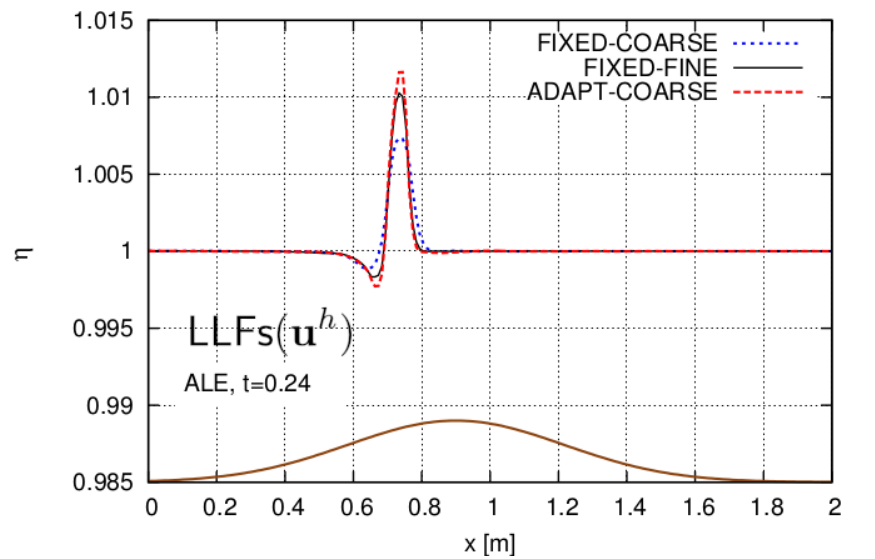}
\end{center}
\end{minipage}\hfill
\begin{minipage}{0.5\textwidth}
\begin{center}
\includegraphics[width=0.8\textwidth,valign=c]{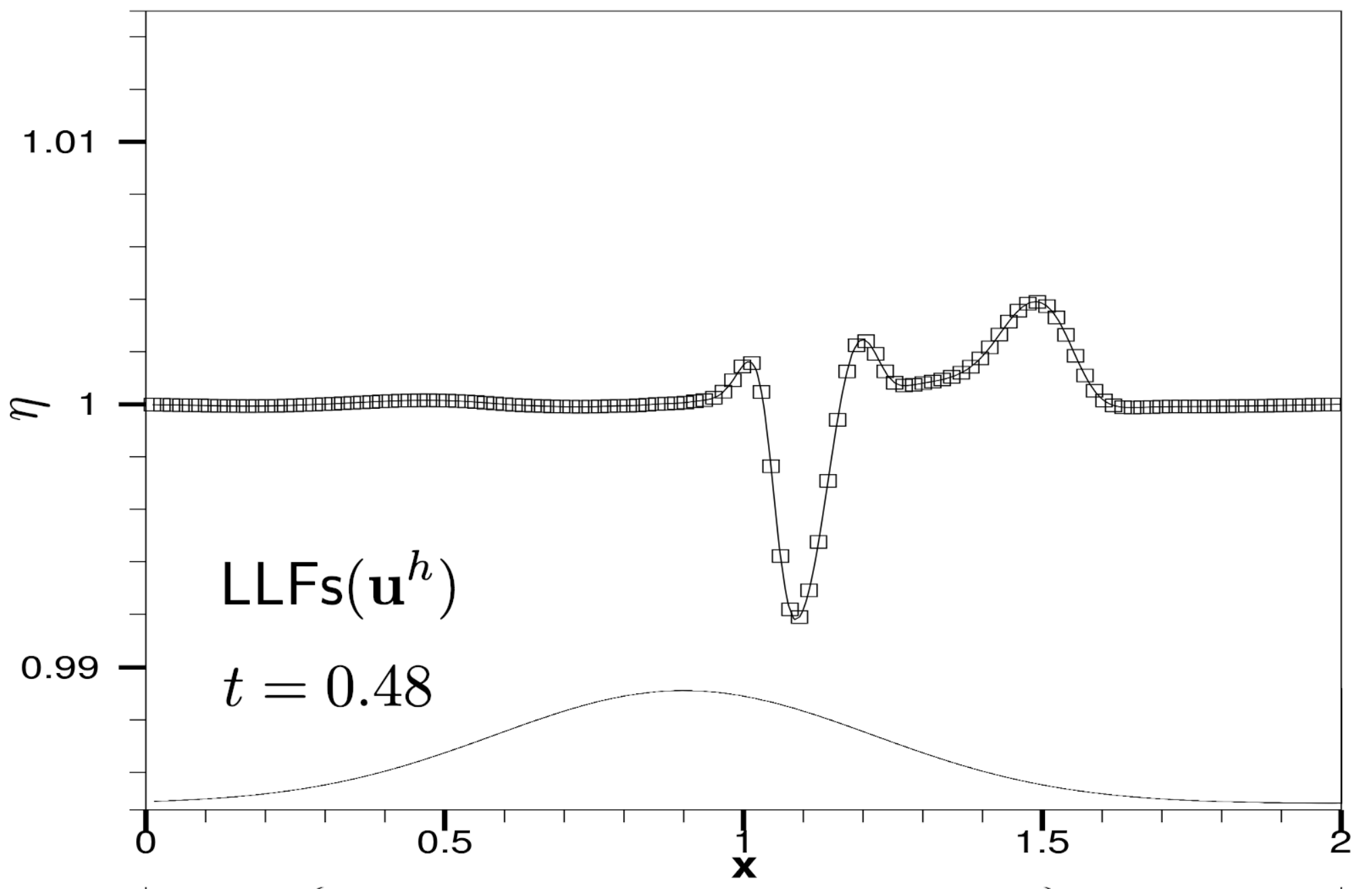}
\includegraphics[width=0.8\textwidth,valign=c]{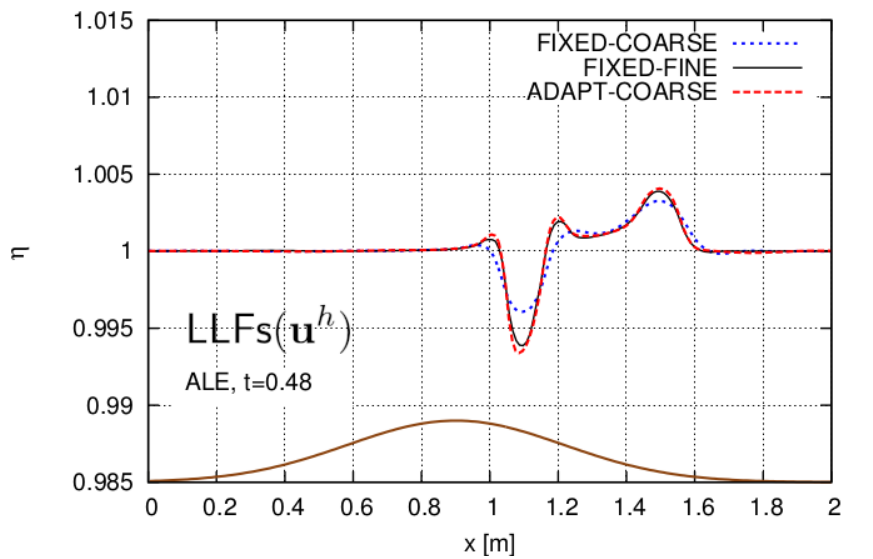}
\end{center}
\end{minipage}
\caption{\label{fig:lake2} Perturbation of  the lake at rest over a smooth topography: centerline free surface levels.  Top: implicit scheme of \cite{RB09}.  Bottom : explicit adaptive moving  mesh approach by \cite{ArR:18} }
\end{figure}

On the left in figure \ref{fig:lake3}  we report  the total mass error   for different quadrature strategies in the bathymetric remap. For this smooth case the error levels are already small 
with second order quadrature ($r$ in the figure denotes the degree of exactness of the integration formulas).  For higher order formulas, the error remains practically at machine zero.
Finally, the table on the right reports CPU times for the explicit computations, which show a  computational gain of almost 35\% in time for the adaptive method compared to the fine mesh results.
This figure could be further improved with some local or partitioned time marching method. To conclude on this aspect, we refer to \cite{R15} for similar comparisons between the explicit and implicit residual
schemes on fixed meshes, shown that for this type of problems, for a given resolution  the explicit  approach can be 3 to 6 time faster.

\begin{figure}
\begin{minipage}{0.6\textwidth}
\includegraphics[width=0.8\textwidth]{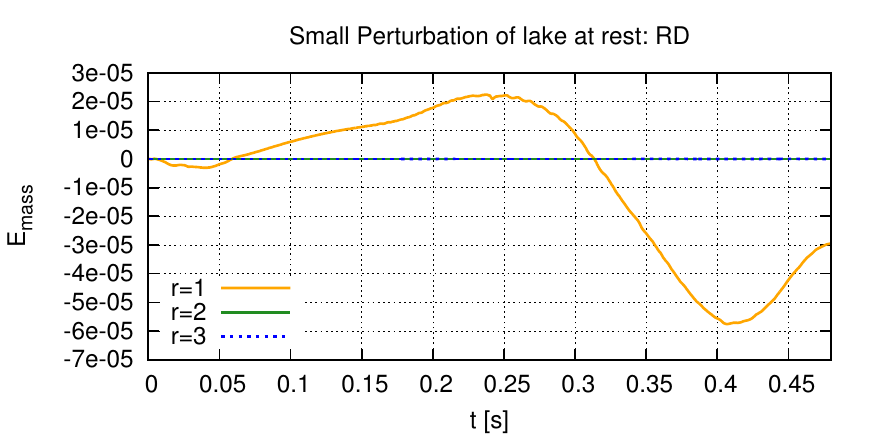}
\end{minipage}\hfill
\begin{minipage}{0.4\textwidth}
\footnotesize
\noindent \begin{center}
\begin{tabular}{|c|c|c|}
\hline 
Mesh      &  Nodes    & CPU  [s]  \tabularnewline
 \hline 
 \hline 
Fix & 12142  & 73.60    \tabularnewline
 \hline 
Fix& 50631  & 711.08   \tabularnewline
 \hline 
Adapt+ALE &  12142 & 204.96   \tabularnewline
 \hline  
\end{tabular}
\par\end{center}
\end{minipage}
\caption{\label{fig:lake3}  Mass conservation in function of the quadrature order for the  topography remap, and CPU times for the explicit schemes} 
\end{figure}

\subsection{Shallow water and moving equilibria}

We consider here three applications involving moving equilibria  in shallow water flows with constant total energy \eqref{lcef}.
The first two applications involve a very classical situations with a $C^0$ bathymetry defined as 
$$
b(x)=\left\{\begin{array}{ll}
\dfrac{1}{5}\left(1 - \dfrac{(x-10)^2}{4}\right)\quad&\text{if }x\in[8,12]\\[10pt]
0&\text{otherwise}
\end{array}\right..
$$
The prescribed values of total energy and mass flux  are :  
$$[\mathcal{E}_0,\,q_0]=[22.06605 m^2/s^2,\,4.42m^2/s  ].$$ 
The spatial domain  has a horizontal length of 25$m$.\\

The first test, \remi{see figure \ref{fig:nrg1d}}, consists of perturbing the 1D steady state  corresponding to the above choices within the slice $x\in[6.5,7.5]$. 
The perturbation is added to the  free surface, and has a magnitude of  0.05$m$.  \mario{ We use a  structured triangulation  containing with
spacing $25/200$, and  
periodic boundary conditions in the vertical direction. The scheme used is the implicit
 LLFs scheme of \cite{RB09},  with 
 approximation done in terms of the steady invariants $\bv=[\mathcal{E},h\bbv]$, and with a }
 piecewise linear approximation  of the 
bathymetry.  We refer to \cite{nxs07,R15} for the practical implementation of this choice,
which requires non-linear iterations to  locally invert the mapping $(\bv,b)\mapsto \bu$. 

 {The results are in excellent agreement, at least qualitatively with similar results of the litterature, using dG or WENO schemes.}
 
\begin{figure}[h]
\begin{center}
\hspace{-1cm}\begin{tabular}{cccc}
\includegraphics[width=0.215\textwidth]{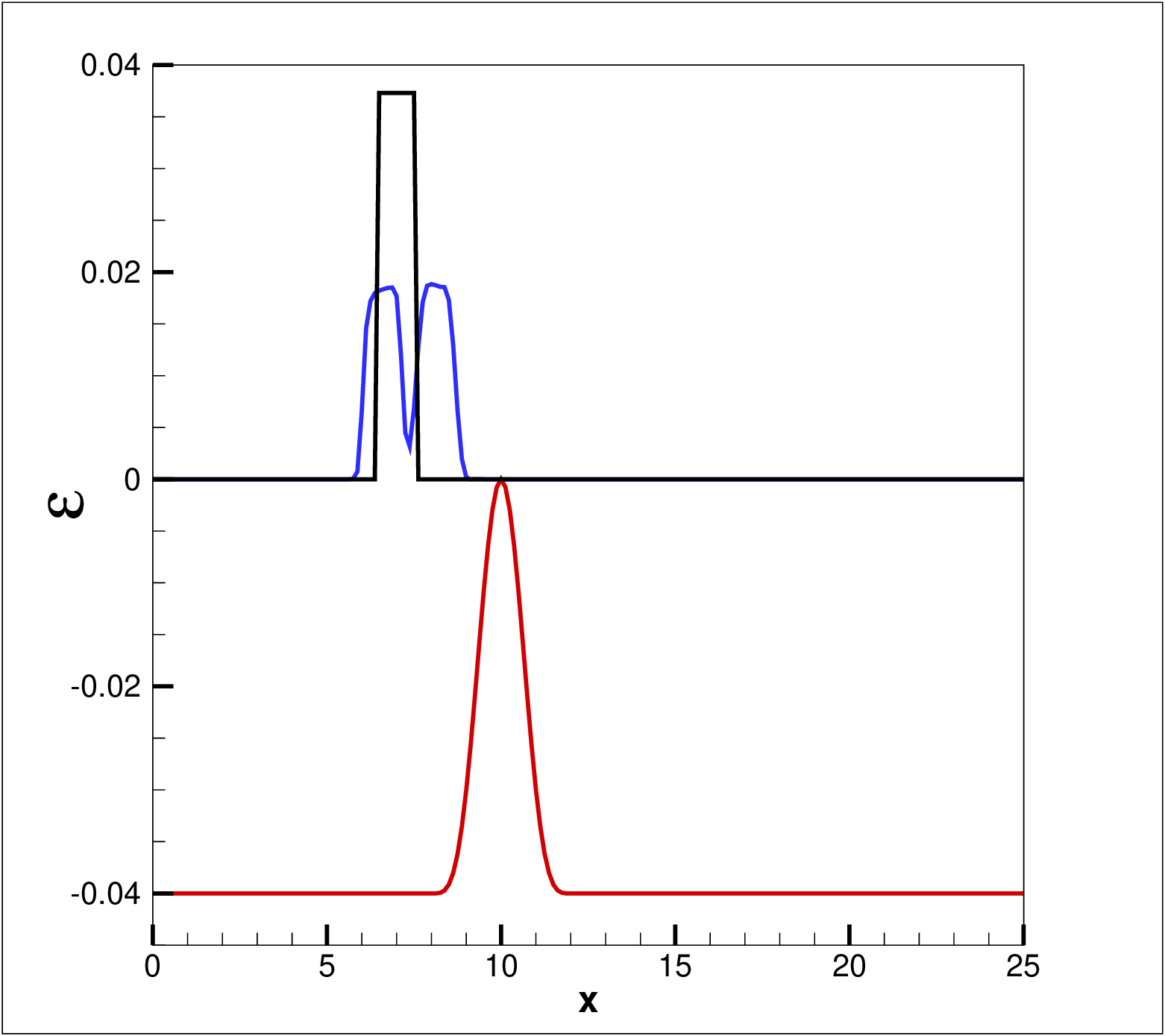}&
\includegraphics[width=0.215\textwidth]{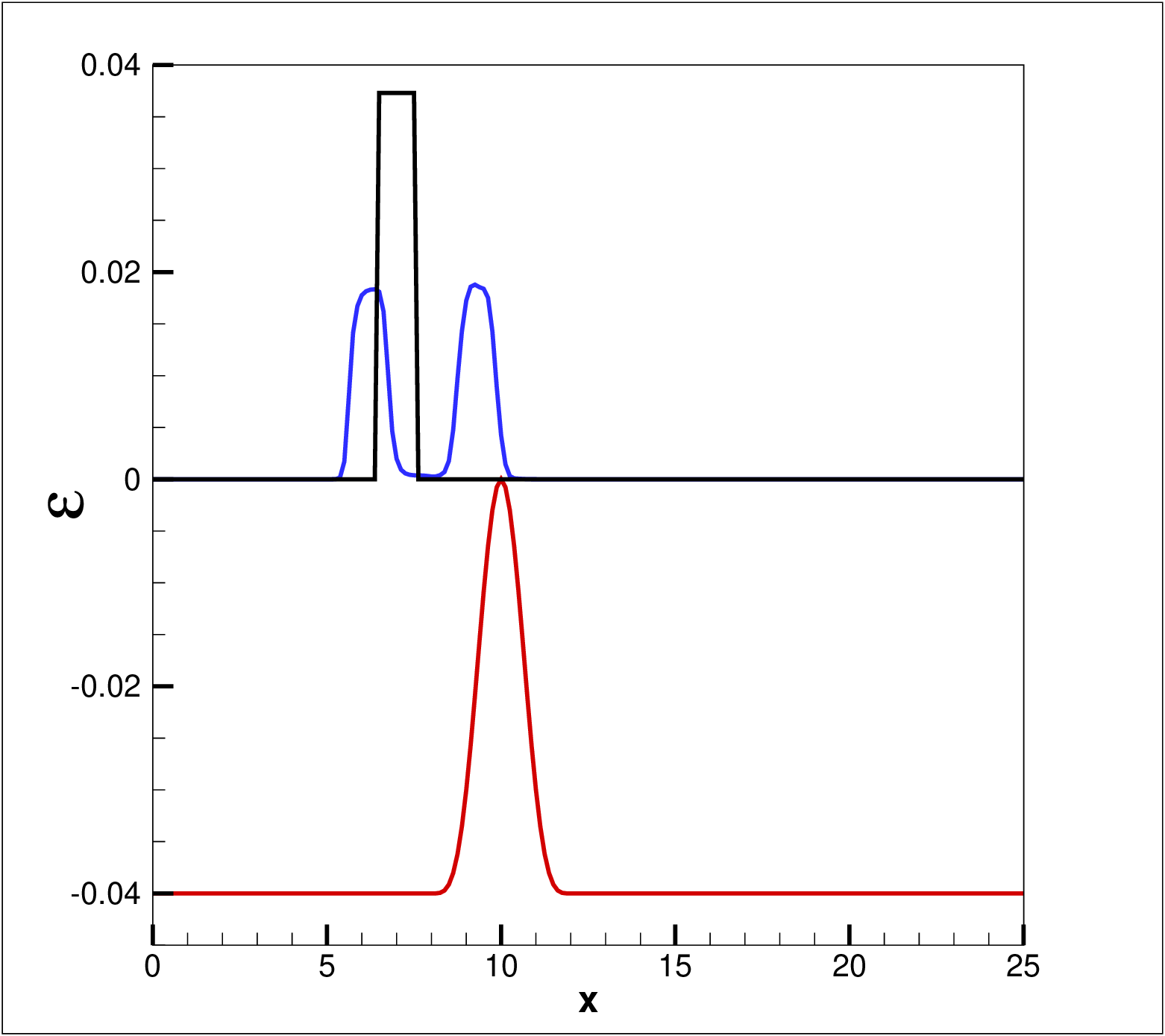}&
\includegraphics[width=0.215\textwidth]{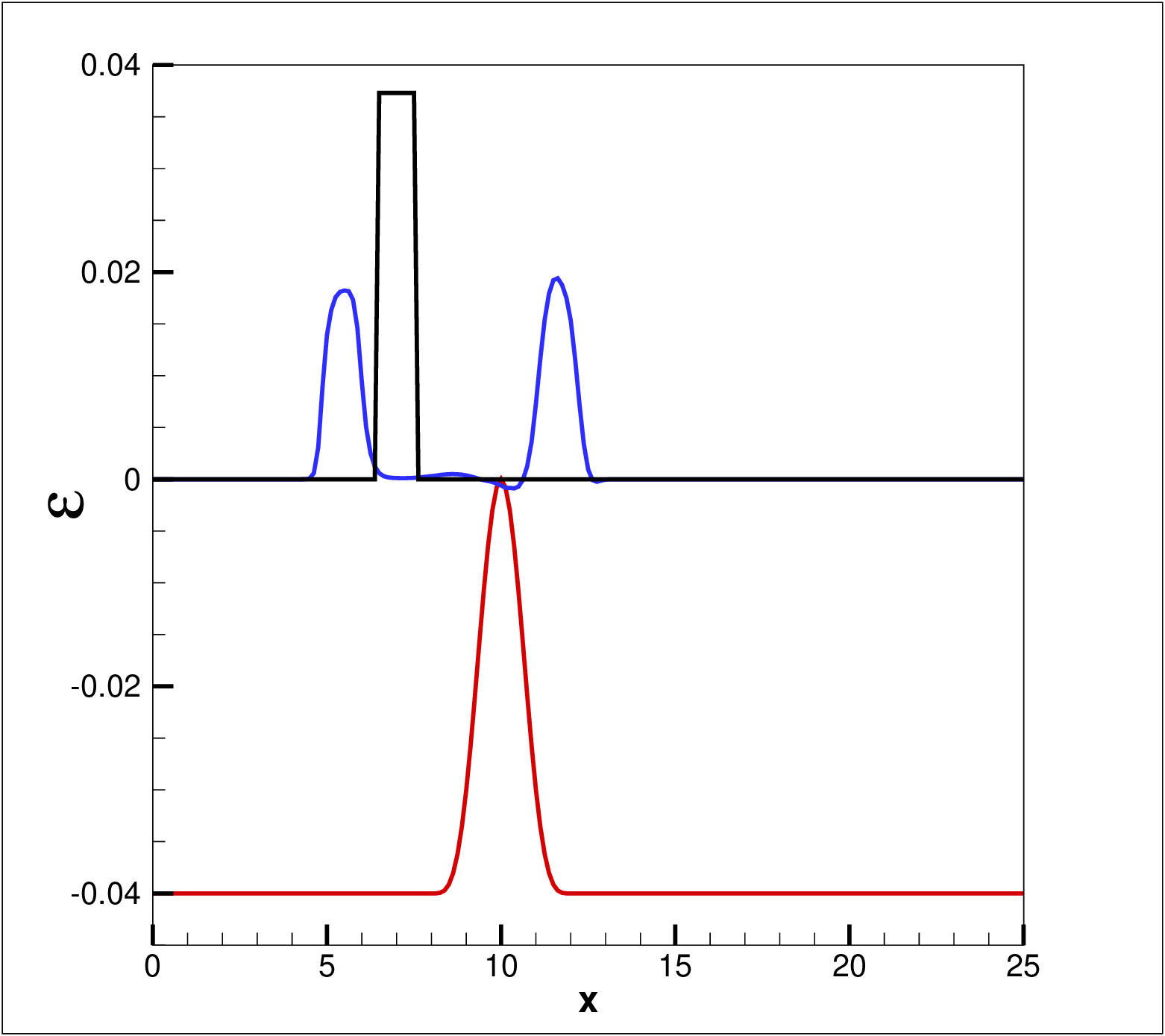}&
\includegraphics[width=0.215\textwidth]{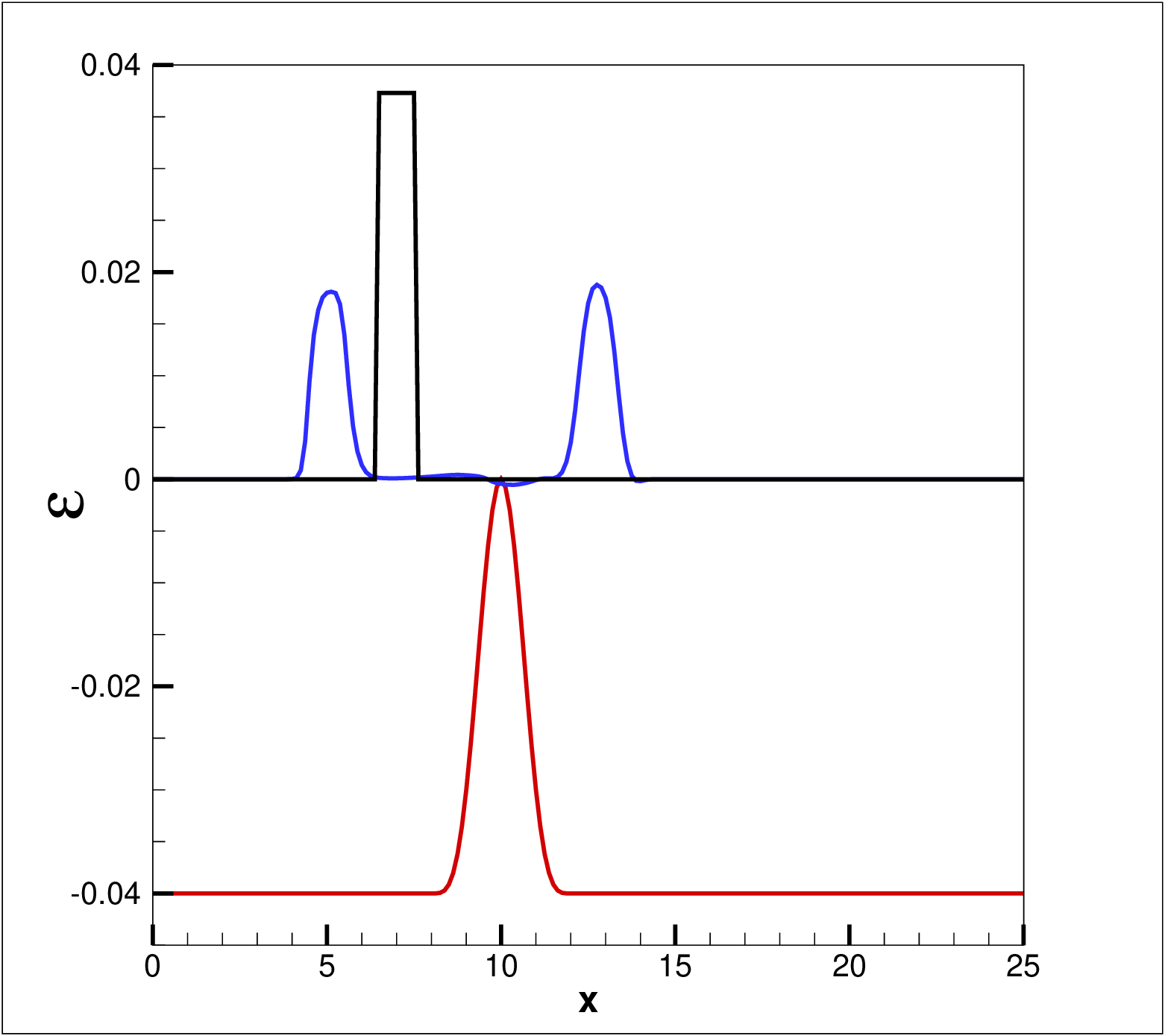}\\
\includegraphics[width=0.215\textwidth]{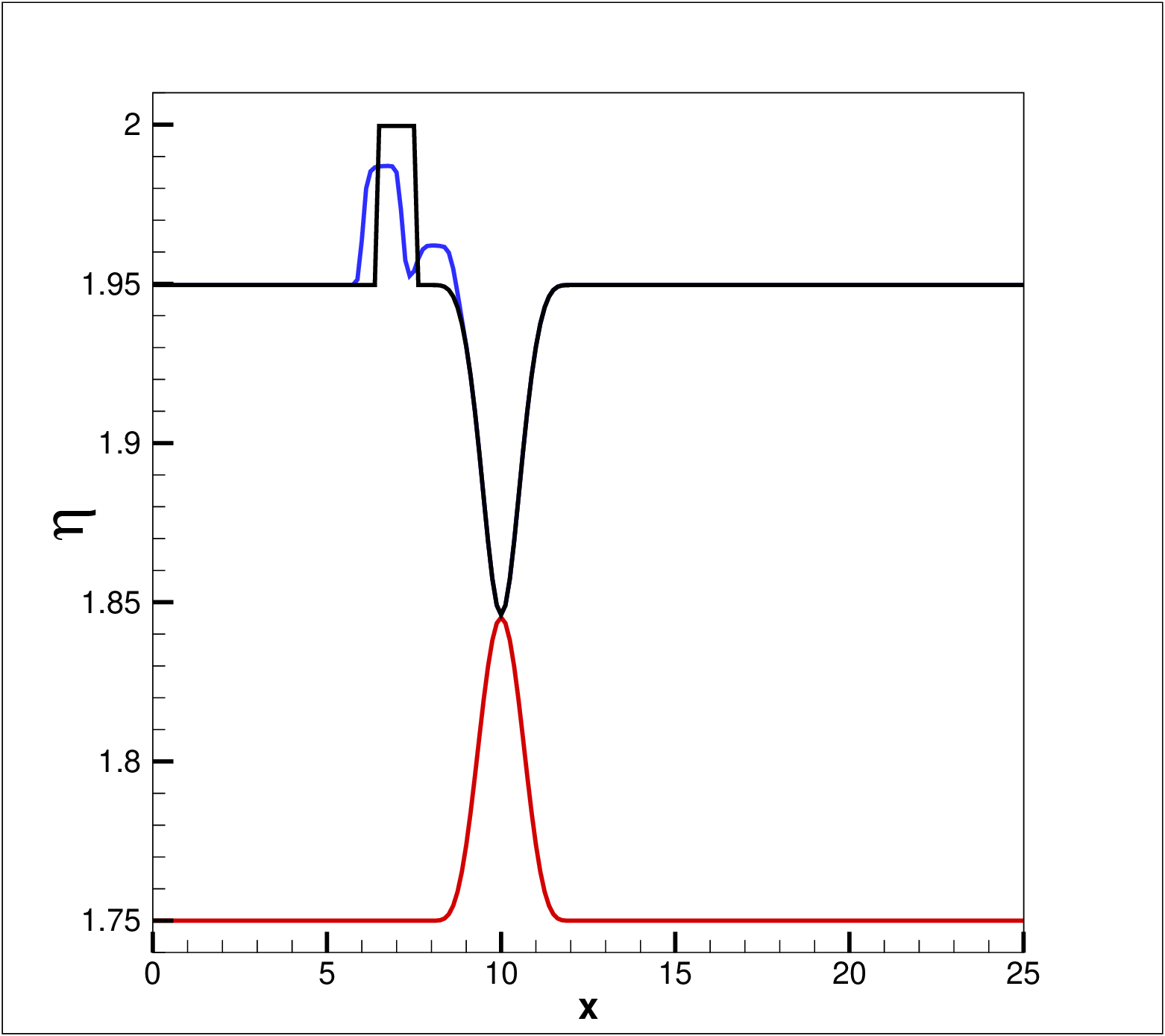}&
\includegraphics[width=0.215\textwidth]{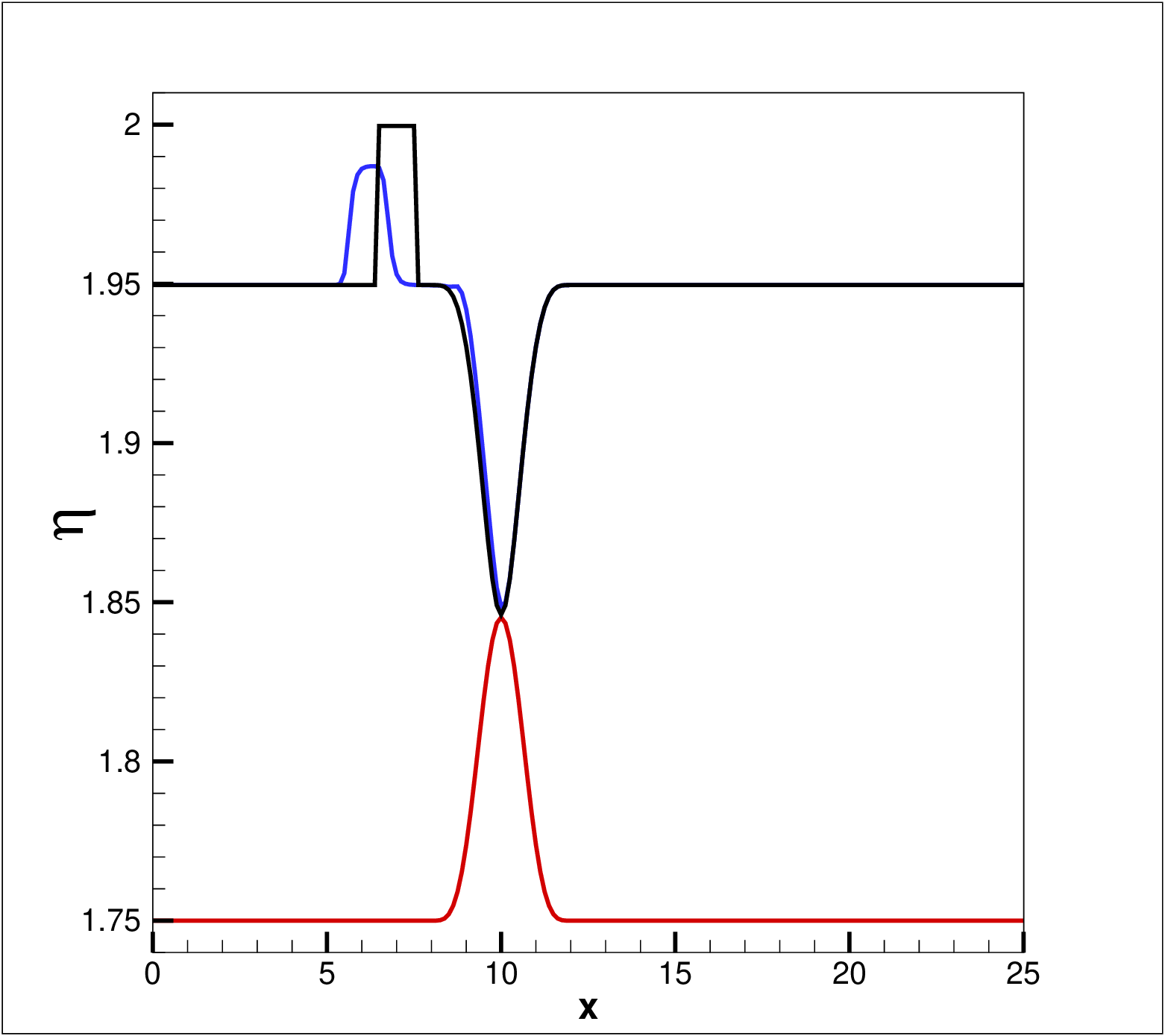}&
\includegraphics[width=0.215\textwidth]{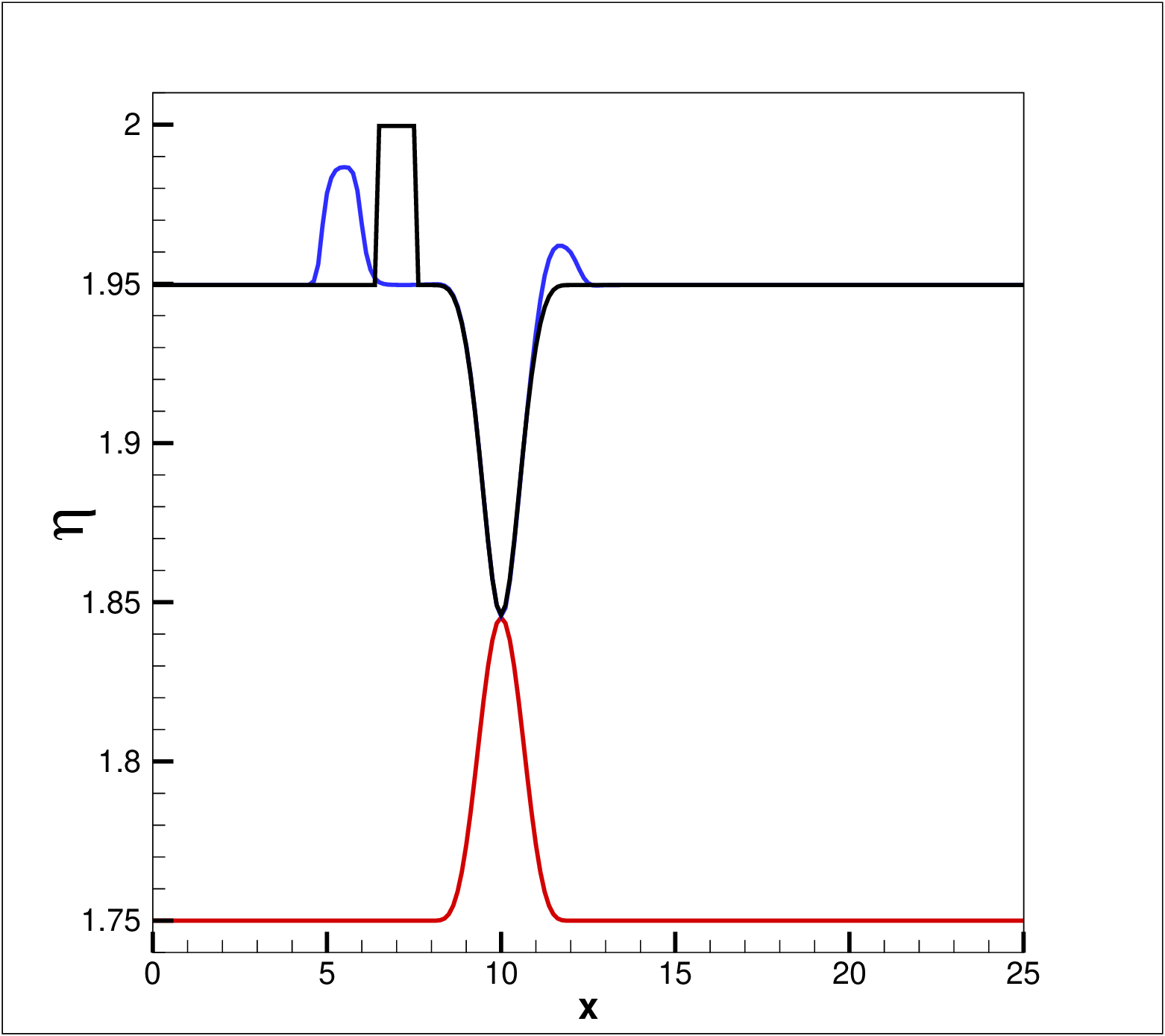}&
\includegraphics[width=0.215\textwidth]{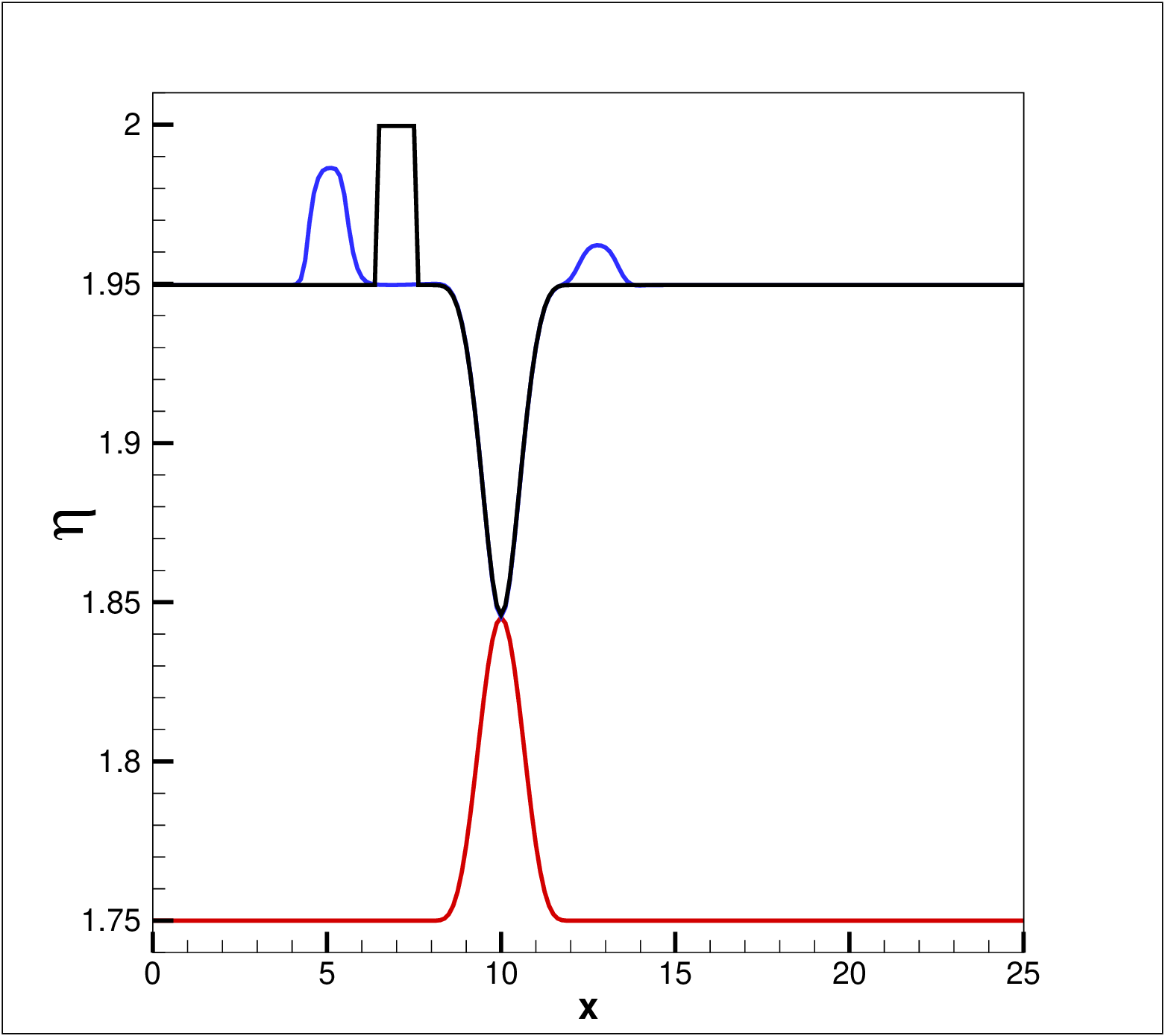}
\end{tabular}
\end{center}
\caption{\label{fig:nrg1d}{Example of a 1D perturbation in a constant energy steady flow.s}}
\end{figure}

 \medskip
 
 We then consider a similar test, but on the 2D domain $[0, 25]^2$ (see \cite{R15} concerning the relevance of \eqref{lcef} in 2D).
 We add a 2D perturbation to the same one-dimensional steady state obtained with the choices described above. As before, a perturbation  of 0.05$m$ is added to the free surface $\eta$,
only  this time in the subdomain $[6.5,7.5]\times[12,13]$. We \mario{report on figure \ref{fig:nrg2d} snapshots of }
the  perturbation $\delta\eta:=\eta-\eta_{\text{steady}}(x)$, with  $\eta_{\text{steady}}(x)$ the free surface level 
corresponding to the exact 1D steady  solution. 
We compare the results obtained \mario{ on a structured triangulation with spacing $25/200$}
 with the straightforward application
of the  LLFs scheme, denoted by LLFs$(\bu^h)$  \mario{exactly well balanced only for the lake at rest state)},
with the 
same scheme  based on the approximation of the steady invariants  (denoted by LLFs$(\bv^h)$). In both cases the standard linear approximation of the bathymetry $b^h$ is used.

\begin{figure}[h]
\begin{center}
 \includegraphics[width=0.4\textwidth]{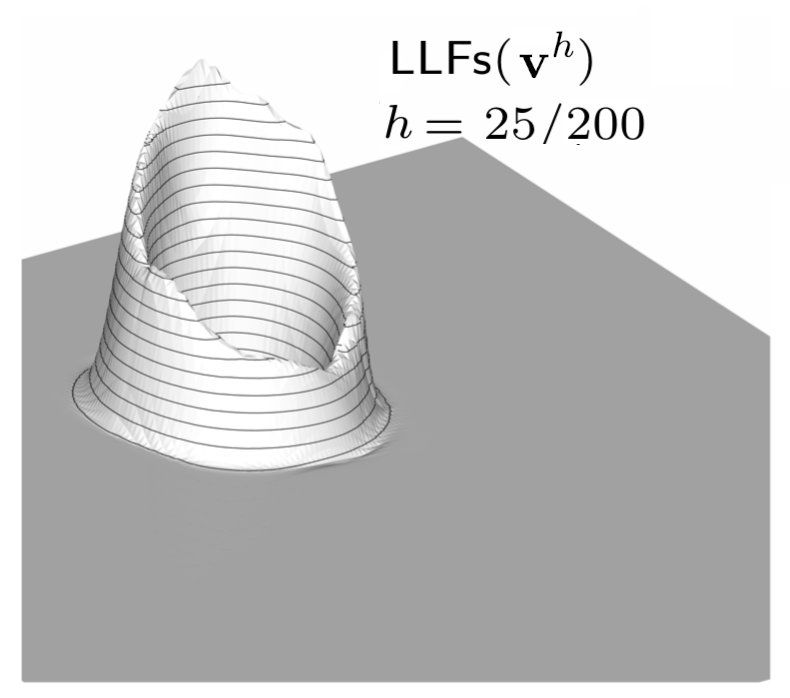} \includegraphics[width=0.4\textwidth]{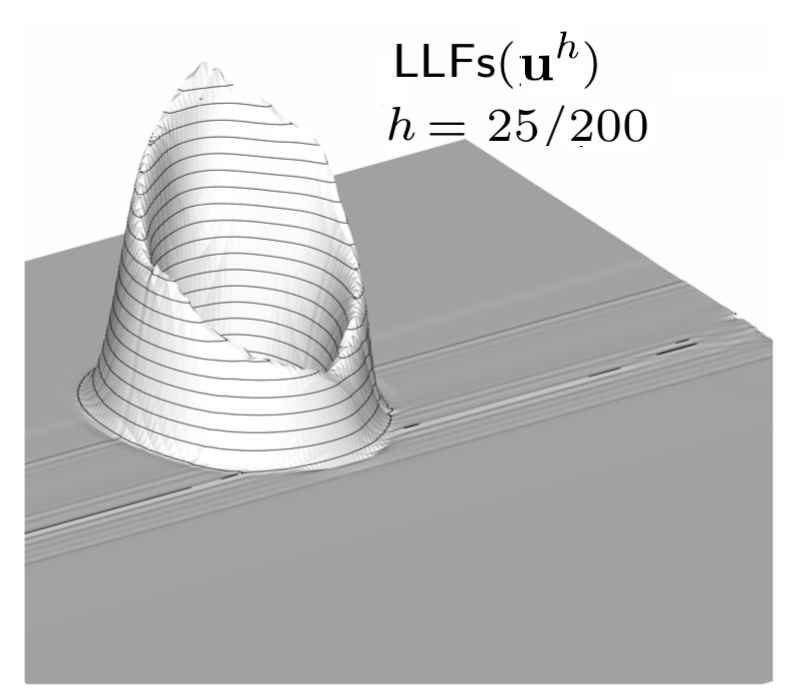}  \includegraphics[width=0.4\textwidth]{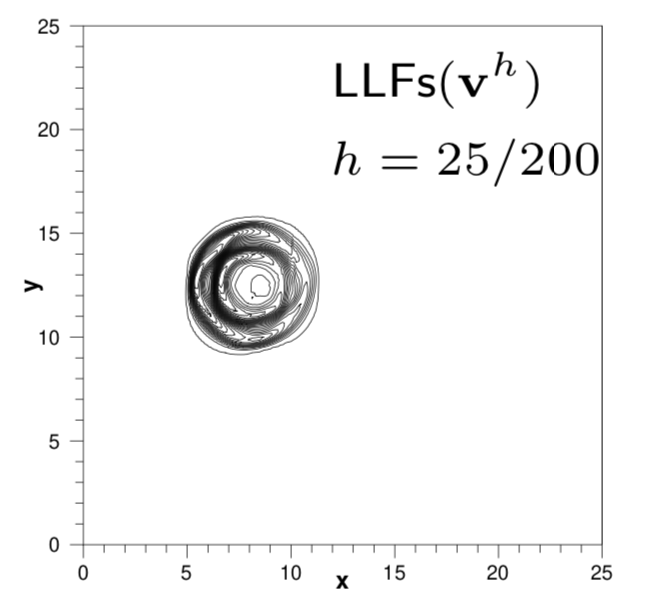}\includegraphics[width=0.4\textwidth]{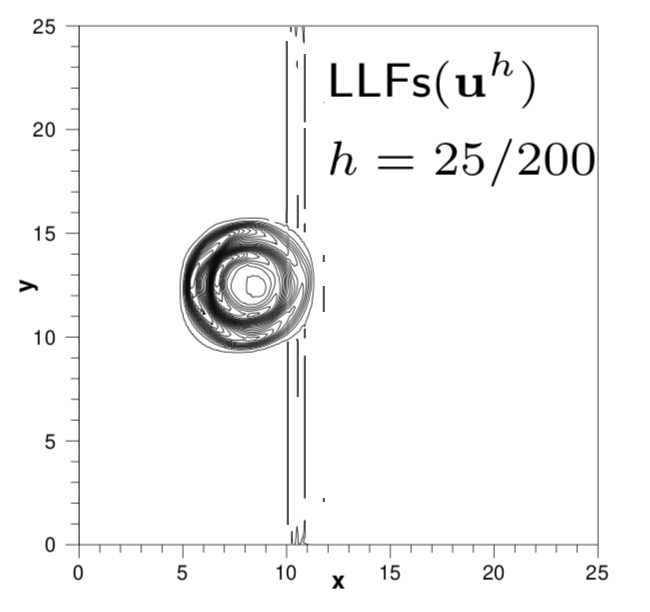}
 \end{center}
\caption{\label{fig:nrg2d}{Example of a 2D perturbation in  a constant energy steady flow.}}
\end{figure}

\mario{Also in in 2D  the LLFs$(\bv^h)$ provides a perfect evolution of the perturbation with no visible spurious effects.}
The LLFs$(\bu^h)$ scheme shows such perturbations. However, we stress very strongly that on such meshes these are of the order of $10^{-4}m$, thus several orders of magnitude smaller  than the  initial 
one added to the free surface. Of course, these would become relevant had we  reduced the magnitude of the initial perturbation. This is  a very nice result for obtained with the ``straight out of the box'' 
application of the residual distribution method.  \mario{This kind of result is not obtained by a straightforward application of  finite volumes or dG schemes.} 

Finally, we \mario{consider a}
a genuinely 2D configuration  obtained by replacing he structured triangulation with an unstructured one.  The bathymetry is 
 now is defined by a series of $C^1$ sinusoidal ribs (see  \cite{R15} for details).  
  We \mario{compare the 
   LLFs$(\bv^h)$ scheme with approximation in steady invariants, and analytical data $b(x)$ and $b'(x)$ used in the residual evaluations,
   and the standard LLFs$(\bu^h)$ method using piecewise linear solution and data. Unstructured triangulations are used}.
\mario{The first scheme  fits the}
 hypotheses of Proposition 2 (see section \S \ref{sec: bien roule}.2), and we can  \mario{see on the leftmost picture of figure \ref{fig:nrg2d-1}  that}
  indeed convergence to the exact solution w.r.t quadrature accuracy is obtained. 

\begin{figure}[h]
\hspace{-0.5cm}\includegraphics[width=0.35\textwidth,valign=c]{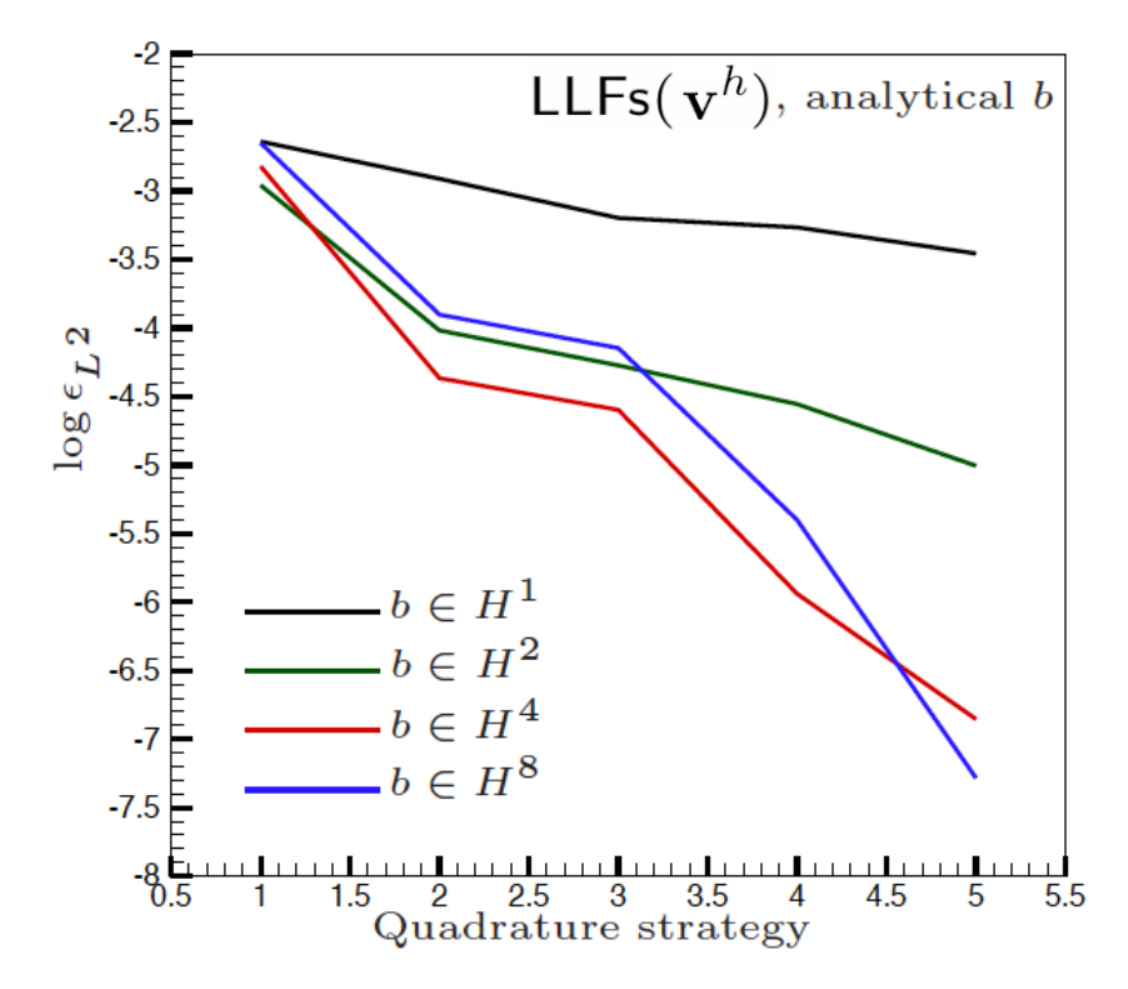}\hspace{0.1cm}\includegraphics[width=0.3\textwidth,valign=c]{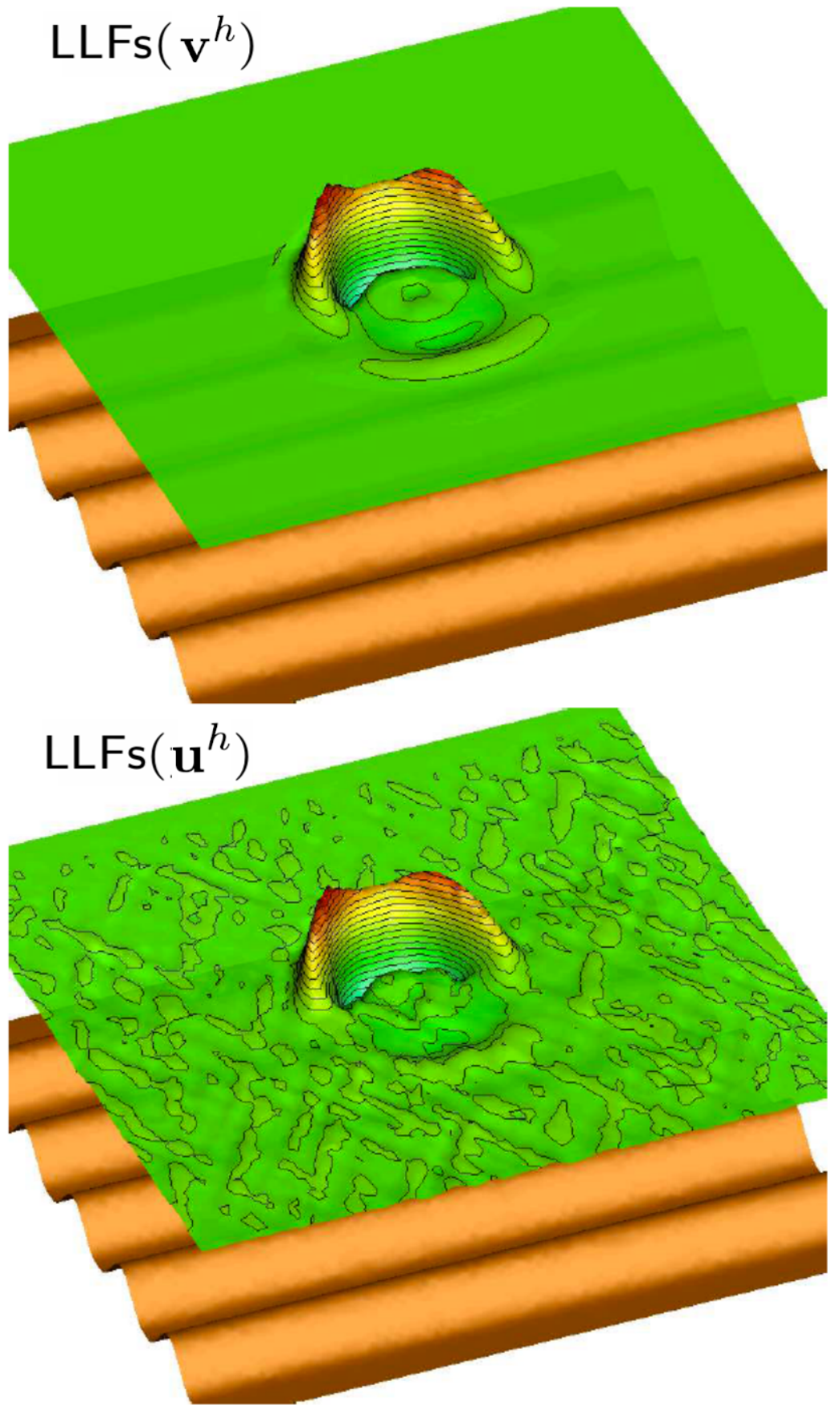}\includegraphics[width=0.38\textwidth,valign=c]{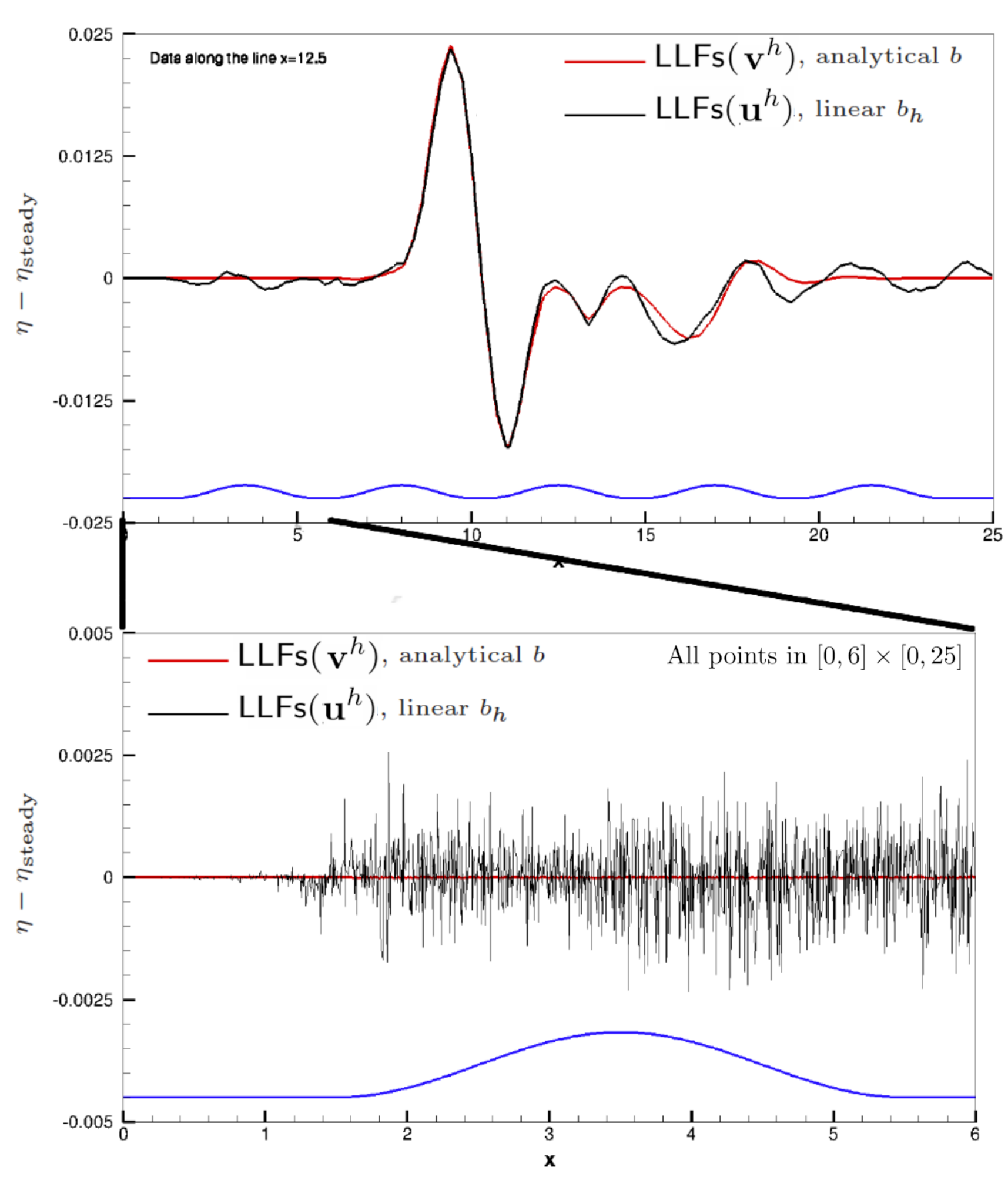}
\caption{\label{fig:nrg2d-1}{Example of a 2D perturbation in a constant energy steady flow.}}
\end{figure}

  Then, a perturbation  is added to   the initial total energy level, \mario{corresponding to an increase in  
  free surface level 
  of the order of 0.05$m$ (see  \cite{R15}  for details). }
  The evolution of the free surface perturbation $\delta\eta$
  obained by removing the exact steady solution from the results \mario{is reported in the middle  column } of  figure \ref{fig:nrg2d-1}.
   The   LLFs$(\bv^h)$  with analytical expressions for the topography provides a perfectly clean evolution of the perturbation.
   \mario{A somewhat noisy  result is obtained with the standard method, with spurious effects still relatively small compared to the actual perturbation, which 
   would be hardly obtainable with other approaches. However, this shows how sensitive
   the results may be to the mesh, and that the multidimensional case may still require some improvements.} 


  \begin{figure}[h]
   \centering
\hspace{-0.5cm}\includegraphics[scale=0.175,angle=0]{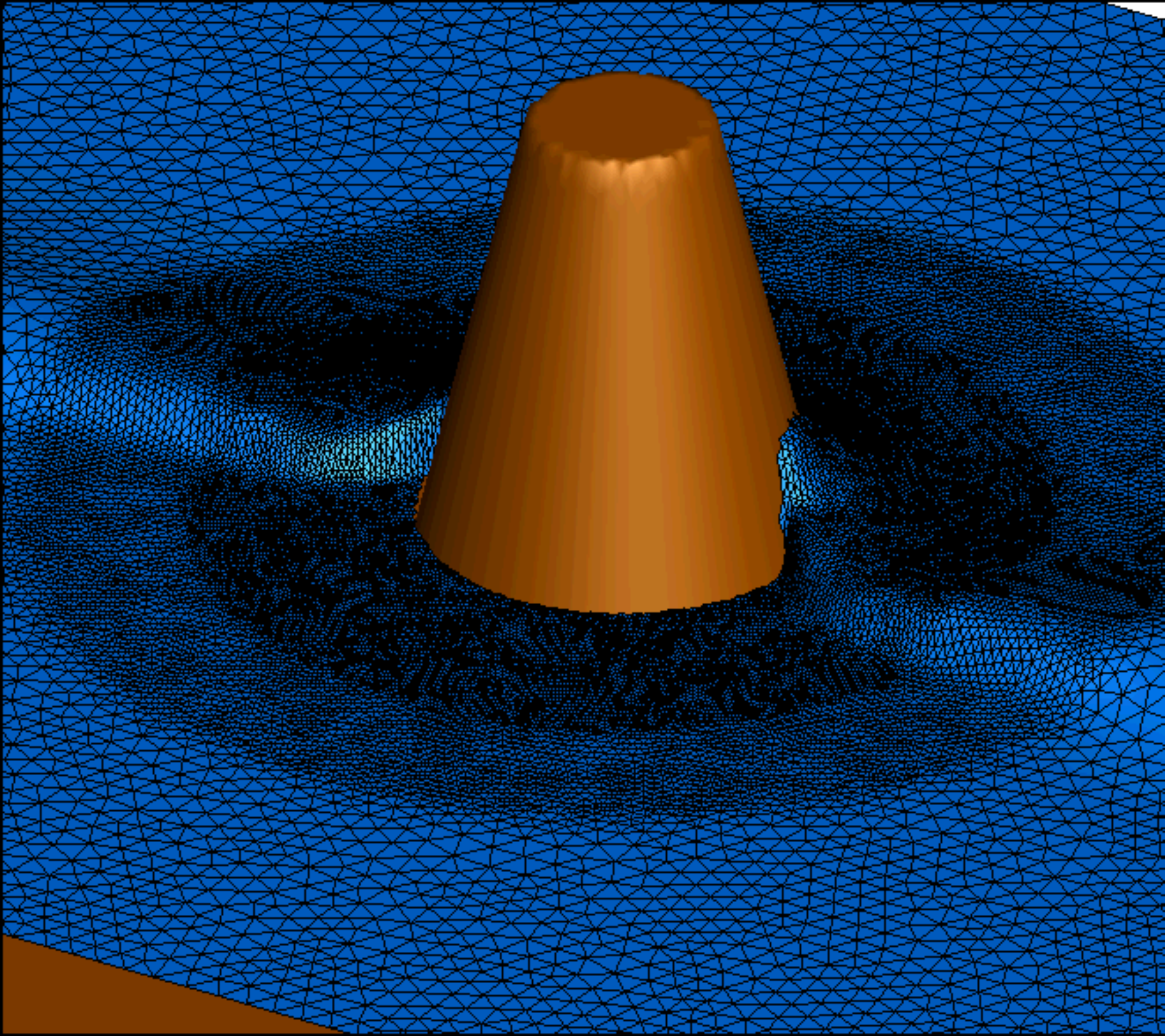}\hspace{-0.1cm}\includegraphics[scale=0.175,angle=0]{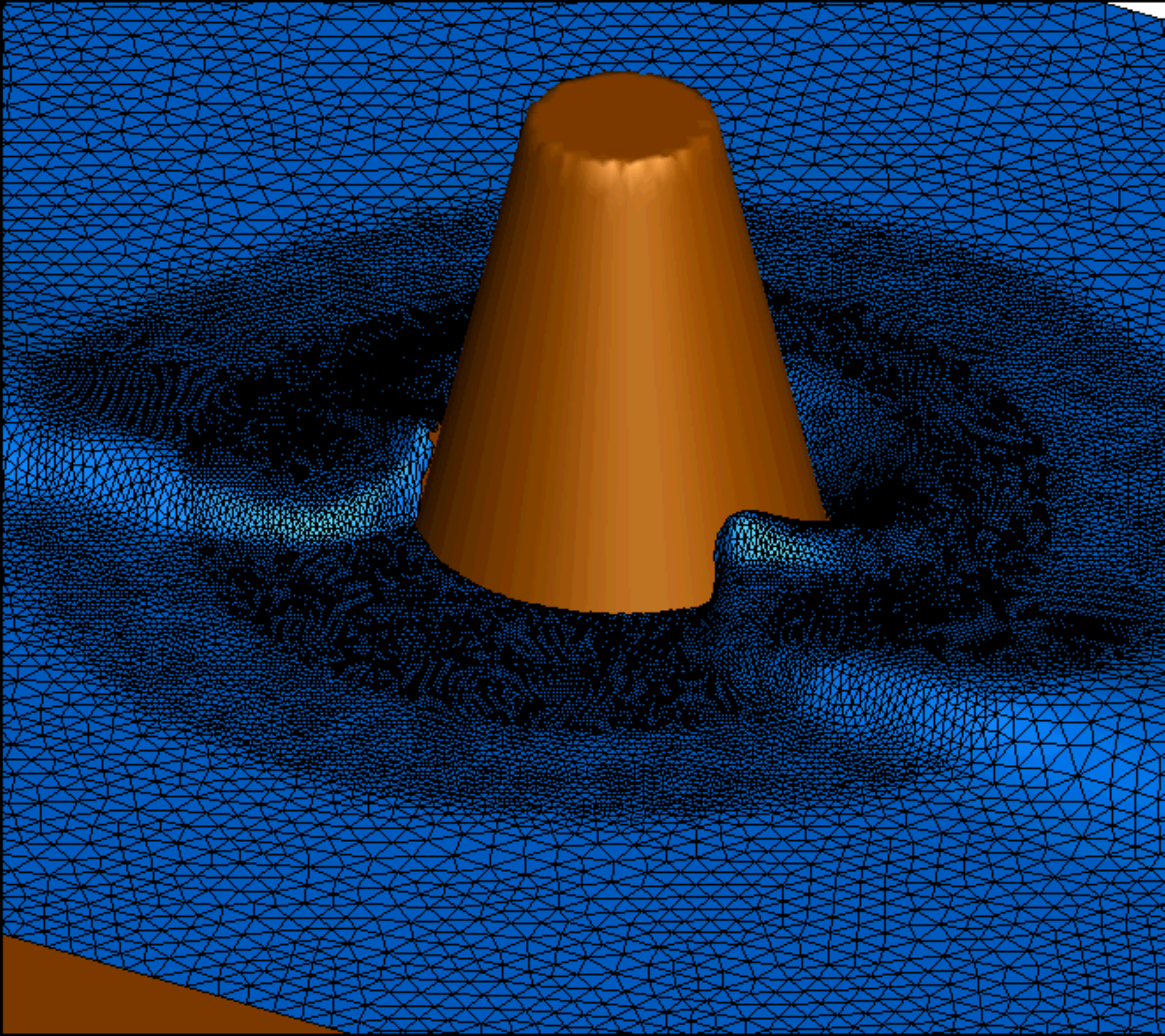}\hspace{-0.1cm}\includegraphics[scale=0.175,angle=0]{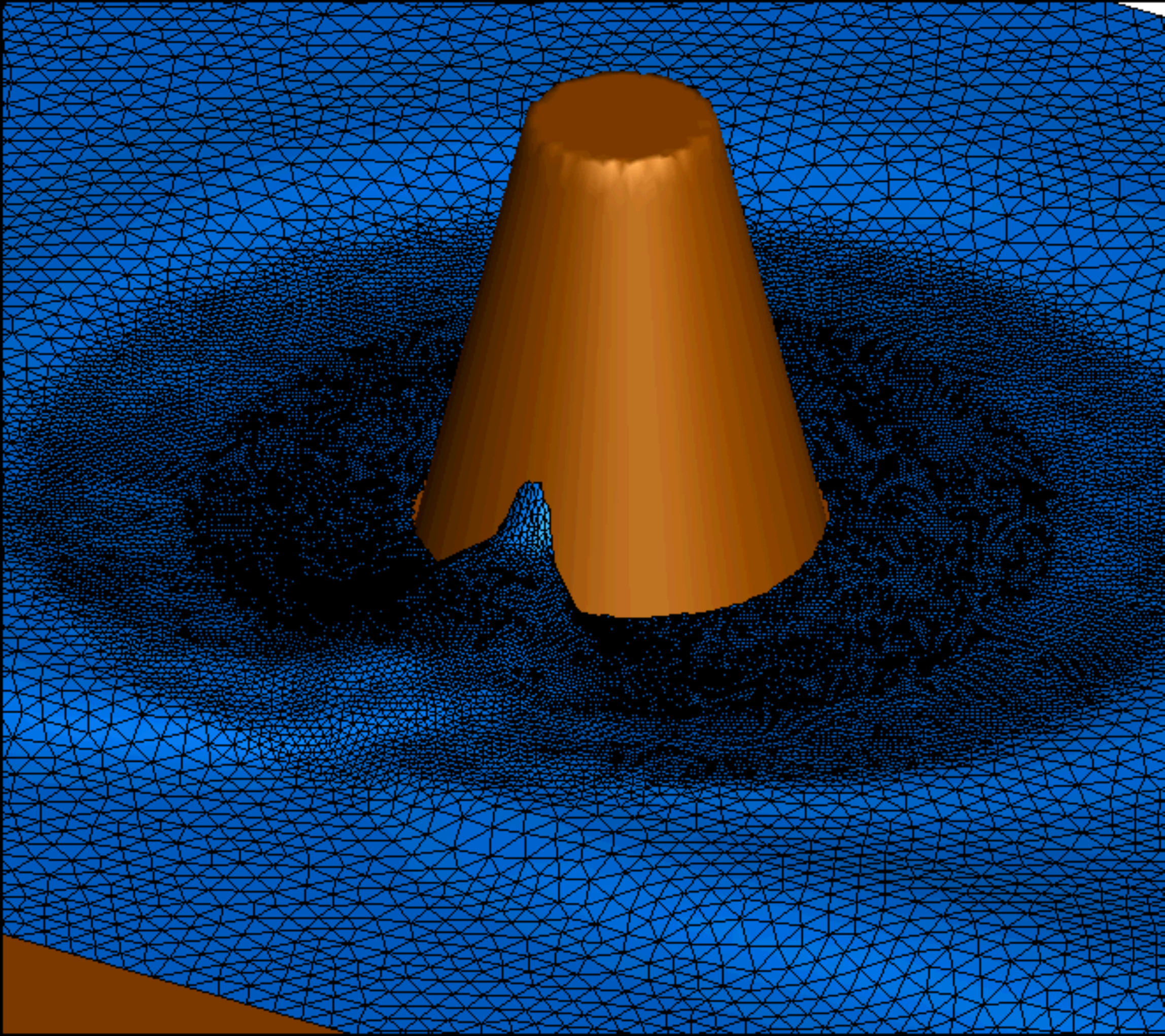}

 \caption{\small  Conical Island. 3D view of the wave runup process, and gauge locations.\label{fig:briggs_RD_o}}
 \end{figure}

\subsection{Shallow water with dry areas}

This is a  very classical benchmark    reproducing   the experiments of \cite{Br:95} on the runup of solitary waves on a conical island. We refer to the above reference, and to \cite{ArR:18}  for details. 
 We report here  results  on moving adaptive meshes \mario{based on the 
 use of the Lax-Friedrichs' based  distribution which}
allows to have control on the non-negativity of the water depth (cf. \cite{R15,ArR:18}). Figure \ref{fig:briggs_RD_o} shows a reference solution obtained on a relatively fine mesh
(uniformly refined in the interaction region). 

  \begin{figure}[h]
 \centering
 \includegraphics[scale=0.2,angle=-90]{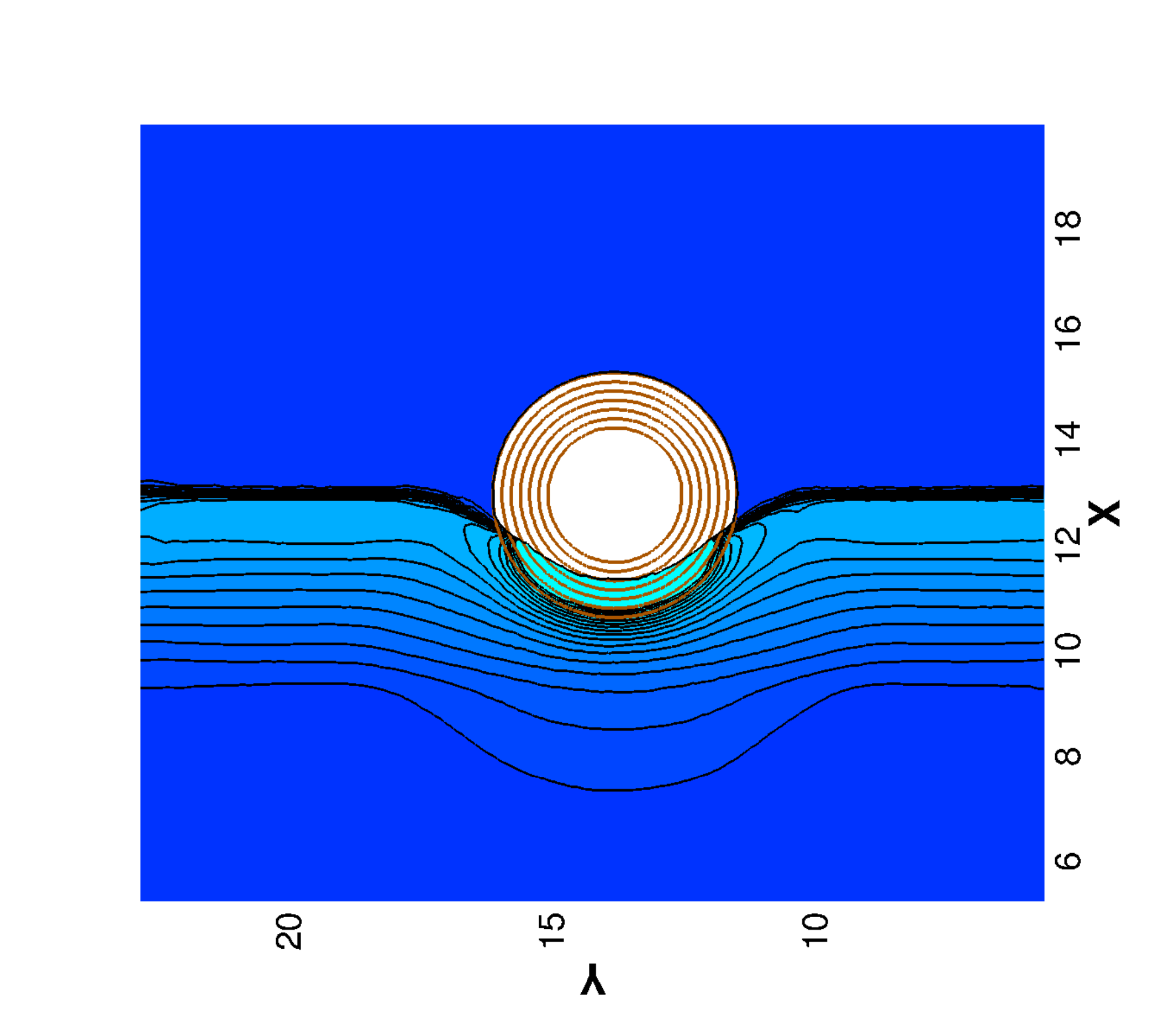}\hspace{-0.1cm}\includegraphics[scale=0.2,angle=-90]{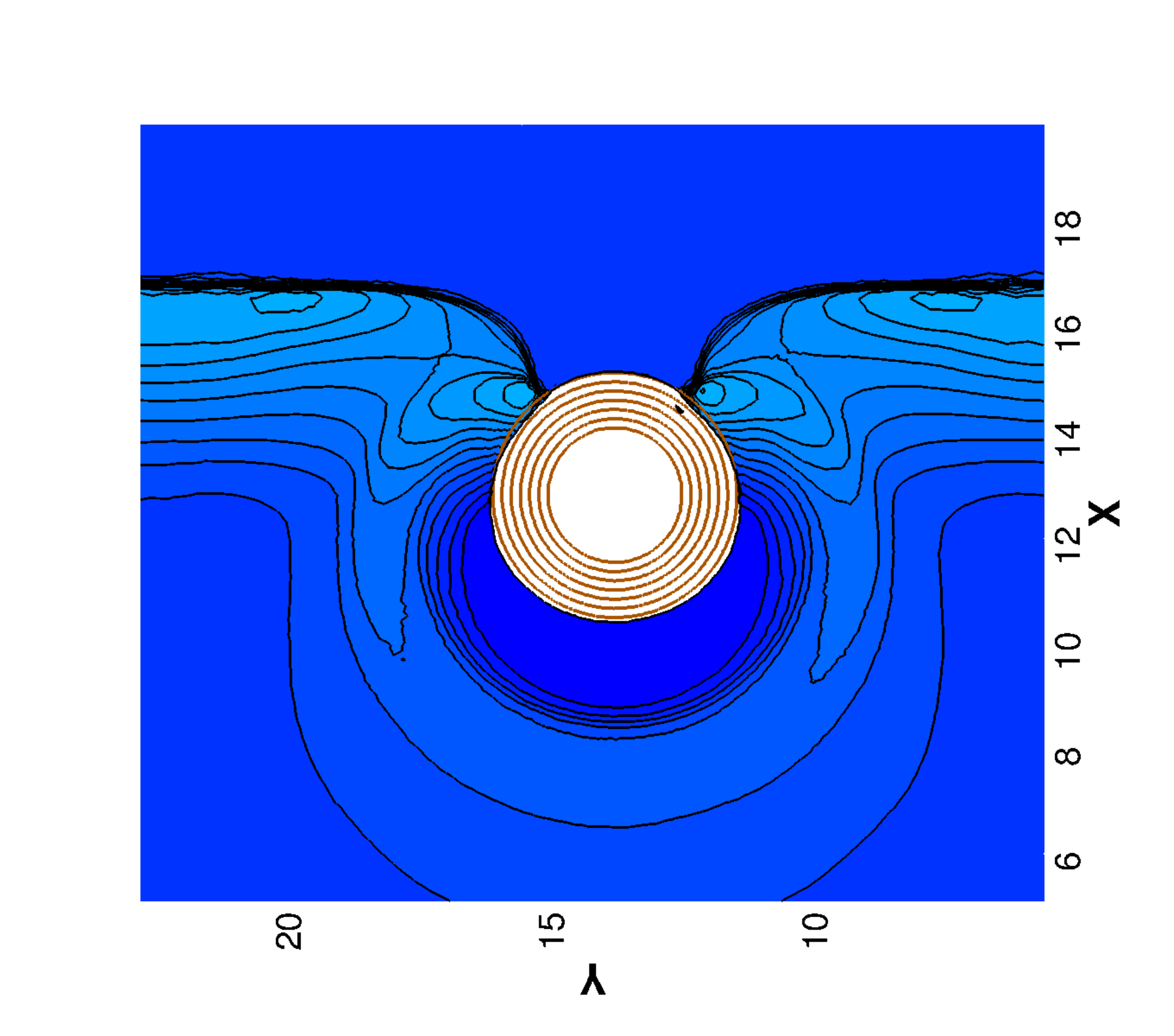}\hspace{-0.1cm}\includegraphics[scale=0.2,angle=-90]{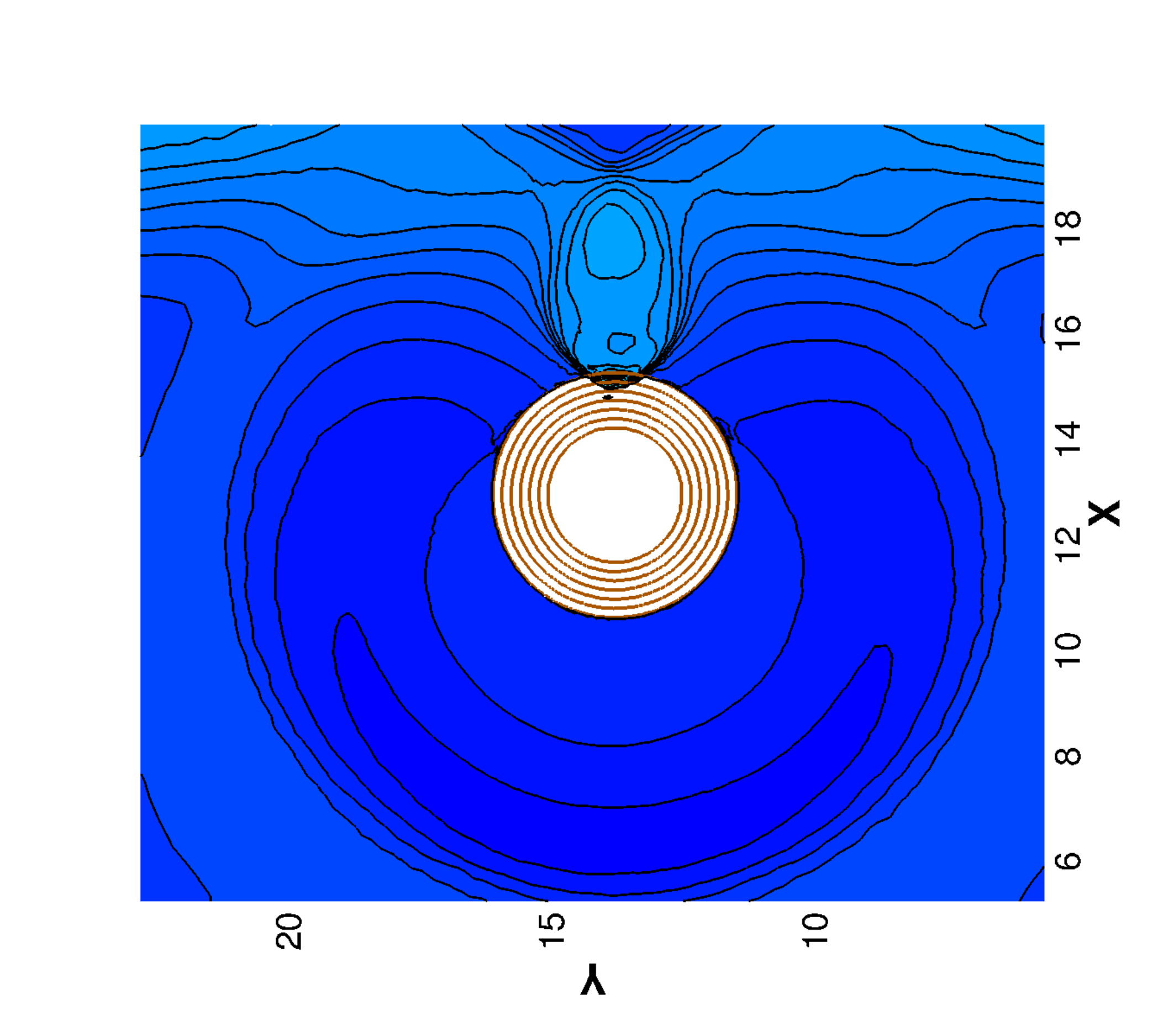}

 \includegraphics[scale=0.2,angle=-90]{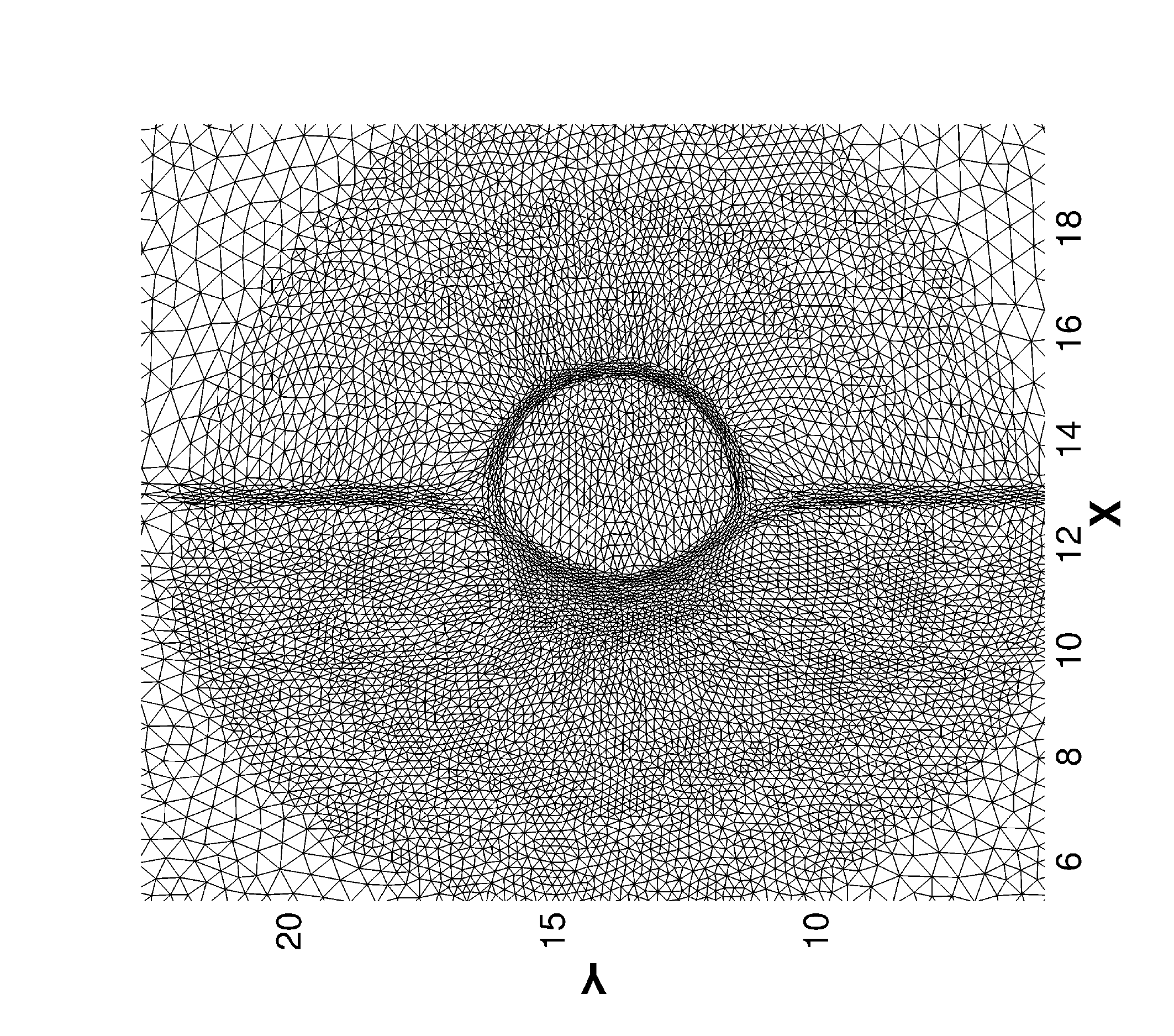}\hspace{-0.1cm}\includegraphics[scale=0.2,angle=-90]{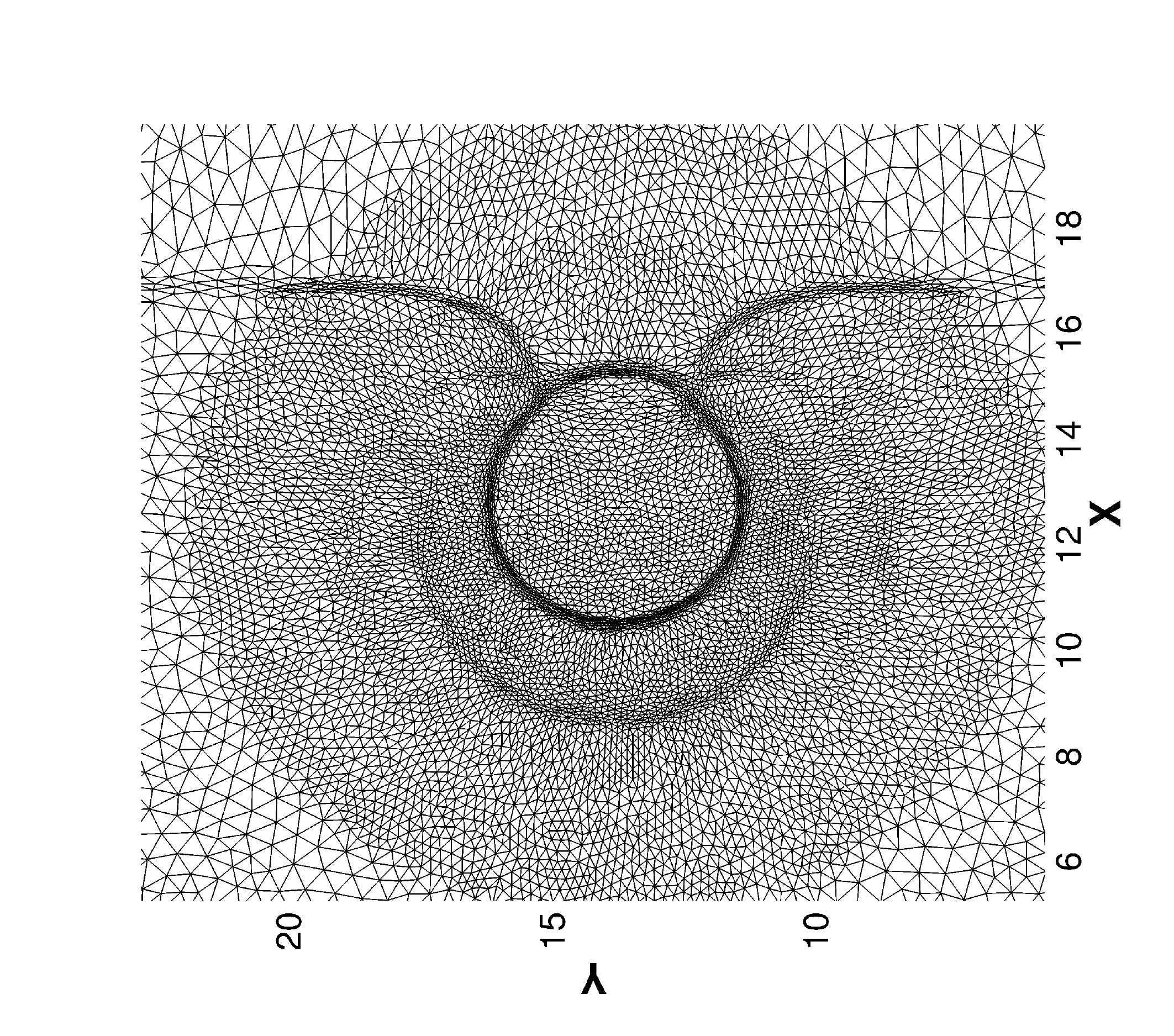}\hspace{-0.1cm}\includegraphics[scale=0.2,angle=-90]{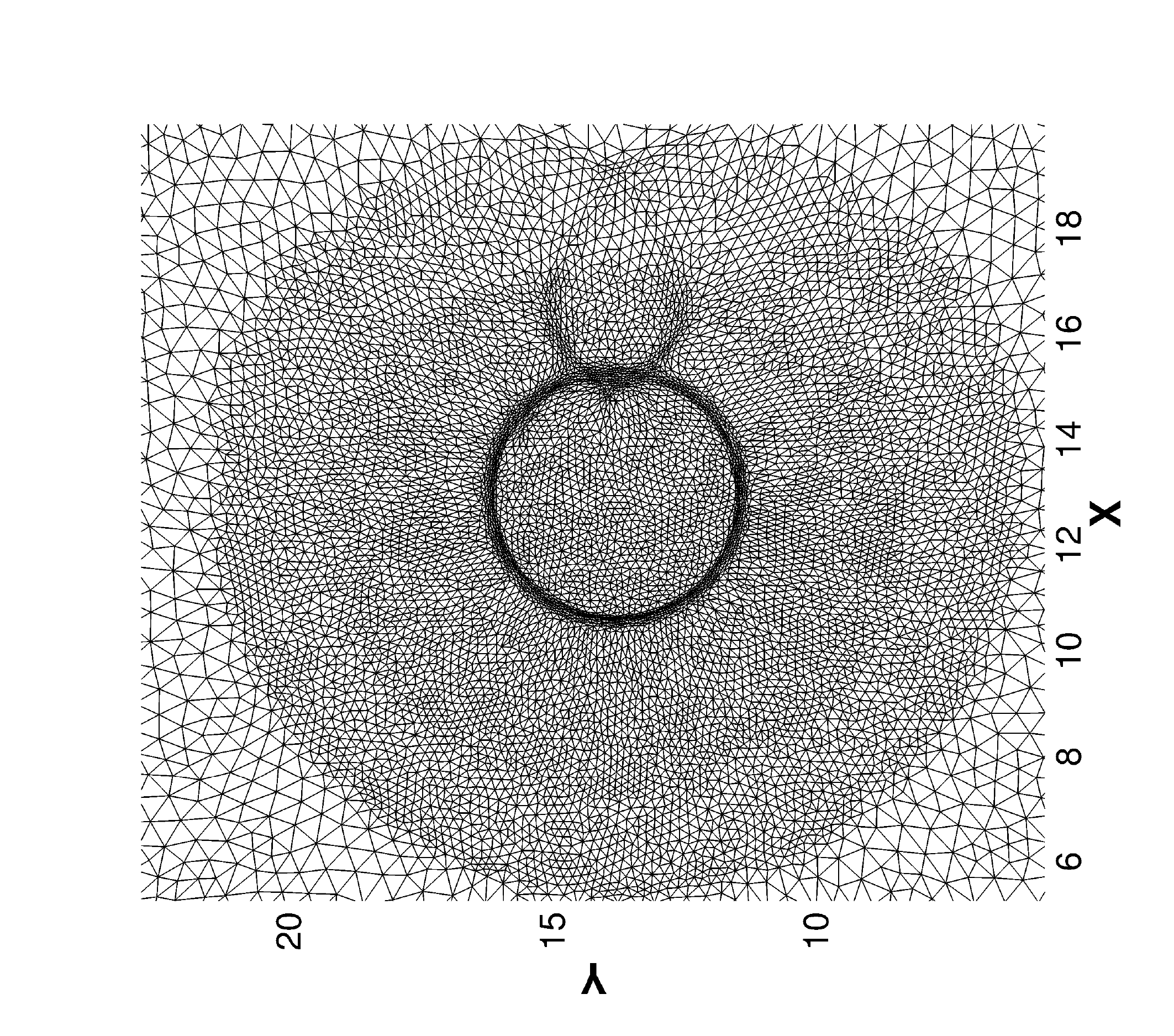} 

 \caption{ \label{fig:briggs_RD_mesh}\small  Conical Island;   free surface  contours  at  $t=6.0,8.0,10.0 \, [s]$, and   adapted meshes  }
 \end{figure}

Solution contours and adapted meshes 
 obtained  
 with the ALE adaptive method discussed in \cite{ArR:18} are shown  on figure \ref{fig:briggs_RD_mesh}. 
\mario{To  check the improvement }
brought by adaptation \mario{we report on figure \ref{fig:briggs1-GAUGE} the
 time series of the water elevation at the front and rear of the island (top row), and the CPU times. 
The coarse mesh  (dashed blue curves) clearly fails to provide correct values of the maximum runup levels, which are much better 
predicted with the adaptive ALE approach (dashed red curves), with CPU  times  of less than 30\% those of the fine mesh.}

\mario{The bottom-right plot on the same figure proves the effectiveness of the
%
%
  mass conservative bathymetry remap of  \cite{ArR:18} also in presence of dry areas.}

 \begin{figure}  
 \hspace{-1.5cm}
 \includegraphics[scale=0.85]{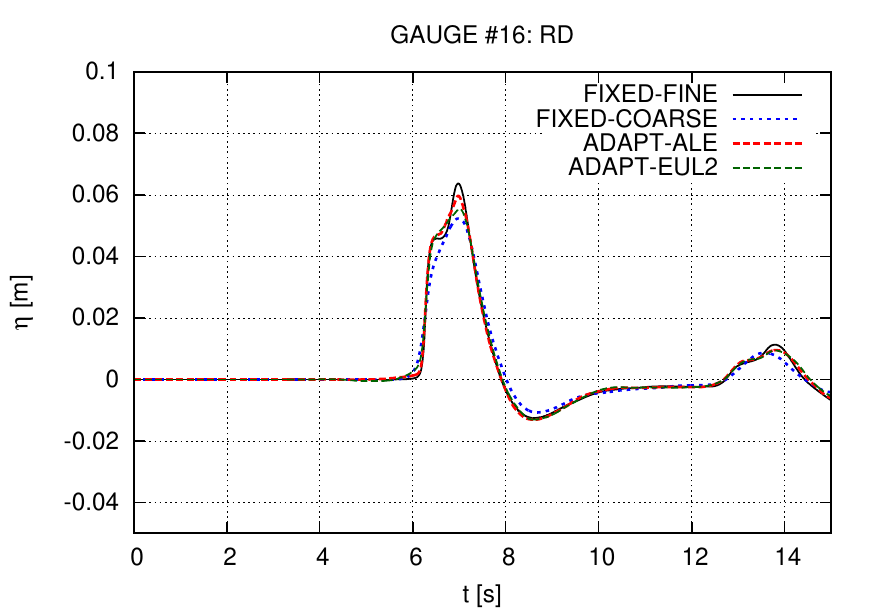}\hspace{-0.2cm}
\includegraphics[scale=0.85]{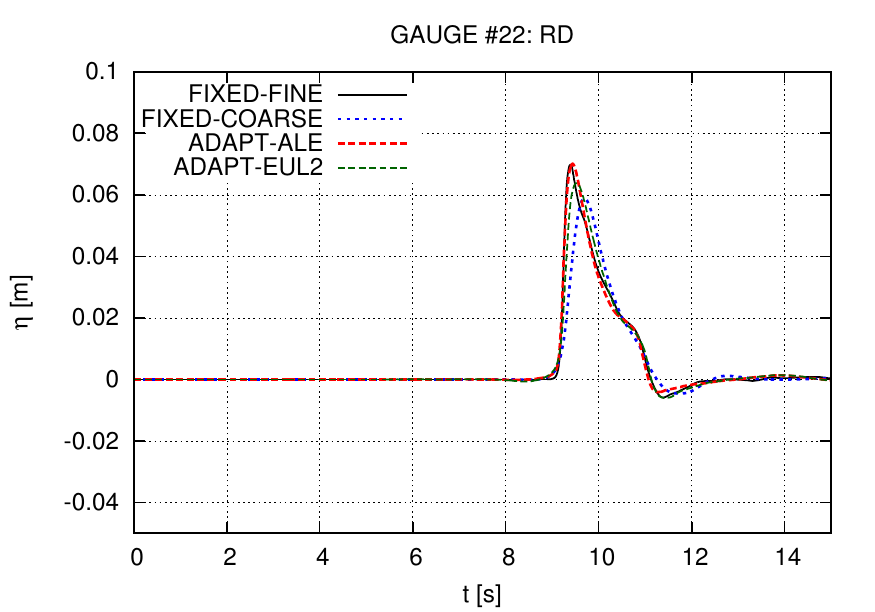}
 \begin{minipage}{0.4\textwidth}
 \begin{tabular}{|c|c|c|c|c|}
\hline 
 Mesh      & Nodes & CPU [s] (\%MMPDE) \tabularnewline
\hline 
\hline 
Fix  & 10401  & 171.30    \tabularnewline
\hline 
Fix    & 37982  & 1785.96   \tabularnewline
\hline 
Adapt+ALE   & 10401  & 510.65 (38.8\%)  \tabularnewline
 \hline
\end{tabular}
 \end{minipage} \hfill
\begin{minipage}{0.6\textwidth}
 \hspace{0.75cm}
  \includegraphics[scale=0.9]{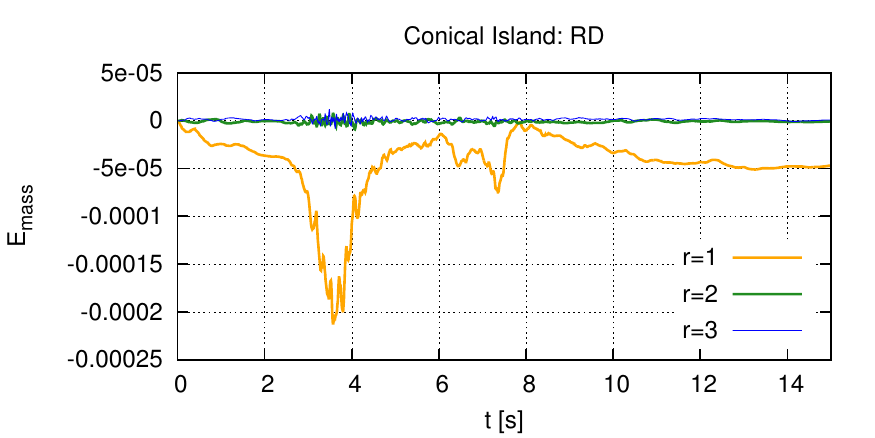}
 \end{minipage}
 \caption{ \label{fig:briggs1-GAUGE}\small Conical Island. Top: free surface time series at the front and rear runup gauges.
 Bottom: mass conservation (left) and CPU times (right). 
 }
 \end{figure}
%
%
 

%



\section{Conclusion and outlook}
In this paper, we have described the numerical framework known as   Residual Distribution.  This formalism encompass most, if not all, the known schemes (at least known by the authors), both on  structured and unstructured meshes. In particular, we have shown that despite not  being formulated in terms of numerical fluxes, one can easily provide explicit definitions of local numerical fluxes, with a consistency defined in the classical sense.
This result applies to continuous finite elements functional approximation, for example, showing that they are locally conservative. Examples of  explicit formula for the flux have been provided in the multi-dimensional case. 
For non-homogenous problems,  these methods are by their nature well-balanced. We have in particular shown that the classical weighted residual formulation  boils down to a global flux method in one space dimension,
allowing to embed a more general notion of consistency better suited for balance laws.  These properties naturally generalize to arbitrary order of accuracy via appropriate  choices for the underlying polynomial approximation and
quadrature formulas. Other extensions of the method have not been discussed here, as e.g. the consistent treatment of higher order derivatives 
 \cite{abgralldeSantisSISC,AbgralldeSantisNS,MAZAHERI2015455,rf14,MAZAHERI2016593} and  the treatment of non-conservative models \cite{ABGRALL201810,abgrall2020simulation}.

 \remi{As the reader might have seen, there are some similarities between the RD schemes and the wave propagation algorithm by R. LeVeque, see e.g. \cite{LeVequeWave}, where the emphasis is put on the upwind nature of the algorithm as well as a non linear stabilisation on the flux. There is also similarities with  the methods developed by A. Lerat and his co-authors, see e.g. \cite{MR3534868,MR2787945}.  They are residual methods, formally the look like the SUPG scheme with the difference that they introduce a non linear stabilisation to deal with the flow discontinuities.  There are more differences with the recently developed active flux method however, see \cite{MR4002758,MR3434930,MR4176953}: the approximation of the solution is completely different, and the evolution operator uses in depth the structure of the exact one. These schemes are probably in their infancy and further development will certain occur, see \cite{MR3995974,MR4176953,MR4002758,abgrall2020combination}.}
 
The   topics discussed are the basis for many future possible developments aiming in particular at further generalizing some of the properties discussed,
as e.g. the notion of well balanced in multiple dimensions, as well as further combining them with adaptive strategies. 
%
%

\section*{Acknowledgements}
We  acknowledge the work of many PhD  and post-docs students: 
M. Mezine, C. Tav\'e, A. Larat, G. Baurin, D. de Santis, L. Arpaia, A. Filippini, P. Bacigaluppi, D. Torlo, and S. Tokareva.
Discussions with S. Karni and P. Roe (both from the University of Michigan) are also warmly acknowledged.

RA has been partially financed by SNF grant 200020\_175784 and an International Chair of Inria.

\bibliographystyle{plain}
\bibliography{papier}
\end{document}